%% file: perfG2.tex
\documentstyle[epic,eepic]{article}
\input{preamble}

\input{kigou} 
\input{picture}
\begin{document}
\begin{center}
{\LARGE PERFECT CRYSTALS OF $U_q(G^{(1)}_2)$} \vspace{3mm} \\
{\large  Shigenori Yamane}\vspace{3mm}\\
Department of Mathematical Science,Faculty of Engineering Science,\\
Osaka University, Toyonaka, Osaka 560, Japan \vspace{3mm}\\
yamane@sigmath.es.osaka-u.ac.jp
\end{center}

\section{Introduction}
\input{intro}
\section{Results}
\input{result}

\section{Preliminaries}
\label{section:prelim}
\input{crystal}

\input{s_perf}
\section{Level-one representation and fusion construction }
\label{section:existance}
\input{fusion}
\section{Perfectness of the crystals}
\label{section:perfectness}
\input{prelim}

\input{const5}
\input{commuta}

\input{minimal}

\newpage
\input{appendix}
\clearpage

\input{biblio}
\end{document}

%% file: preamble.tex
\makeatletter
\@addtoreset{equation}{subsection}
\makeatother

\newtheorem{theorem}{Theorem}[section]
\newtheorem{proposition}[theorem]{Proposition}
\newtheorem{definition}[theorem]{Definition}
\newtheorem{lemma}[theorem]{Lemma}
\newtheorem{corollary}[theorem]{Corollary}
\newtheorem{remark}[theorem]{Remark}
\newtheorem{example}[theorem]{Example}

\newfont{\frak}{eufm10}
\font\frak=eufm10\relax
\newfont{\frakbigger}{eufm10 scaled \magstep2}
\newfont{\fraks}{eufm5}
\font\fraks=eufm5\relax

\newcommand{\vminiD}[1]{\rule{0mm}{4mm}%
\begin{minipage}{4.1mm}
\hrule width 4mm\hbox{\vrule height 3.2mm depth 0.8mm%
\hbox to 3.75mm{\small\hss$#1$\hss}\vrule}\hrule width 4mm
\end{minipage}} 
\newcommand{\vminimini}[1]{\rule{0mm}{3mm}%
\begin{minipage}{3.1mm}
\hrule width 3mm\hbox{\vrule height 2.4mm depth 0.6mm%
\hbox to 2.75mm{\scriptsize\hss$\mojihenkan{#1}$\hss}\vrule}\hrule width 3mm
\end{minipage}} 
\newcommand{\vyokoD}[2]{\rule{0mm}{4mm}%
\begin{minipage}{8.4mm}
\hrule width 8.37mm\hbox{\vrule height 3.2mm depth 0.8mm%
\hbox to 4mm{\small\hss$#1$\hss}\vrule%
\hbox to 4mm{\small\hss$#2$\hss}\vrule}\hrule width 8.37mm
\end{minipage}} 
\newcommand{\vminilD}[1]{\rule{0mm}{4mm}%
\begin{minipage}{6.1mm}
\hrule width 6mm\hbox{\vrule height 3.2mm depth 0.8mm%
\hbox to 5.75mm{\small\hss$#1$\hss}\vrule}\hrule width 6mm
\end{minipage}} 
\newcommand{\vmini}[1]{\vminiD{\mojihenkan{#1}}}
\newcount\bboxnumber
\def\mojihenkan#1{
\ifnum #1<0 \bboxnumber=-#1 \overline{\mojihenkan{\the\bboxnumber}}
\else\ifnum #1<7 #1
\else \bboxnumber=#1 \advance \bboxnumber by -7 
\mojihenkanshift{\the\bboxnumber} \fi\fi}
\def\mojihenkanshift#1{\ifcase#1 0_1\or 0_2\or \cdot \fi}
\newcommand{\vminiz}[1]{%
\begin{minipage}{4.1mm}
\hrule width 4mm\hbox{\vrule height 3.2mm depth 0.8mm%
\hbox to 3.75mm{\small\hss$\mojihenkan{#1}$\hss}\vrule}\hrule width 4mm
\end{minipage}}

\newcommand{\vtate}[2]{\begin{minipage}{4.1mm}
\hrule width 4mm\hbox{\vrule height 3.2mm depth 0.6mm%
\hbox to 3.75mm{\small\hss$#1$\hss}\vrule}\hrule width 4mm
\hbox{\vrule height 3.2mm depth 0.6mm%
\hbox to 3.75mm{\small\hss$#2$\hss}\vrule}\hrule width 4mm
\end{minipage}} 

\newcommand{\sayouf}[4]{$f_#1\:\vmini{#2}\:=\:
\if#41\else#4\fi\:\vmini{#3}\:$}
\newcommand{\sayoue}[4]{$e_#1\:\vmini{#2}\:=\:
\if#41\else#4\fi\:\vmini{#3}\:$}
\newcommand{\sayout}[3]{$t_#1\:\vmini{#2}\:=\:
\if#30\else q_{#1}^{#3}\fi\:\vmini{#2}\:$}
\newcommand{\sayounokori}[3]{$\displaystyle \:+ 
\if#21 \if#11 \: else \:#1\: \fi \else \:\frac{#1}{#2}\: \fi
\vmini{#3}\:$}
\newcommand{\sayouft}[4]{$\tf{#1}\:\vmini{#2}\:=\:
\if#41\else#4\fi\:\vmini{#3}\:$}
\newcommand{\sayouet}[4]{$\te{#1}\:\vmini{#2}\:=\:
\if#41\else#4\fi\:\vmini{#3}\:$}

%% file: kigou.tex
\def\U#1{U_q ( #1 )}
\def\Up#1{U_q^{\prime}( #1 )}
\def\G{{\mbox{\frak g}}}
\def\h{{\mbox{\frak h}}}
\def\fW{{\mbox{\frak W}}}
\def\fWs{{\mbox{\fraks W}}}
\def\fN{{\mbox{\frak N}}}
\def\fNs{{\mbox{\fraks N}}}
\def\gone{G_2^{(1)}} 
\def\te#1{\tilde{e}_{#1}}
\def\tf#1{\tilde{f}_{#1}}
\def\teg#1{\tilde{e}^{\CG}_{#1}}
\def\tfg#1{\tilde{f}^{\CG}_{#1}}
\def\tea#1{\tilde{e}^{\CA}_{#1}}
\def\tfa#1{\tilde{f}^{\CA}_{#1}}
\def\ea#1#2{\left ( \tea{#1} \right )^{#2}}
\def\fa#1#2{\left ( \tfa{#1} \right )^{#2}}
\def\eg#1#2{\left ( \teg{#1} \right )^{#2}}
\def\fg#1#2{\left ( \tfg{#1} \right )^{#2}}
\def\ve#1{\varepsilon_{#1}} 
\def\vp#1{\varphi_{#1}} 
\def\sayouA#1#2{\fa{1}{#1}\fa{0}{#2}}

\def\hwv#1#2#3#4{\overline{b}^{\,#1,#2}_{(#3,#4)}}

\def\hwvl#1#2#3{\overline{b}^{\,l,#1}_{(#2,#3)}}

\def\hwvs{\overline{b}^{\,l,i}_{(l-i,j)}}
\def\lwv#1#2#3{\underline{b}^{\,#1,#2}_{#3}}
\def\syme#1{\mbox{\rm C}_{\CA}\left ( #1 \right )}
\def\symg#1{\mbox{\rm C}_{\CG}\left ( #1 \right )}

\def\Z{{\bf Z}} 

\def\Q{{\bf Q}}
\def\CA{{\cal A}}
\def\CC{{\cal C}}
\def\CG{{\cal G}}
\def\CW{{\cal W}}
\def\CWb{\overline{\cal W}}
\def\mod{{\rm mod}}
\def\wt{{\rm wt}}

\def\proof{\noindent{\bf Proof \hspace{2mm}}}

\def\bB#1#2{%
\begin{minipage}{1mm}
\hrule width 1mm\hbox{\vrule height 4mm depth 0mm\hbox to 0.25mm{}\vrule}
\hrule width 1mm
\end{minipage}
\raisebox{-0.4mm}{$#1$}
\begin{minipage}{1mm}
\hrule width 1mm\hbox{\hbox to 0.5mm{}\vrule height 4mm depth 0mm
\hbox to 0.25mm{}\vrule}\hrule width 1mm
\end{minipage}
^{\;#2\;}\rule{0mm}{6mm}
}
\def\BB#1{%
\begin{minipage}{1mm}
\hrule width 1mm\hbox{\vrule height 3.5mm depth 1.5mm\hbox to 0.25mm{}\vrule}
\hrule width 1mm
\end{minipage}
\mbox{ $#1$}
\begin{minipage}{1mm}
\hrule width 1mm\hbox{\hbox to 0.5mm{}\vrule height 3.5mm depth 1.5mm
\hbox to 0.25mm{}\vrule}\hrule width 1mm
\end{minipage}\;
}
\def\bb#1#2{\bB{\makebox[3mm]{\bf$\mojihenkan{#1}$}}{#2}}

\def\uu#1{\Psi_{#1}}
\def\saba#1{{\rm (a#1)}} 

%% file: picture.tex
\def\younglambda{
\begin{minipage}{45mm}
\hrule width 47.9mm\hbox{\vrule height 4mm depth 1.5mm%
\hbox to 9mm{\hss\small$b^{m\!+\!n}_1$\hss}\vrule
\hbox to 10mm{\hss$\cdots$\hss}\vrule
\hbox to 8mm{\hss\small$b^{m\!+\!1}_1$\hss}\vrule
\hbox to 5mm{\hss\small$b^m_1$\hss}\vrule
\hbox to 10mm{\hss$\cdots$\hss}\vrule
\hbox to 5mm{\hss\small$b^1_1$\hss}\vrule}
\hrule width 47.9mm\hbox{\vrule height 4mm depth 1.5mm%
\hbox to 9mm{\hss\small$b^{m\!+\!n}_2$\hss}\vrule
\hbox to 10mm{\hss$\cdots$\hss}\vrule
\hbox to 8mm{\hss\small$b^{m\!+\!1}_2$\hss}\vrule}
\hrule width 27.5mm
\end{minipage}
\rule{0mm}{8mm} 
}
\def\youngHigh{
\begin{minipage}{41mm}
\hrule width 40.9mm\hbox{\vrule height 4mm depth 1.5mm%
\hbox to 5mm{\hss\small$1$\hss}\vrule
\hbox to 10mm{\hss$\cdots$\hss}\vrule
\hbox to 5mm{\hss\small$1$\hss}\vrule
\hbox to 5mm{\hss\small$1$\hss}\vrule
\hbox to 10mm{\hss$\cdots$\hss}\vrule
\hbox to 5mm{\hss\small$1$\hss}\vrule}
\hrule width 40.9mm\hbox{\vrule height 4mm depth 1.5mm%
\hbox to 5mm{\hss\small$2$\hss}\vrule
\hbox to 10mm{\hss$\cdots$\hss}\vrule
\hbox to 5mm{\hss\small$2$\hss}\vrule}
\hrule width 20.5mm
\end{minipage}
\rule{0mm}{8mm} 
}
\newcount\graphloc
\newcount\edgehight
\makeatletter
\newcommand{\graphstraight}[3]{%
\graphloc 0
\@for\vertexnumber:=#1\do{%
\put(\graphloc,0){\makebox(0,3){\vmini{\vertexnumber}}}%
\advance \graphloc by 4\put(\graphloc,0){\vector(1,0){12}}%
\advance \graphloc by 16}
\put(\graphloc,0){\makebox(0,3){\vmini{#2}}}%
\graphloc 10
\@for\vertexnumber:=#3\do{%
\put(\graphloc,4){\makebox(0,0){$\vertexnumber$}}%
\advance \graphloc by 20}}%
\makeatother
\newcommand{\crystalgraphGtwo}[1]{%
\put(5,30){\graphstraight{1,2,3}{4}{1,2,2}}
\put(69,30){\vector(#1,1){12}}\put(69,30){\vector(#1,-1){12}}
\if#11\edgehight 42\else \edgehight 34\fi
\put(85,\the\edgehight){\graphstraight{5,7}{-5}{2,2}}
\put(129,\the\edgehight){\vector(#1,-1){12}}
\if#11\edgehight 18\else \edgehight 26\fi
\put(85,\the\edgehight){\graphstraight{6,8}{-6}{1,1}}
\put(129,\the\edgehight){\vector(#1,1){12}}
\put(145,30){\graphstraight{-4,-3,-2}{-1}{2,2,1}}}

\newcommand{\hakomoji}[1]{\hbox to 4mm{\small\hss$#1$\hss}\vrule}
\newcommand{\hakomojilong}{\hbox to 6mm{\small\hss$\cdots$\hss}\vrule}
\newcommand{\crystalbase}[2]{%
\rule[-2.5mm]{0mm}{7mm}
\begin{minipage}{14mm}
\hrule width 14.45mm \hbox{\vrule height 3.2mm depth 0.8mm%
\hakomoji{#1}\hakomojilong\hakomoji{#2}}
\hrule width 14.45mm
\end{minipage}}

\def\crystalGtwoMoto{
\put(22,50){
\put(0,40){\makebox(0,0){\vminiz{1}}}
\put(-2,38.66){\line(-3,-2){8}}
\put(-12,32){\makebox(0,0){\vminiz{2}}}
\put(-10,30){\line(1,-1){4}}
\put(-4,24){\makebox(0,0){\vminiz{3}}}
\put(-2,22){\line(1,-1){4}}
\put(4,16){\makebox(0,0){\vminiz{4}}}
\put(6,14){\line(1,-1){4}}
\put(2,14.66){\line(-3,-2){8}}
\put(12,8){\makebox(0,0){\vminiz{6}}}
\put(-8,8){\makebox(0,0){\vminiz{5}}}
\put(10,6.4){\line(-5,-4){5.5}}
\put(-6.5,6){\line(3,-4){2.3}}
\put(2,0){\makebox(1,2){\vminiz{8}}}
\put(-2,0){\makebox(-1,2){\vminiz{7}}}
\put(6,-6.4){\line(-5,4){7}}
\put(-10,-6.86){\line(7,4){10.5}}
\put(8,-8){\makebox(0,0){\vminiz{-5}}}
\put(-12,-8){\makebox(0,0){\vminiz{-6}}}
\put(6,-9.33){\line(-3,-2){8}}
\put(-10,-10){\line(1,-1){4}}
\put(-4,-16){\makebox(0,0){\vminiz{-4}}}
\put(-2,-18){\line(1,-1){4}}
\put(4,-24){\makebox(0,0){\vminiz{-3}}}
\put(6,-26){\line(1,-1){4}}
\put(12,-32){\makebox(0,0){\vminiz{-2}}}
\put(10,-33.33){\line(-3,-2){8}}
\put(0,-40){\makebox(0,0){\vminiz{-1}}}}
\put(62,50){\makebox(0,0){$\phi$}}
\put(0,0){\line(1,0){80}}
\put(0,0){\line(0,1){100}}
\put(80,100){\line(-1,0){80}}
\put(80,100){\line(0,-1){100}}
}

\def\crystallevelone{
\crystalGtwoMoto
\put(22,50){
\put(-12,-6){\vector(0,1){36}}
\put(-5.4,-14){\vector(0,1){36}}
\put(5.4,-22){\vector(0,1){36}}
\put(12,-30){\vector(0,1){36}}
}
\put(24,48){\arc{76}{0}{1.57}}
\put(24,52){\arc{76}{-1.57}{0}}
\put(62,48){\vector(0,1){0}}
\put(24,90){\vector(-1,0){0}}
\put(73,0){\makebox(3,7)[r]{$U_q(G_2^{(1)})$}}
\put(64,52){\makebox(0,0){$*$}}
\put(21,54){\makebox(0,0){$*$}}
\put(53,90){\makebox(0,0)[r]{$\uparrow \tilde{f}_0$}}
\put(76,90){\makebox(0,0)[r]{$\swarrow \tilde{f}_1$\hspace{3mm}$\searrow 
\tilde{f}_2$}} 
\put(76,85){\makebox(0,0)[r]{$*$ means minimal}}
}

\def\crystallevelone{
\crystalGtwoMoto
\put(22,50){
\put(-12,-6){\vector(0,1){36}}
\put(-5.4,-14){\vector(0,1){36}}
\put(5.4,-22){\vector(0,1){36}}
\put(12,-30){\vector(0,1){36}}
}
\put(24,48){\arc{76}{0}{1.57}}
\put(24,52){\arc{76}{-1.57}{0}}
\put(62,48){\vector(0,1){0}}
\put(24,90){\vector(-1,0){0}}
\put(73,0){\makebox(3,7)[r]{$U_q(G_2^{(1)})$}}
\put(64,52){\makebox(0,0){\small$*$}}
\put(21,54){\makebox(0,0){\small$*$}}
\put(53,90){\makebox(0,0)[r]{\small$\uparrow \tilde{f}_0$}}
\put(76,90){\makebox(0,0)[r]{\small$\swarrow \tilde{f}_1$\hspace{3mm}$\searrow 
\tilde{f}_2$}} 
\put(76,85){\makebox(0,0)[r]{\small$*$ is minimal}}
} 

\def\crystalAtwoA{
\put(20,90){\line(-1,-1){10}}
\put(20,90){\circle*{2}}
\put(40,70){\line(-1,-1){20}}
\put(20,50){\circle*{2}}
\put(20,50){\line(-1,-1){5}}
\put(40,70){\circle*{2}}
\put(45,65){\circle*{2}}
\put(10,11){\makebox(0,0){$\swarrow \tfa{1}$}}
\put(23,11){\makebox(0,0){$\searrow \tfa{0}$}}

}
\def\crystalAtwoB{
\put(20,90){\line(1,-1){50}}
\put(70,40){\line(-1,-1){35}}
\put(70,40){\circle*{2}}
\put(45,65){\line(-1,-1){29}}
\put(20,50){\line(1,-1){30}}
\put(25,45){\circle*{2}}
\put(37.5,32.5){\circle*{2}}
\put(50,20){\circle*{2}}
\put(32.5,27.5){\circle*{2}}
\put(20,40){\line(1,-1){25}}
\put(20,40){\circle*{2}}
\put(45,15){\circle*{2}}
\put(28,23){\line(1,1){13}}
\put(23,30){\makebox(0,0){\small$r$}}
}
\def\crystalAtwoC{
\put(20,90){\circle*{2}}
\put(20,90){\line(1,-1){27}}
\put(20,90){\line(-1,-1){10}}
\put(35,75){\circle*{2}}
\put(35,75){\line(-1,-1){30}}
\put(20,60){\circle*{2}}
\put(5,45){\circle*{2}}
\put(5,45){\line(-1,-1){5}}
\put(40,70){\line(-1,-1){40}}
\put(40,70){\circle*{2}}
\put(20,50){\circle*{2}}
\put(8,38){\circle*{2}}
\put(3,33){\circle*{2}}
\put(3,33){\line(1,-1){30}}
\put(5,45){\line(1,-1){18.5}}
\put(8,38){\line(1,-1){13.5}}
\put(23.5,26.5){\line(1,-5){1}}
\put(21.5,24.5){\line(6,-1){6}}
\put(24.5,21.5){\line(1,-1){13.5}}
\put(27.5,23.5){\line(1,-1){8.5}}
\put(15.5,20.5){\circle*{2}}
\put(12.5,17.5){\line(1,1){9.5}}
\put(20.5,25.5){\circle*{2}}
\put(20.5,15.5){\circle*{2}}
\put(17.5,12.5){\line(1,1){9.5}}
\put(25.5,20.5){\circle*{2}}
\put(30,0){\line(1,1){9.5}}
\put(33,3){\circle*{2}}
\put(38,8){\circle*{2}}
\put(36,15){\circle*{2}}
\put(34.5,13.5){\line(1,1){4.5}}
\put(8.6,32.4){\circle*{0.5}}
\put(11.7,29.3){\circle*{0.5}}
\put(14.8,26.2){\circle*{0.5}}
\put(26.1,14.9){\circle*{0.5}}
\put(29.2,11.8){\circle*{0.5}}
\put(32.3,8.7){\circle*{0.5}}
\put(7,7){\makebox(0,0){$\swarrow \tfa{1}$}}
\put(20,7){\makebox(0,0){$\searrow \tfa{0}$}}
}

\newcommand{\marumoji}[1]{%
{\ooalign{\hfill$\scriptstyle#1$\hfill\crcr$\bigcirc$}}}
\newcommand{\vmmyoko}[2]{\rule{0mm}{3mm}%
\begin{minipage}{6.1mm}
\hrule width 6mm\hbox{\vrule height 2.4mm depth 0.6mm%
\hbox to 2.77mm{\scriptsize\hss$\mojihenkan{#1}$\hss}\vrule%
\hbox to 2.77mm{\scriptsize\hss$\mojihenkan{#2}$\hss}\vrule}\hrule width 6mm
\end{minipage}} 
\def\Rone{\put(-2,-2){\line(-1,-1){5.8}}}
\def\Rtwo{%
\put(-0.5,-2){\line(-1,-1){5.8}}
\put(-10,-12){\line(-4,-3){7.8}}
}
\def\Rtwos{\put(0,0){\Rone}\put(-10,-10){\Rone}}
\def\Rthree{%
\put(-0.5,-2){\line(-1,-1){5.8}}
\put(-9,-12){\line(-4,-3){7.8}}
\put(-16,-12){\line(-4,-3){7.8}}
\put(-19,-22){\line(-4,-3){7.8}}
}
\def\Lone{\put(2,-2){\line(1,-1){5.8}}}
\def\Ltwo{%
\put(1,-2){\line(2,-3){4}}
\put(7,-12){\line(5,-3){9.8}}
}
\def\Ltwos{\put(0,0){\Lone}\put(10,-10){\Lone}}
\def\Lthree{%
\put(1,-2){\line(2,-3){3.8}}
\put(13,-12){\line(4,-3){7.8}}
\put(7,-12){\line(4,-3){7.8}}
\put(17,-22){\line(5,-3){9.8}}
}
\def\Cone#1#2{\put(0,0){\makebox(0,0.8){\vmmyoko{#1}{#2}}}}
\def\Ctwo#1#2#3{%
\put(-4.5,0){\makebox(0,0.8){\vminimini{#1}}}
\put(2.5,0){\makebox(0,0.8){\vmmyoko{#2}{#3}}}
}
\newcommand{\RRone}[4]{%
\put(0,0){\Rone}
\put(0,0){\Cone{#1}{#2}}
\put(-10,-10){\Cone{#3}{#4}}
}
\newcommand{\RRtwos}[6]{%
\put(0,0){\Rtwos}
\put(0,0){\Cone{#1}{#2}}
\put(-10,-10){\Cone{#3}{#4}}
\put(-20,-20){\Cone{#5}{#6}}
}
\newcommand{\RRtwo}[7]{%
\put(0,0){\Rtwo}
\put(0,0){\Cone{#1}{#2}}
\put(-10,-10){\Ctwo{#3}{#4}{#5}}
\put(-20,-20){\Cone{#6}{#7}}
}
\newcommand{\RRthree}[8]{%
\put(0,0){\Rthree}
\put(-10,-10){\Ctwo{#1}{#2}{#3}}
\put(-20,-20){\Ctwo{#4}{#5}{#6}}
\put(-30,-30){\Cone{#7}{#8}}
}

\newcommand{\gdynkin}{%
\put(10,5){\circle{10}}
\put(15,5){\line(1,0){10}}
\put(30,5){\circle{10}}
\put(33.65,2){\line(1,0){11.5}}
\put(35,5){\line(1,0){13}}
\put(33.65,8){\line(1,0){11.5}}
\put(48,5){\line(-1,1){5}}
\put(48,5){\line(-1,-1){5}}
\put(53,5){\circle{10}}
}

%% file: intro.tex
The theory of {\it crystal base} was introduced by Kashiwara in \cite{K1}. He 
proved the existence and uniqueness of this base for a quantized 
universal enveloping algebras $\U\G$, where $\G$ is an arbitrary Kac-Moody 
Lie algebra with symmetrizable Cartan matrix. 
The theory of crystal base provides a powerful combinatorial 
tool for studying the structure of the integrable highest weight modules
and the decomposition of tensor products of them.

Let $\G$ be an arbitrary Kac-Moody Lie algebra with a symmetrizable Cartan 
matrix.
Let $I$ be its index set. Denote the set of simple roots of $\G$ by  
$\{ \alpha_i \mid i\in I\}$ .  
Let $P$ be the weight lattice of $\G$.
Take an integrable $\U\G$-module $M$ and let $\tf{i}$, $\te{i}$ be 
the Kashiwara operators (cf. \S\ref{def:kashiwara}) on $M$.
Let ${\bf A}$ be the subring of ${\bf Q}(q)$ consisting of 
$f \in {\bf Q}(q)$ that is regular at $q=0$.
A crystal lattice $L$ of an integrable $\U\G$-module $M$ is free 
${\bf A}$-submodule of $M$ such that $M \cong {\bf Q}(q)\otimes_{\bf A}L$,
$L=\bigoplus_{\lambda\in P}L_\lambda$, where  
$L_\lambda=L\cap M_\lambda$ and $\te{i}L\subset L$, $\tf{i}L \subset L$.
A crystal base of the integrable $\U\G$-module $M$ is a pair $(L,B)$ such that 
(i) $L$ is a cyrstal lattice of $M$, (ii) $B$ is a ${\bf Q}$-base of $L/qL$,
(iii) $B=\sqcup_{\lambda \in P}B_\lambda$ where $B_\lambda=B\cap
(L_\lambda/qL_\lambda)$, (iv)$\te{i}B \subset B \sqcup\{0\}$, and (v) for 
$b,b' \in B$, $b'=\tf{i}b$ if and only if $b= \te{i}b'$ for $i\in I$. 
We sometimes replace condition (ii) by : $B_{ps}=B'\sqcup(-B')$ where $B'$
is a ${\bf Q}$-base of $L/qL$. In this case, w e call $(L, B_{ps})$ a 
{\it crystal pseudo base}. The quotient $B_{ps}/\{\pm 1\}$ is called 
an {\it associated crystal} of $(L,B_{ps})$.
 
In \cite{KN}, Kashiwara and Nakashima gave an explicit construction of 
crystal bases for all finite demensional 
irreducible modules over $\U\G$, where $\G$ is of type $A_n$, $B_n$, $C_n$, 
$D_n$, in terms of generalized Young tableaux. 
Kang and Misra gave a construction of crystal bases of 
$\U{G_2}$-module by a method similar to Kashiwara and Nakashima's 
in \cite{KM}. 

The notion of {\it perfect crystals} was introduced in \cite{KMN1} and 
\cite{KMN2} in order to compute the one-point functions of vertex model
in 2 dimensional lattice statistical models.
If a $\U\G$-module $M$ has perfect crystal base,     
we have the notion of paths. 
Using paths, we can compute one-point functions of vertex models. 
We define classical part of the weight lattice by $P_{cl}=\sum\Z\Lambda_i$.
A $P_{cl}$-weighted crystal is called a {\it classical crystal} and 
$P$-weighted crystal is called an {\it affine crystal}.
Let $c$ be the central element. Further we define $\ve{i}(b)$, $\vp{i}(b)$ by
$\max\{ k \geq 0 \mid \te{i}^kb \neq 0 \}$, $\max \{ k \geq 0 \mid \tf{i}^kb 
\neq 0 \}$ respectively. We define $\ve{}(b)$, $\vp{}(b)$ by $\sum_i\ve{i}(b) 
\Lambda_i$, $\sum_i\vp{i}(b) \Lambda_i$ respectively.
A classical crystal $B$ is called a {\it perfect crystal of level $l$} 
if $B$ satisfies:
(i) $B \otimes B$ is connected, (ii) there exists $\lambda_0 \in P_{cl}$ 
such that $\wt(B) \in \lambda_0 + \sum_{i \neq 0}\Z_{\leq 0}cl(\alpha_i)$ 
and that $\sharp(B_{\lambda_0}) = 1$, (iii) there is a finite dimensional 
$U_q^{\prime}(\G)$-module with a crystal pseudo base $(L,B_{ps})$ such that 
$B \cong B_{ps} / \pm 1$, (iv) for any $b \in B$, we have $\langle c, 
\ve{}(b) \rangle \geq l$, and (v) the maps $\ve{}$, $\vp{}\::\:B_l=\{b \in B 
\:|\: \langle c, \ve{}(b) \rangle = l \} \rightarrow (P^+_{cl})_l=
\{\lambda\in\sum\Z_{\geq 0}\Lambda_i \:|\:\langle c,\lambda\rangle = l \}$
are bijective.

 In \cite{KMN1}, they give examples of perfect crystals of arbitrary levels  
for algebras of the following types:
$A_n^{(1)}$, $B_n^{(1)}$, $C_n^{(1)}$, $D_n^{(1)}$, $A_{2n}^{(2)}$, 
$A_{2n-1}^{(2)}$, $D_{n+1}^{(2)}$. 
Up to now we do not know simple criterion for perfectness of crystals. 
Further, perfect crystals with arbirary level for $E_n^{(1)}$, $F_4^{(1)}$, 
$G_2^{(1)}$, $E_6^{(2)}$ are not known. 
In this paper, we give a series of perfect crystals of $\U\gone$,
following method in \cite{KMN1}.

Let $\Up\gone$ be the subalgebra of $\U\gone$, generated by $e_i$, $f_i$,
$q^h$ $(h \in P_{cl}^* = \sum_{i \in \{ 0,1,2\}} \Z h_i)$.  
We proceed in the following way.
\begin{enumerate}
\item  As in \cite{KMN1}, for a given level $l$, our strategy is to emply the 
fusion procedure in order to construct finite dimensional 
irreducible module $V_l$ of $\Up\gone$ with a crystal pseudo base 
such that its associated crystal is perfect of level $l$.
As a starting point, we take the following data. 
Let $V$ be the direct sum of the 14-dimensional module with highest weight  
$\Lambda_1$ and the trivial module over $\U{G_2}$.
The module $V$ has a crystal base which is characterized by a polarization on 
$V$. A polarization on $\U\G$-module $M$ is the symmetric bilinear such that
for $u,v \in M$ (i)$(q^hu,v)=(u,q^hv)$, 
(ii) $(e_iu,v)=(u,q_i^{-1}t_i^{-1}f_iv)$, (iii) $(f_iu,v)=(u,q_i^{-1}t_ie_iv)$ 
 and (iv) positive definite. 
Let $B^1$ be the associated crystal of $\U{G_2}$. 
We define the actions of $e_0$,
$f_0$, and $q^{h_0}$ on $V$ to make it an irreducible module for $\Up\gone$.
Then we show that there exists a polarization for $\Up\gone$.
We see that $B^1$ is a perfect crystal of level $1$ by direct calculation.
Using the global base \cite{K3}, we compute the R-matrix for V explicitly.
By the fusion construction, we obtain a certain finite-dimensional 
submodule $V_l$ of $V^{\otimes l}$. The existence of polarization on $V_l$ 
ensures that the submodule $V_l$ has a crystal pseudo base.
The associated crystal of $V_l$ is isomorphic to a crystal of $\U{G_2}$ and a
crystal of $\U{A_2}$ as crystals of $\U{G_2}$ and $\U{A_2}$ respectively. 
These isomorphisms show that the associated crystal is a crystal of 
$\Up\gone$. It is a perfect crystal of level $l$.
This is shown in the section \ref{section:existance}. 

\item In section \ref{section:perfectness}, 
we construct a perfect crystal $B^l$ of $\U\gone$ in a combinatorial way. 

The Dynkin diagram of $\gone$ is
\unitlength 0.3mm
\begin{picture}(60,15)(0,3)
\gdynkin
\put(10,13){\makebox(0,0){\scriptsize $0$}}
\put(30,13){\makebox(0,0){\scriptsize $1$}}
\put(53,13){\makebox(0,0){\scriptsize $2$}}
\end{picture}.
Here, let $J_0=\{1,2\}$ and $J_2=\{1,0\}$ be the index set of $G_2$ and 
$A_2$, respectively. We define $\imath_i:J_i \rightarrow I$ by 
$\imath_i(j)=j$ $(i=0,2)$.
We define crystal $\CG^l$ to be 
\[
\CG^l  =  \displaystyle\bigoplus_{n = 0}^{l} B^{G_2}(n \Lambda_1) 
\mbox{ (as crystals for $\U{G_2}$ )}.
\]
  
The affine crystal $b^l$ for $\U\gone$ is constructed with $\CG^l$ and $\tf{0}$.
To define $\tf{0}$ on $\CG^l$, we use crystals of $\U{A_2}$.
Here, we define $\CA$ by
\begin{eqnarray*}
\CA_i & = & \displaystyle\bigoplus_{i \leq j_1,j_0 \leq l-i} B^{A_2}
(j_1 \Lambda_1 + j_0 \Lambda_0 ) \mbox{ (as crystals for $\U{A_2}$)},
 \vspace{3mm}\\
\CA & = & \displaystyle\bigoplus_{i = 0}^{\left [ \frac{l}{2} \right ]}
\CA_i.  
\end{eqnarray*} 


We define operators $E_\CA$, $F_\CA$ on $\CA$ such that for $b, b' \in \CA$,
(C1) $E_\CA b=b'$, if and only if $F_\CA b'=b$,
(C2) $E_\CA \tfa{0}(b) = \tfa{0}E_\CA(b)$,
(C3) $\mbox{max}\{m \mid F_\CA^mb \neq 0\}-
\mbox{max}\{m' \mid E_\CA^{m'}b \neq 0 \}=-2\wt_1^\CA(b) -\wt_0^\CA(b)$.
Operators $E_\CA$ and $F_\CA$ are counterpart of $\teg{2}$, $\tfg{2}$ on $\CG^l$.
The crystal base $\CA$ has an involution $\mbox{C}_\CA$ such that
\[
\syme{\fa{0}{r}\fa{1}{q}\fa{0}{p}\hwvl{i}{k}{j}}=
\fa{0}{j+q-2p-r}\fa{1}{k+j-q}\fa{0}{k-q+p}\hwvl{i}{j}{k},
\] 
where $0 \leq i \leq \left [ \frac{l}{2} \right ]$, $i \leq k,j \leq l-i$,
$\hwvl{i}{k}{j}$ is the highest weight element on $\CA_i$ with weight  
$j\Lambda_1+k\Lambda_0$, $0 \leq p \leq j$, $p \leq q \leq p+k$, 
$0 \leq r \leq j+q-2p$.
By calculation, we see $\sharp \CA = \sharp \CG^l$. In view of this, we 
construct  one-to-one map $\Phi\,:\,\CA \rightarrow \CG^l$ such that 
for $b\in \CA$
(E1) $\teg{1}\Phi(b) = \Phi\left (\tea{1}b \right )$, 
(E2) $\tfg{2}\Phi(b) = \Phi \left ( F_\CA b \right )$, 
(E3) $\wt_1^\CG\left (\Phi(b) \right)=\wt_1^\CA(b)$,
(E4) $\wt_2^\CG\left (\Phi(b) \right )= -2\wt_1^\CA(b)-\wt_0^\CA(b)$,
(E5) $\teg{0}\Phi(b)=0$ (resp. $\tfg{0}\Phi(b)=0$) if and only if 
$\tea{0}(b)=0$ (resp. $\tfa{0}(b)=0$), for $b \in \CA$.
Then we define $\tf{0}$ on $\CG^l$ by $\Phi \tfa{0} \Phi^{-1}$. By weight 
consideration, (C2),
involution, etc., we see that the action of $\tf{0}$ is unique.
Therefore we have crystal base  of $\U\gone$.
\item We show that perfectness are satisfied by direct calculation 
(\S\ref{section:minimal},\S\ref{section:connect}).
\end{enumerate}
Graphs of the perfect crystal of level 1 and level 2 are given by Figure 
\ref{picture:level1} and \ref{picture:level2} respectively.
In Figure \ref{picture:level2}, for simplicity, we omit arrows of $E_\CA(b)$ 
such that $\ve{0}^\CA(b) \neq 0$.   
By definition,  we have $B^2 = B^{G_2}(2\Lambda_1) \oplus B^{G_2}(\Lambda_1) 
\oplus B^{G_2}(\phi)$ and $\CA_0 = \oplus_{0 \leq j_1, j_0 \leq 2}B^{A_2}
(j_1\Lambda_1+j_0\Lambda_0)$, $\CA_1 = B^{A_2}(\Lambda_1+\Lambda_0)$,
$\CA = \CA_0 \oplus \CA_1$. In this graph, the element $b\in B^2$ from which we draw 
arrows of $E_\CA$ satisfies $\Phi(b) \in B^{G_2}(2\Lambda_1)$, 
$\Phi \left ( \tfa{0}b \right ) \not\in B^{G_2}( \Lambda_1 )$ and 
$\tea{0}b=0$.
For $b\in \CA$ such that $\Phi(b) \in B^{G_2}(\Lambda_1)\oplus B^{G_2}(\phi)$,
the action of $E_\CA$ is same as the case of level $1$. Using 
commutativity (C2), we obtain the completed graph of the perfect crystal of 
level 2.
We obtain the perfect crystal $B^1$ of level $1$ by restricting $B^2$ to
$B^{G_2}(\Lambda_1)\oplus B^{G_2}(\phi)$. From these examples, we expect that 
we can obtain the perfect crystal $B^{l+1}$ by extending $B^l$.


To find perfect crystals of level $l$, we wrote a C program which judges 
whether a given representation has perfect crystal on a combinatorial level. 
In this program, we construct crystal bases $B$
using a Lakshmibai-Seshadri path \cite{Li}.
A Lakshmibai-Seshadri path is a piecewise linear paths
$\pi : [0,1] \rightarrow (\sum \Z\Lambda_i)\otimes_Z {\bf R}$. Operators $e_\alpha$ 
and $f_\alpha$ on the path are modification by simple reflection $\sigma$ with respect to 
the simple root $\alpha$. 
When translated into the language of crystal base, a path $\pi$ is a base and operators 
$e_{\alpha_i}$ and $f_{\alpha_i}$ are Kashiwara operators $\te{i}$ and $\tf{i}$.
Let $\delta = \sum a_i\alpha_i$ be the null root. 
We define the operator $\tf{0} b$ $(b\in B)$, by searching bases with weight 
$\wt(b) + \sum_{i\neq 0} a_i\alpha_i$ and taking into
account the conditions for perfectness. The conditions are
$\langle c,\vp{}(b) \rangle \geq l$, and the map $\varphi : B_l =
\{ b \in B \mid \langle c, \vp{}(b) \rangle = l\}
\rightarrow \left (P_{cl}^+ \right)_l =
\{ \lambda \in \sum\Z_{\geq 0}\Lambda_i \mid \langle c, \lambda \rangle =l\}$
are bijective ($c$ is the element of the center). 
If we find suitable action of $\tf{0}$, $B$ is perfect crystal on a 
combinatorial level. 

Using this program we find that $\oplus_{k=0}^lB(k\Lambda_1)$ is 
a candidate of perfect crystals of  level $l$ also for $\U{F_4^{(1)}}$, 
$\U{E_6^{(1)}}$, etc, on a combinatorial level. 
We hope to report on these in a future publication. 

The author would like to acknowledge professor Etsuro Date and T. Nakashima for 
their helpful advices.

\unitlength 1mm
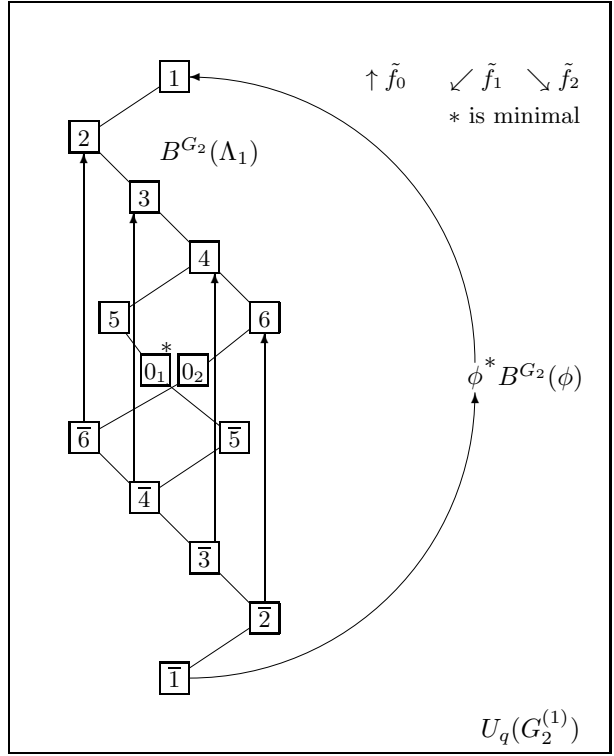
\begin{figure}[htb]
\begin{center}
\begin{picture}(80,100)
\crystallevelone
\put(20,80){\makebox(0,0)[l]{$B^{G_2}( \Lambda_1 )$}}
\put(65,50){\makebox(0,0)[l]{$B^{G_2}( \phi )$}}
\end{picture}
\end{center}
\caption{Perfect crystal of level 1 } \label{picture:level1}
\end{figure}

%% file: result.tex
Denote the index set of the simple roots of $\gone$ by $I=\{0,1,2\}$.
We can take $J_0=\{1,2\}$ and $J_2=\{1,0\}$ as the index set of $G_2$ and 
$A_2$, respectively. We define $\imath_i:J_i \rightarrow I$ by 
$\imath_i(j)=j$ $(i=0,2)$. \vspace{5mm}

\begin{proposition} \label{proposition:construction}
\rm For any integer $l \geq 1$, there exists the crystal $B^l$ for $\gone$
such that
\[
\imath_0^*(B^l) \cong  \bigoplus_{k=0}^l 
B^{G_2}(k\Lambda_1) \mbox{ } ( \mbox{ as crystals of } \U{G_2} ), 
\]
and
\[
 \imath_2^*(B^l) \cong 
\bigoplus_{i=0}^{[ \frac{l}{2} ]}
\bigoplus_{i \leq j_1,j_0 \leq l-i}
B^{A_2}(j_1 \Lambda_1 + j_0 \Lambda_0) \mbox{ } ( \mbox{ as crystals of } 
\U{A_2} ).\vspace{5mm}
\]
\end{proposition}

\begin{theorem} \label{theorem:perfectness}
The crystal $B^l$ is a perfect crystal of level $l$ for $\U\gone$.
\end{theorem}

%% file: crystal.tex
\subsection{Crystal base}
We recall definitions of quantum enveloping algebras $\U{A_2}$, $\U{G_2}$, 
$\U{\gone}$ and crystal bases of $\U{A_2}$, $\U{G_2}$.

\subsubsection{Definition of $\U{\G}$ }
Let $\G$ be a semi simple (resp. affine ) Lie algebra generated by $e_i$, $f_i$ 
$(i \in I = \{1, \ldots, n\} \:(resp. \{ 0, 1, \ldots, n\})\,)$ and {\frak h}
the Cartan subalgebra  over \Q. Let $\{ \alpha_i \mid i \in I \} \subset 
\mbox{\frak h}^*$ and $\{ h_i \mid i \in I \} \subset \mbox{\frak h}$ 
denote the simple roots and simple coroots, respectively.  
We denote a non-degenerate  invariant symmetric
bilinear form on $\mbox{\frak h}^*$  by $( , )$. This bilinear form $( , )$ 
satisfies $( \alpha_i ,\alpha_i ) \in {\bf Z}_{>0}$. 
Let $P$ be the weight lattice 
and $P^*$ be the dual lattice. 

The quantized universal enveloping algebra $\U\G$ is the 
$\Q(q)$-algebra generated by the symbols $e_i$, $f_i$ $(i \in I)$ 
and $q^h$ $(h \in P^*)$ with the following defining relations:
\[
q^0=1, q^hq^{h'} = q^{h+h'} \mbox{ for all } h, h' \in P^*,
\]
\[
q^he_iq^{-h} = q^{ \langle h, \alpha_i \rangle } e_i,
q^hf_iq^{-h} = q^{ -\langle h, \alpha_i \rangle } f_i \mbox{ for all }
i \in I, h\in P^*,
\]
\[
[e_i, f_i] = \delta_{ij} \frac{t_i-t_i^{-1}}{q_i-q_i^{-1}}, 
\mbox{ where } q_i = q^{(\alpha_i,\alpha_i)} \mbox{ and }
t_i = q^{(\alpha_i, \alpha_i)h_i} = q_i^{h_i},
\]
\[
\sum^b_{k=0}(-1)^ke_i^{(k)}e_je_i^{(b-k)} = 0
\mbox{ for } i \neq j \mbox{ and } b = 1 - \langle h_i, \alpha_j \rangle. 
\]
\[
\sum^b_{k=0}(-1)^kf_i^{(k)}f_jf_i^{(b-k)} =0
\mbox{ for } i \neq j \mbox{ and } b = 1 - \langle h_i, \alpha_j \rangle. 
\]
Here we use the following notations: 
$[m]_i = \frac{q_i^m-q_i^{-m}}{q_i-q_i^{-1}}$,
$[k]_i!=\raisebox{-1mm}{$\mathop{\mbox{\Large$\Pi$}}_{m=1}^k$}[m]_i$, 
and $e_i^{(k)} = e_i^k/[k]_i!$, $f_i^{(k)} = f_i^k/[k]_i!$. 
For $k < 0$, we put $e_i^{(k)} = f_i^{(k)} = 0$.

It is well-known that $\U\G$  has a Hopf algebra srtucture with a 
comultiplication $\Delta$ defined by
\begin{eqnarray*}
\Delta (e_i) &=& e_i \otimes t_i^{-1} + 1 \otimes  e_i,\\
\Delta (f_i) &=& f_i \otimes 1 + t_i \otimes f_i, \\
\Delta (q^h) &=& q^h \otimes  q^h.
\end{eqnarray*}
for all $i \in I$, $h \in P^*$. The tensor product of two $\U\G$-modules
has a stucture  of $\U\G$-module by this comultiplication.

\subsubsection{Definition of crystal bases[3]}
\label{def:kashiwara}
For $\U\G$-module $M$ and $\lambda \in P$, the $\lambda$-weight space of $M$ 
is defined by $M_{\lambda} = \{ u \in M \mid  q^hu = q^{\langle h, \lambda 
\rangle} u \mbox{ for all } h \in P^* \}$.
For $J \subset I$, let $\U{\G_J}$ be the subalgebra of  
$\U\G$ generated by $e_i$, $f_i$, $t_i$ and $t_i^{-1}$ $(i \in J)$.
We say that $M$ is {\it integrable} if $M = \bigoplus_{\lambda \in P}
M_\lambda$ and $M$ is a union of finite-dimensional $\U{\G_{\{i\}}}$-modules
for each $i \in I$.  By the representation theory of  $\U{\mbox{\frak sl }_2}$
 any element $u \in M_{\lambda}$ can be uniquely written as 
\[
u = \sum_{k \geq 0}f_i^{(k)}u_k,
\]
where $u_k \in \ker e_i \cap M_{\lambda+k\alpha_i}$.
We define the Kashiwara operators $\te{i}$ and $\tf{i}$ on $M$ as follows. 
\[
\te{i}u = \sum f_i^{(k-1)}u_k,
\]
\[
\tf{i}u = \sum f_i^{(k+1)}u_k.
\]

Let ${\bf A}$ be the subring of $\Q(q)$ consisting of the rational functions 
regular at $q = 0$. A {\it crystal lattice} of integrable $\U\G$-module $M$  
is a free $A$-submodule of $M$ such that (1) $M \cong \Q(q) \otimes_A L$,
(2) $L = \oplus_{\lambda\in P}L_\lambda$ where $L_\lambda = L \cap M_\lambda$,
(3) $\te{i}L \subset L$, $\tf{i}L \subset L$.
A {\it crystal base} of the integrable $\U\G$-module $M$ is a pair $(L,B)$ 
such that  
\renewcommand{\labelenumi}{(\makebox[3mm]{\roman{enumi}}) }
\begin{enumerate}
\item $L$ is a crystal lattice of $M$, 
\item $B$ is a \Q-base of $L / qL$,
\item $B = \sqcup_{\lambda \in P} B_{\lambda}$ where 
$B_\lambda = B \cap ( L_\lambda / qL_\lambda )$,
\item $\te{i}B \subset B \sqcup \{ 0 \}$, $\tf{i}B \subset B \sqcup \{ 0 \}$,
\item for $b$, $b' \in B$, $b' \in \tf{i}b$ if and only if $b = \te{i}b'$ for 
$i \in I$.
\end{enumerate}
We sometimes replace (ii) by: 
$B_{ps} = B' \sqcup (-B')$ where $B'$  is a \Q-base of $L/qL$.
We call $(L, B_{ps})$ a {\ crystal psuedo-base} and $B_{ps}/\{\pm1\}$ the 
associated crystal of $(L, B_{ps})$.

For a dominant integrable weight $\lambda\in P_+ = 
\{\lambda \in P \mid \langle h_i,\lambda \rangle \leq 0
\mbox{ for all } i \}$, let $V(\lambda)$ denote the irreducible integrable 
$\U\G$-module with highest weight $\lambda$. Let $u_\lambda$ be the highest 
weight vector of $V(\lambda)$, and let $L(\lambda)$ be the smallest A-submodule
of $V(\lambda)$ containing $u_\lambda$ and stable under $\tf{i}$'s. 
Set $B(\lambda) = \{b \in L(\lambda) / qL(\lambda) \mid b=\tf{i_1} \cdots
\tf{i_l} u_{\lambda} \mbox{ \rm mod } qL(\lambda)\} \setminus \{0\}$.
Than $(L(\lambda),B(\lambda))$ is crystal  base of $V(\lambda)$.
Crystal graph is an oriented colored ( by I ) graph with $B(\lambda)$ as the 
set of vertices and colored arrow $b \stackrel{i}{\longrightarrow} b'$ 
if and only if $b' = \tf{i}b$ (hence $\te{i}b'=b$). 
Then graph completely describes the actions of $\te{i}$
and $\tf{i}$ on $B(\lambda)$. For $b \in B(\lambda)$, we set
\[
\ve{i}(b) = \max \{ k \geq 0 \mid \te{i}^kb \neq 0 \},
\]
\[
\vp{i}(b) = \max \{ k \geq 0 \mid \tf{i}^kb \neq 0 \}.
\]
Note that $B(\lambda) = \sqcup_{\mu \in P}B(\lambda)_{\mu}$ and for 
$b \in B(\lambda)_\mu$, we have $\langle h_i, \wt(b) \rangle = 
\vp{i}(b)-\ve{i}(b)$, where $\wt(b) =\mu$ denotes the weight of b.
\begin{proposition}\label{proposition:sayou}
Let $(L_j,B_j)$ be crystal bases of integrable $\U\G$-modules $M_j$ 
$(j=1,2)$. Set $L = L_1 \otimes _{\bf A} L_2 \subset M_1 \otimes M_2$  and
$B = \{ b_1 \otimes b_2 \mid b_j \in B_j \; (j=1,2) \} \subset L / qL\}$.
Then $(L,B)$ is the crystal base of  $M_1 \otimes M_2$  and the action of
$\tf{i}$, $\te{i}$ is given by  
\begin{eqnarray*}
\tf{i}(b_1 \otimes b_2) &=& \left \{
\begin{array}{ll} 
\tf{i}b_1 \otimes b_2 & \mbox{ if } \vp{i}(b_1) > \ve{i}(b_2),\\
b_1 \otimes \tf{i}b_2 & \mbox{ if } \vp{i}(b_1) \leq \ve{i}(b_2),
\end{array}
\right . \\
\te{i}(b_1 \otimes b_2) &=& \left \{
\begin{array}{ll} 
\te{i}b_1 \otimes b_2 & \mbox{ if } \vp{i}(b_1) \geq \ve{i}(b_2),\\
b_1 \otimes \te{i}b_2 & \mbox{ if } \vp{i}(b_1) < \ve{i}(b_2).
\end{array}
\right . 
\end{eqnarray*}
\end{proposition}
\begin{corollary}
For $b_j \in B_j$ $(j=1,2)$, we have
\begin{eqnarray*}
\ve{i}(b_1 \otimes b_2) &=& 
\max(\ve{i}(b_1),\ve{i}(b_2)+\ve{i}(b_1)-\vp{i}(b_1)),\\
\vp{i}(b_1 \otimes b_2) &=& 
\max(\vp{i}(b_2),\vp{i}(b_1)+\vp{i}(b_2)-\ve{i}(b_2)).
\end{eqnarray*}
\end{corollary}
\begin{proposition}\label{proposition:politensor}
Let $(L_j,B_j)$ be crystal bases of integrable $\U\G$-modules $M_j$ 
$(j=1,\ldots,N)$, and let $\{a_k\}$ be the sequence defined by 
$a_1=0$, $a_{k+1}=a_k+\vp{i}(b_k)-\ve{i}(b_{k+1})$.
Then we have 
\[
\tf{i}(b_1 \otimes \cdots \otimes b_N)=b_1 \otimes \cdots \otimes b_{k-1} 
\otimes \tf{i}b_k \otimes b_{k+1} \otimes \cdots \otimes b_N,
\]
where $k=\mbox{max}\{j \mid a_j=\mbox{min}\{a_1, \ldots, a_N\}\}$,
\[
\te{i}(b_1 \otimes \cdots \otimes b_N)=b_1 \otimes \cdots \otimes b_{l-1} 
\otimes \te{i}b_l \otimes b_{l+1} \otimes \cdots \otimes b_N,
\]
where $l=\mbox{min}\{j \mid a_j=\mbox{min}\{a_1, \ldots, a_N\}\}$.
\end{proposition}
We regard integrable $\U\G$-module $M$ as a union of finite dimensional 
$\U{\G_{\{i\}}}$-module
for any $i \in I$. 
Let $u_0$ be the base of trivial representation of 
$\U{\mbox{\frak sl }_2}$-module, and let $u_+$, $u_-$ be a base of 
2 dimensional $\U{\mbox{\frak sl }_2}$-module which satisfy 
\[
fu_+ = u_-,
\]
where $f$ is one of generators of $\U{\mbox{\frak sl }_2}$. \\
In order to calculate the action of $\te{i}$ and  $\tf{i}$ on tensor products, 
we will use $u_+$, $u_-$ and $u_0$.
For that purpose we prepare the following. 
Let $B$ be a crystal base. 
Let $\fW$ be the set of words generated by $u_+$, $u_-$, $u_0$. 
\begin{definition}
\label{def:uui}
For $i \in I$  we define a map 
\[
\uu{i}: B \rightarrow \fW 
\]
as follows:
\begin{enumerate}
\renewcommand{\labelenumi}{\rm (\arabic{enumi}) }
\item for $b\in B$ such that $\te{i}b=\tf{i}b=0$, we define  
\[
\uu{i}(b)=u_0,
\]
\item if $b_k$ is an element of a maximally connected i-colored crystal subgraph
$b_1 \stackrel{i}{\longrightarrow} b_2 \stackrel{i}{\longrightarrow} \cdots
\stackrel{i}{\longrightarrow} b_n$
where $n \in {\bf Z}_{\geq 2}$, we define
\[
\uu{i}(b_k)=u_-^{k-1} u_+^{n-k}.
\]
\end{enumerate}
We extend $\uu{i}$ to $b'_1 \otimes \cdots \otimes b'_N \in 
B^{\otimes N}$, $N\in {\bf Z}_{>0}$, by  
\[
\uu{i}(b'_1\otimes \cdots \otimes b'_N)
=\uu{i}(b'_1) \cdots \uu{i}(b'_N) = u_1 \cdots u_N,
\] 
where, $\uu{i}(b'_k)=u_k\in \fW$.
\end{definition}
For $v=v_1\cdots v_{N'} \in \fW$, we define the length of the word $v$ by 
$l(v)=N'$. Let $\fN$ be the set of words generated by ${\bf Z}_{>0}$.
\begin{definition}
We define a map $\uu{\fNs}:B^{\otimes N} \rightarrow \fN$ by 
\[
\uu{\fNs}(b'_1\otimes \cdots \otimes b'_N)=
1^{l(\uu{i}(b_1))}2^{l(\uu{i}(b_2))}\cdots N^{l(\uu{i}(b_N))}
\]
\end{definition}
For $r=r_1\cdots r_{N''} \in \fN$, we also define the length of the word $r$ by 
$l(r)=N''$. 

We put $\fW_j=\{v\in \fW \mid l(v)=j \}$, $\fN_j=\{r\in \fN \mid l(r)=j\}$,
where $j\in \Z_{\geq 0}$.

Take words $v=v_1 \cdots v_M \in \fW_M$, $r=r_1 \cdots r_M \in \fN_M$,
where $M\in \Z_{\geq 0}$.
We put 
\[
s=\sharp\{s'\mid v_{s'}=u_0, 1 \leq s' \leq M\},
\]
\[
\{m_1, m_2, \ldots, m_s\}=\{s' \mid v_{s'}=u_0, 1 \leq s' \leq M\},
\]
where $1 \leq m_1 < m_2 < \cdots < m_s \leq M$.
\begin{definition}
\label{def:red0}
We define a map 
\[
\mbox{\rm Red}_0 : \bigcup_{j\in \Z_{\geq 0}}\left (\fW_j \times 
\fN_j\right ) \rightarrow \bigcup_{j\in \Z_{\geq 0}}\left ( \fW_j \times \fN_j
\right ).
\]
The map $\mbox{Red}_0$ delete all $v_{m_i}$ which equal to $u_0$ in $v$ and 
$r_{m_i}$ in $r$.
Namely we can express
\[
\mbox{\rm Red}_0(v,r)= (
v_1\cdots v_{{m_1}-1} v_{{m_1}+1} \cdots v_{{m_s}-1}v_{{m_s}+1}\cdots v_M,
r_1\cdots r_{{m_1}-1} r_{{m_1}+1} \cdots r_{{m_s}-1}r_{{m_s}+1}\cdots r_M).
\]
\end{definition}
We define 
\[
s_0=\left \{
\begin{array}{ll}
\mbox{\rm min}\{m \mid v_m=u_+, v_{m+1}=u_-\}& 
\mbox{if there exists $m$ such that } v_m=u_+,\;v_{m+1}=u_-,\\
0& \mbox{otherwise}.
\end{array}
\right .
\]
\begin{definition}
We define a map 
\[
\mbox{\rm Red}_+ : \bigcup_{j\in \Z_{\geq 0}}\left ( \fW_j \times \fN_j \right ) 
\rightarrow \bigcup_{j\in \Z_{\geq 0}}\left ( \fW_j \times \fN_j \right ),
\]
for $v=v_1\cdots v_M\in \fW_M$, $r=r_1\cdots r_M\in \fN_M$ by
\[
\mbox{\rm Red}_+(v,r)=\left \{
\begin{array}{ll}
(v_1\cdots v_{s_0-1} v_{s_0+2}\cdots v_M,
r_1\cdots r_{s_0-1} r_{s_0+2}\cdots r_M) &\mbox{if } s_0\geq 1,\\
(v,r)&\mbox{if } s_0=0.
\end{array} \right .
\]
\end{definition}
\begin{definition}
\label{def:reduce}
We define a map 
\[
\mbox{\rm Red} : \bigcup_{j\in \Z_{\geq 0}}\left ( \fW_j \times \fN_j \right ) 
\rightarrow \bigcup_{j\in \Z_{\geq 0}}\left ( \fW_j \times \fN_j \right )
\]
by 
\[
\mbox{\rm Red}=\mbox{\rm Red}_+^n\mbox{\rm Red}_0,
\]
where $n\in {\bf Z}_{\geq 0} \mbox{ such that } n=
\mbox{\rm min}\{n'\mid\mbox{Red}_+^{\,n'+1}(\mbox{Red}_0(v,r))=
\mbox{Red}_+^{\,n'}(\mbox{Red}_0(v,r))\}$.
\end{definition}

By Definition \ref{def:reduce}, we can denote
\[
\mbox{Red}\left ( \uu{i}(b),\uu{\fNs}(b) \right ) = 
(v'_1 v'_2 \cdots v'_{M'}, r'_1 r'_2 \cdots r'_{M'})
\]
where $M' \in {\bf Z}_{\geq 0}$.
Let $\mbox{\rm Pr}_\fWs$, $\mbox{\rm Pr}_\fNs$ be a projection of 
$\cup_{j\in Z}(\fW_j \times \fN_j)$ to $\fW$, $\fN$ respectively. 
We define a map
\[
\mbox{\rm Red}_i : B^{\otimes N} \rightarrow \fW,
\]
by
\[
\mbox{\rm Red}_i=\mbox{\rm Pr}_\fWs \cdot \mbox{\rm Red}(\uu{i}, \uu{\fNs})
\]
We define a map
\[
\mbox{\rm Red}_\fNs : B^{\otimes N} \rightarrow \fN,
\]
by
\[
\mbox{\rm Red}_\fNs=\mbox{\rm Pr}_\fNs \cdot \mbox{\rm Red}(\uu{i}, \uu{\fNs})
\]
We set  
\[
s_+ = \left \{ 
\begin{array}{ll}
\mbox{min}\{s' \mid v'_{s'}=u_+\}& 
\mbox{if there exists $s'$ such that } v'_{s'}=u_+,\\
0& \mbox{otherwise},
\end{array}
\right .
\]
\[
s_- = \left \{ 
\begin{array}{ll}
\mbox{\rm max}\{s' \mid v'_{s'}=u_-\} & 
\mbox{if there exists $s'$ such that } v'_{s'}=u_-\\
0 & \mbox{otherwise}.
\end{array}
\right .
\]
\begin{proposition}
\label{proposition:vplus1}
For $b=b_1 \otimes \cdots \otimes b_N \in B^{\otimes N}$, we have
\[
\tf{i}(b)= 
\left \{
\begin{array}{ll}
b_1 \otimes \cdots \otimes \tf{i}b_{r_{s_+}} \otimes \cdots \otimes b_N& 
\mbox{ if } s_+>0,\\
0& \mbox{ if } s_+=0,
\end{array}
\right .
\]
\[
\te{i}(b)= 
\left \{
\begin{array}{ll}
b_1 \otimes \cdots \otimes \te{i}b_{r_{s_-}} \otimes \cdots \otimes b_N& 
\mbox{ if } s_->0,\\
0& \mbox{ if } s_-=0.
\end{array}
\right .
\]
\end{proposition}  
For $b=b_1 \otimes \cdots \otimes b_N \in 
B^{\otimes N}$, we denote 
$\uu{i}(b_k)=u_-^{\ve{i}(b_k)}u_+^{\vp{i}(b_k)}$.
Then we see that Proposition \ref{proposition:politensor} and 
Proposition \ref{proposition:vplus1} are equivalent.  

\begin{remark}
\label{rem:uveuvp}
For $b=b_1\otimes \cdots \otimes b_N$, we denote 
$\mbox{\rm Red}_i(b) = u_-^{m_1} u_+^{m_2}$ $(m_1, \:m_2 \in {\bf Z}_{\geq 0})$.
Then by Definition \ref{def:uui} and Proposition \ref{proposition:vplus1}, we see 
\[
\ve{i}(b)=m_1,
\]
\[
\vp{i}(b)=m_2.
\] 
\end{remark}
\begin{example}
We assume $b = b_1 \otimes b_2 \otimes b_3$, $\uu{i}(b_1)= u_1=u_- u_+ u_+$, 
$\uu{i}(b_2)=u_2=u_0$ and $\uu{i}(b_3)=u_3=u_- u_+$.
We calculate $\tf{i}b$ as follows:   
\begin{eqnarray*} 
\uu{i}(b) &=& \uu{i}(b_1) \uu{i}(b_2) \uu{i}(b_3)\\
&=& u_1 u_2 u_3\\
&=& u_- u_+ u_+ u_0 u_- u_+,\\
\uu{\fNs}(b)&=&111233,\\
\mbox{\rm Red}_0(\uu{i}(b),\uu{\fNs}(b))&=& 
(u_- u_+ u_+ u_- u_+,11133), \\
\mbox{\rm Red}(\uu{i}(b),\uu{\fNs}(b))&=& (u_- u_+ u_+,113), \\
\mbox{\rm Red}_i(b)&=&u_-u_+u_+.
\end{eqnarray*}
Then we have $s_+=2$, and that $s_+$-th integer of the word $113$ is $1$.
Therefore,
\[
\tf{i} b = \left ( \tf{i} b_1 \right ) \otimes b_2 \otimes b_3.
\]
\end{example}

\subsubsection{Crystal bases for $\U{A_2}$-module[6]}
Consider the finite  dimensional simple Lie algebra of type $A_2$ with Cartan
matrix 
\raisebox{2mm}{
$\left ( 
\begin{array}{cc}
2 & -1\\
-1 & 2
\end{array}
\right )$},
Let us denote by $\{ \alpha_1, \alpha_2\}$ the set of simple roots and 
$\{h_1,h_2\}$ the set of simple coroots. 
Define $\Lambda_i \in \mbox{\frak h}^*$ $(i = 1,2)$
by  $\Lambda_i(h_j) = \delta_{ij}$, $(j=1,2)$. We put $P=\Z\Lambda_1 \oplus 
\Z\Lambda_2$, $P^* = \Z h_1 \oplus \Z h_2$.
Thus any dominant integral weight $\Lambda \in P_+$ is of the form 
$\Lambda = m\Lambda_1 + n\Lambda_2$, $(n,m \in \Z_+)$.
We define a non-degenerate symmetric bilinear form $( , )$ on 
$\mbox{\frak h}^*$ by $( \alpha_1, \alpha_1 ) = ( \alpha_2, \alpha_2 ) = 1$,
$( \alpha_1, \alpha_2 ) = ( \alpha_2, \alpha_1 ) = -\frac{1}{2}$,
$\langle h_i, \alpha_j \rangle$ is the $(i,j)$-th element of Cartan matrix and 
set $t_i = q^{h_i}$ for $i=1,2$.
Then the corresponding quantized universal enveloping algebra $\U{A_2}$ is the
associative algebra over $\Q(q)$ generated by $e_i$, $f_i$, $t_i$, 
$t_i^{-1}$ $(i = 1,2)$ satisfying the relations:
\[
t_it_j = t_jt_i, \hspace{5mm} t_i^{}t_i^{-1}= t_i^{-1}= t_i^{}=1,
\]
\[
t_i^{}e_j^{}t_i^{-1}=q^{\langle h_i, \alpha_j \rangle } e_j,
t_i^{}f_j^{}t_i^{-1}=q^{-\langle h_i, \alpha_j \rangle } f_j,
\]
\[
[e_i,f_j] = \delta_{ij}\frac{t_i^{}-t_i^{-1}}{q-q^{-1}},
\]
\[
e_i^{(2)}e_j - e_ie_je_i + e_je_i^{(2)} = 
f_i^{(2)}f_j - f_if_jf_i + f_jf_i^{(2)} = 0 \; (i \neq j).
\] 

Let $V$ be the 3-dimensional  $\Q(q)$-vector space with a basis 
$\left \{ \vmini{1}, \vmini{2}, \vmini{3} \right \}$. 
We define a $\U{A_2}$-module structure on $V$ as follows:
\[
\begin{array}{lll}
e_1\vmini{1}=0, & e_1\vmini{2} = \vmini{1}, & e_1\vmini{3} = 0,\\
e_2\vmini{1}=0, & e_2\vmini{2} = 0, & e_2\vmini{3} = \vmini{2},\\
f_1\vmini{1}=\vmini{2}, & f_1\vmini{2} = 0, & f_1\vmini{3} = 0,\\
f_2\vmini{1}=0, & f_2\vmini{2} = \vmini{3}, & f_2\vmini{3} = 0,\\
t_1\vmini{1}=q\;\vmini{1}, & t_1\vmini{2} = q^{-1}\:\vmini{2}, & 
t_1\vmini{3} = \vmini{3},\\
t_2\vmini{1}=\vmini{1}, & t_2\vmini{2} = q\;\vmini{2}, & 
t_2\vmini{3} = q^{-1}\:\vmini{3}.
\end{array}
\]
The module $V$ is isomorphic to the integrable highest weight $\U{A_2}$-module 
$V(\Lambda_1)$ with highest weight $\Lambda_1$ and highest weight vector \vmini{1}.
\begin{proposition}
Put 
\[
L(\Lambda_1) = {\bf A}\vmini{1} \oplus {\bf A}\vmini{2} \oplus {\bf A}
\vmini{3},
\]
\[
B(\Lambda_1) = \left \{ \vmini{1} , \: \vmini{2} , \: \vmini{3} \right \}.
\]
Then $(L(\Lambda_1),B(\Lambda_1))$ is the crystal base of $V(\Lambda_1)$. 
The crystal graph $B(\Lambda_1)$ is given by:
\[
\vmini{1} \stackrel{1}{\longrightarrow} \vmini{2} \stackrel{2}{\longrightarrow}
\vmini{3}.
\]
\end{proposition}
Define an ordering on the set $B(\Lambda_1)$ by
\[
\vmini{1} \prec \vmini{2} \prec \vmini{3}.
\]
We write \vtate{a}{b} for $\vminiD{a}\otimes \vminiD{b}$ $(a \prec b )$. 
Similar to the $B(\Lambda_1)$-case, the crystal graph $B(\Lambda_2)$ is 
given by:
\[
\vtate{1}{2} \stackrel{2}{\longrightarrow} \vtate{1}{3} 
\stackrel{2}{\longrightarrow} \vtate{2}{3}.
\]
\begin{proposition}
\label{prop:A2crystal}
The crystal bases of $B(m\Lambda_1+n\Lambda_2)$ is given by the set of 
semi-standard Young tableaux of the same shape 
\[
\begin{array}{l}
B(m\Lambda_1+n\Lambda_2)\\
=\left \{ \left . \younglambda \; \right |
b^j_i \in \{1,2,3\},b^j_i \leq b_i^{j-1}, b_1^{j} < b_2^{j} \right \}.
\end{array}
\]
In particular, highest weight element 
$\hat{b} \in B(m\lambda_1+ n\Lambda_2)$ is given by 
\[
\hat{b}=\youngHigh.
\]
\end{proposition}
By Definition \ref{def:uui} we have
\[
\uu{i}\left ( \vminiD{j} \right )=
\left \{
\begin{array}{ll}
u_+& \mbox{ if } i=j,\\
u_-& \mbox{ if } i+1=j,\\
u_0& \mbox{otherwise}.
\end{array}
\right .
\]
An element $b$ of $B(m\Lambda_1 + n\Lambda_2)$ can be 
expressed in the form 
\[
\vminiD{b^1_1} \otimes \cdots \vminiD{b^{\!m}_1} \otimes \left ( \; 
\vminilD{b^{\!m\!+\!1}_1} \otimes \vminilD{b^{\!m\!+\!1}_2} \; \right ) 
\otimes \cdots \otimes \left ( \; \vminilD{b^{\!m\!+\!n}_1} \otimes
\vminilD{b^{\!m\!+\!n}_2} \; \right ) \in B(\Lambda_1)^{\otimes m+2n}.
\]
Then using Proposition \ref{proposition:vplus1}, we have the action of $\te{i}$
and $\tf{i}$.

\subsubsection{Crystal bases for $\U{G_2}$-module[7]}
Consider the finite  dimensional simple Lie algebra of type $G_2$ with Cartan
matrix 
\raisebox{2mm}{
$\left ( 
\begin{array}{cc}
2 & -1\\
-3 & 2
\end{array}
\right )$}.\\
Let $\{ \alpha_1, \alpha_2\}$, $\{h_1,h_2\}$ be the set of simple root and simple coroot 
respectively. Let $\mbox{\frak h}$, $\Lambda_i$ $(i=1,2)$, $P$ and $P^*$ be Caltan 
subalgebra, fundermental weight, weight lattice and dual lattice, respectively. 
We define a non-degenerate symmetric bilinear form $( , )$ on $\mbox{\frak h}^*$
by $( \alpha_1, \alpha_1 ) = 3$, $( \alpha_2, \alpha_2 ) = 1$,
$( \alpha_1, \alpha_2 ) = ( \alpha_2, \alpha_1 ) = -\frac{3}{2}$,
$\langle h_i, \alpha_j \rangle$ is $(i,j)$-element of Cartan matrix and 
set $q_i = q^{(\alpha_i,\alpha_i)}$, $t_i = q^{h_i}$ for $i=1,2$.
Then the corresponding quantized universal enveloping algebra $\U{G_2}$ is the
associative algebra over $\Q(q)$ generated by $e_i$, $f_i$, $t_i$, 
$t_i^{-1}$ $(i = 1,2)$ satisfying the relations:
\[
t_it_j = t_jt_i, \hspace{5mm} t_i^{}t_i^{-1}= t_i^{-1}= t_i^{}=1,
\]
\[
t_i^{}e_j^{}t_i^{-1}=q_i^{\langle h_i, \alpha_j \rangle } e_j,
t_i^{}f_j^{}t_i^{-1}=q_i^{-\langle h_i, \alpha_j \rangle } f_j,
\]
\[
[e_i,f_j] = \delta_{ij}\frac{t_i^{}-t_i^{-1}}{q_i-q_i^{-1}},
\]
\[
e_1^{(2)}e_2 - e_1e_2e_1 + e_2e_1^{(2)} = 
f_1^{(2)}f_2 - f_1f_2f_1 + f_2f_1^{(2)} = 0, 
\] 
\[
\sum_{n=0}^4(-1)^ne_2^{(n)}e_1e_2^{(4-n)} =
\sum_{n=0}^4(-1)^ne_2^{(n)}e_1e_2^{(4-n)} = 0.
\]
Let $V$ be the 14-dimensional  $\Q(q)$-vector space with a basis 
\[
\left \{ \left . \vminiD{i}, \vminiD{\overline{i}} \: \right |\: i=1,\ldots,6 
\right \} \cup \left \{ \vmini{7}, \vmini{8} \right \}.
\] 
We define a $\U{G_2}$-module structure on $V$ as follows:\\
\sayoue{1}{2}{1}{1}, \sayoue{1}{5}{4}{1}, \sayoue{1}{8}{6}{[2]_2}, 
\sayoue{1}{-6}{8}{1} \sayounokori{1}{[2]_2}{7}, \vspace{2mm}\\ 
\sayoue{1}{-4}{-5}{1}, \sayoue{1}{-1}{-2}{1}, \vspace{2mm}\\  
\sayoue{2}{3}{2}{[3]_2}, \sayoue{2}{4}{3}{[2]_2}, \sayoue{2}{6}{4}{1} ,
\sayoue{2}{7}{5}{[2]_2}, \vspace{2mm}\\ 
\sayoue{2}{-5}{7}{1} \sayounokori{[3]_2}{[2]_1}{8}, \sayoue{2}{-4}{-6}{[3]_2},
\sayoue{2}{-3}{-4}{[2]_2}, \sayoue{2}{-2}{-3}{1}, 
\vspace{2mm}\\
\sayouf{1}{1}{2}{1}, \sayouf{1}{4}{5}{1}, \sayouf{1}{6}{8}{1} 
\sayounokori{1}{[2]_2}{7}, \vspace{2mm}\\
\sayouf{1}{8}{-6}{[2]_2}, \sayouf{1}{-5}{-4}{1}, 
\sayouf{1}{-2}{-1}{1}, \vspace{2mm}\\
\sayouf{2}{2}{3}{1}, \sayouf{2}{3}{4}{[2]_2}, 
\sayouf{2}{4}{6}{[3]_2}, \sayouf{2}{5}{7}{1} \sayounokori{[3]_2}{[2]_1}{8} ,
\vspace{2mm}\\
\sayouf{2}{7}{-5}{[2]_2}, \sayouf{2}{-6}{-4}{1}, \sayouf{2}{-4}{-3}{[2]_2}, 
\sayouf{2}{-3}{-2}{[3]_2},
\vspace{2mm}\\
\sayout{1}{1}{}, \sayout{1}{2}{-1}, \sayout{1}{3}{0}, \sayout{1}{4}{},  
\sayout{1}{5}{-1},
\vspace{2mm}\\
\sayout{1}{6}{-2}, \sayout{1}{7}{0}, \sayout{1}{8}{0}, 
\sayout{1}{-6}{-2},
\vspace{2mm}\\
\sayout{1}{-5}{}, \sayout{1}{-4}{-1}, \sayout{1}{-3}{0}, 
\sayout{1}{2}{}, \sayout{1}{-1}{-1},
\vspace{2mm}\\  
\sayout{2}{1}{0}, \sayout{2}{2}{3}, \sayout{2}{3}{}, \sayout{2}{4}{-1},  
\sayout{2}{5}{2},
\vspace{2mm}\\
\sayout{2}{6}{-3}, \sayout{2}{7}{0}, \sayout{2}{8}{0}, 
\sayout{2}{-6}{3},
\vspace{2mm}\\
\sayout{2}{-5}{-2}, \sayout{2}{-4}{}, \sayout{2}{-3}{-1}, 
\sayout{2}{-2}{-3}, \sayout{2}{-1}{0}.
\vspace{2mm}\\  
The other cases, $b \in B(\Lambda_1)$ is annihilated by $\te{i}$ or $\tf{i}$.

The module $V$ is isomorphic to the integrable highest weight $\U{G_2}$-module 
$V(\Lambda_1)$ with highest weight $\Lambda_1$ and highest weight vector \vmini{1}.
\begin{proposition}
\label{proposition:g2graph}
Put	
\[
L(\Lambda_1) = \bigoplus_{i=1}^6 \left ( {\bf A}\vminiD{i}\:\oplus {\bf A} 
\vminiD{\overline{i}} \right ) \bigoplus {\bf A} \vmini{7} \bigoplus {\bf A} 
\vmini{8},
\]
\[
B(\Lambda_1) = \left \{ \left . \vminiD{i}\:, \: {\bf A} \vminiD{\overline{i}} 
\: \right | i=1,2,\ldots,6 \right ) \cup \left \{ \vmini{7} , \: \vmini{8} 
\right \}.
\]
Then $(L(\Lambda_1),B(\Lambda_1))$ is the crystal base of $V(\Lambda_1)$. 
The crystal graph $B(\Lambda_1)$ is given by:
\begin{center}
\unitlength 0.5mm
\begin{picture}(220,60)
\crystalgraphGtwo 1
\put(75,19){\makebox(0,0){$2$}}\put(75,41){\makebox(0,0){$1$}}
\put(135,19){\makebox(0,0){$2$}}\put(135,41){\makebox(0,0){$1$}}
\end{picture}.
\end{center}
\end{proposition}
Define an ordering on the set $B(\Lambda_1)$ by
\[
\vmini{1} \prec \vmini{2} \prec \vmini{3} \prec \vmini{4} \prec \!\!%
\begin{array}{l}
\vmini{5}\\ \vmini{6}
\end{array}%
\!\!\prec\!\!%
\begin{array}{l}
\vmini{7}\\ \vmini{8}
\end{array}%
\!\!\prec\!\!%
\begin{array}{l}
\vmini{-5}\\ \vmini{-6}
\end{array}%
\!\!\prec \vmini{-4} \prec \vmini{-3} \prec \vmini{-2} \prec \vmini{-1}.
\vspace{5mm}
\]
We define the map sgn, $x_i$, and $\overline{x}_i$. \vspace{3mm}\\
For $a \in {\bf R}$, we define ${\rm sgn}(a)$ by
\[
{\rm sgn}(a) = \left \{
\begin{array}{ll}
\displaystyle \frac{a}{|a|} &  a \neq 0 \vspace{2mm},\\
0 & a = 0.
\end{array}
\right .
\]
For $b = \crystalbase{b_1}{b_n}$, we define $x_i(b)$, $\overline{x}_i(b)$ by 
\[
x_i(b)= \sharp \left \{ \left . \vminiD{b_k} = \vminiD{i} \:\right |
\: k = 1, \ldots ,n \right \} \;\; (i= 1,\ldots,6,0_1,0_2), 
\]
\[
\overline{x}_i(b)=\sharp \left \{ \left . \vminiD{b_k} = \vminiD{\overline{i}}
\:\right |\: k = 1, \ldots ,n \right \} \;\; (i= 1,\ldots,6). 
\]

\begin{proposition}\label{proposition:gtwocrystal}
The crystal bases of $B(n\Lambda_1)$ is given by following set of restrected 
semi-standard Young tableaux:
\[
\left \{\: 
\begin{array}{lcl}
b&=&\crystalbase{b_1}{b_l}\\
&=&\vminiD{b_l} \otimes \cdots \otimes \vminiD{b_1}
\end{array} \: \left |
\begin{array}{l}
b_k \preceq b_{k+1} \: ( k = 1, \cdots, l-1 ) \\
x_5(b)+x_{0_1}(b)+x_{0_2}+\overline{x}_5(b) \leq 1,\\
x_3(b)+x_4(b)+x_5(b) \leq 1, \\
\overline{x}_5(b)+\overline{x}_4(b)+\overline{x}_3(b)\leq 1,\\
x_5(b)+{\rm sgn}(x_6(b))+x_{0_1}(b)\leq 1,\\
x_{0_1}(b)+{\rm sgn}(\overline{x}_6(b))+\overline{x}_5(b) \leq 1
\end{array}
\right . \right \}.
\]
\end{proposition}

By Definition \ref{def:uui}, we have 
\[
\begin{array}{llll}
\uu{1}\left (\vmini{1} \right ) = u_+, & 
\uu{1}\left (\vmini{2} \right ) = u_-, &
\uu{1}\left (\vmini{3} \right ) = u_0, & 
\uu{1}\left (\vmini{4} \right ) = u_+, \vspace{2mm}\\
\uu{1}\left (\vmini{5} \right ) = u_-, & 
\uu{1}\left (\vmini{6} \right ) = u_+^2, &
\uu{1}\left (\vmini{7} \right ) = u_0, &
\uu{1}\left (\vmini{8} \right ) = u_- u_+, \vspace{2mm}\\
\uu{1}\left (\vmini{-6} \right ) = u_-^2, &
\uu{1}\left (\vmini{-5} \right ) = u_+, &
\uu{1}\left (\vmini{-4} \right ) = u_-, & 
\uu{1}\left (\vmini{-3} \right ) = u_0, \vspace{2mm}\\
\uu{1}\left (\vmini{-2} \right ) = u_+, & 
\uu{1}\left (\vmini{-1} \right ) = u_-, 
\end{array}
\]
\[
\begin{array}{llll}
\uu{2}\left (\vmini{1} \right )= u_0, & 
\uu{2}\left (\vmini{2} \right ) = u_+^3, &
\uu{2}\left (\vmini{3} \right ) = u_- u_+^2, &
\uu{2}\left (\vmini{4} \right ) = u_-^2 u_+, \vspace{2mm}\\
\uu{2}\left (\vmini{5} \right ) = u_+^2, &
\uu{2}\left (\vmini{6} \right ) = u_-^3, &
\uu{2}\left (\vmini{7} \right ) = u_- u_+, &
\uu{2}\left (\vmini{8} \right ) = u_0, \vspace{2mm}\\
\uu{2}\left (\vmini{-6} \right ) = u_+^3, &
\uu{2}\left (\vmini{-5} \right ) = u_-^2, &
\uu{2}\left (\vmini{-4} \right ) = u_- u_+^2, &
\uu{2}\left (\vmini{-3} \right ) = u_-^2 u_+, \vspace{2mm}\\
\uu{2}\left (\vmini{-2} \right ) = u_-^3, &
\uu{2}\left (\vmini{-1} \right ) = u_0. 
\end{array}
\]
\begin{equation}
\label{eqn:psig2}
\end{equation}
\subsubsection{Affine crystal}
Here we recall the definition of affine crystal.
Let \G ~ be an affine Lie algebra over $\Q$ with an indecomposable generalized
Cartan Matrix.
Let $c \in \sum_{i} \Z_{\geq 0}\h_i$ be  the canonical central element.
Let $\{ \alpha_i \mid i\in I \} \subset \h^*$ be the set of simple roots, 
and let $\{ h_i \mid i \in I \} \subset \h$ be the set of coroots. 
Let $\delta \in \sum_i \Z_{\geq 0}\alpha_i$ be the generator of null roots. 
Set $\h_{cl} = 
\bigoplus_{i \in I}\Q h_i \subset \h$ and $\h_{cl}^* = ( 
\bigoplus_{i \in I}\Q h_i)^*$. Let $cl: \h^* \rightarrow 
\h^*_{cl}$ denote the canonical morphism. We have an exact sequence
$0 \rightarrow \Q\delta \rightarrow \h^* \rightarrow \h^*_{cl}
\rightarrow 0$.
Fix $i_0 \in I$ and take an integer $d$ such that $\delta - d\alpha_{i_0} \in
\sum_{i \neq i_0}\Z \alpha_i$. For simplicity we write $0$ for $i_0$.
We define a map $af:\h^*_{cl} \rightarrow \h^*$ satisfying: $cl \circ
af = {\rm id}$ and $cl \circ af(\alpha_i) = \alpha_i$ for $i \neq 0$.
Let $\Lambda_i$ be the element of $\h^*_{cl} \subset \h^*$ such that 
$\langle h_j, \Lambda_i \rangle = \delta_{ij}$.
Hence we have $\alpha_i = \sum_j \langle h_j, \alpha_i \rangle af(\Lambda_i)
+ \delta_{i,0}d^{-1} \delta$. We define $P = \sum_i \Z af(\Lambda_i)+\Z d^{-1}
\delta \subset \h^*$ and $P_{cl} = cl(P) \subset \h_{cl}^*$. An element of $P$
is called an {\it affine weight} and an element in $P_{cl}$ is called a 
{\it classical weight}. Let $\Up\G$ be the quantized universal enveloping 
algebra associated with $P_{cl}$. A $P$-weighted crystal is called an 
{\it affine crystal} and a $P_{cl}$-weighted crystal is called a {\it classical
crystal}.
Let ${\rm Mod}^f(\G,P_{cl})$ be the category of $\Up{\G}$-module $M$
satisfing the following conditions:
\[
M \mbox{ has the weight decomposition } M = \oplus_{\lambda \in P_{cl}} 
M_\lambda, 
\]
and
\[
M \mbox{ is finite-dimensional over } \Q(q).
\]
For a $\Up{\G}$-module $M$ in ${\rm Mod}^f(\G,P_{cl})$, we define the 
$\U\G$-module ${\rm Aff}(M)$ by 
\[
{\rm Aff}(M) = \oplus_{\lambda \in P}{\rm Aff}(M)_\lambda,\;
{\rm Aff}(M)_\lambda = M_{cl(\lambda)} \mbox{ for } \lambda \in P.
\]
The actions of $e_i$ and $f_i$ are defined by the commutative diagrams
\[
\begin{array}{ccc}
{\rm Aff}(M)_\lambda&\stackrel{e_i}{\longrightarrow}&
{\rm Aff}(M)_{\lambda + \alpha_i}\\
\wr | && \wr |\\
M_{cl(\lambda)}&\stackrel{e_i}{\longrightarrow}&M_{cl(\lambda + \alpha_i)},
\end{array}
and
\begin{array}{ccc}
{\rm Aff}(M)_\lambda&\stackrel{f_i}{\longrightarrow}&
{\rm Aff}(M)_{\lambda - \alpha_i}\\
\wr | && \wr |\\
M_{cl(\lambda)}&\stackrel{f_i}{\longrightarrow}&M_{cl(\lambda - \alpha_i)}.
\end{array}
\]
We define the $\Up\G$-linear automorphism $T$ of ${\rm Aff}(M)$ by 
\[
\begin{array}{ccc}
{\rm Aff}(M)_\lambda&\stackrel{T}{\longrightarrow}&
{\rm Aff}(M)_{\lambda + \delta}\\
\wr | && \wr |\\
M_{cl(\lambda)}&\stackrel{id}{\longrightarrow}&M_{cl(\lambda + \delta )}.
\end{array}
\]

\subsubsection{Quantized universal enveloping algebra $\U{\gone}$-module}
Consider the finite  dimensional affine Lie algebra $\gone$ with Cartan
matrix 
\raisebox{3mm}{
$\left ( 
\begin{array}{ccc}
2 & -1 & 0\\
-1 & 2 & -1\\
0 & -3 & 2
\end{array}
\right )$}.\\
Let $\{\alpha_0, \alpha_1, \alpha_2\}$, $\{h_0, h_1, h_2\}$ be the set of simple root 
and simple coroot respectively. 
Let $\mbox{\frak h}$ and $\Lambda_i$ $(i=0,1,2)$ be Caltan subalgebra, 
fundermental weight, respectively. 
We define a non-degenerate symmetric bilinear form $( , )$ on $\mbox{\frak h}^*$
by $( \alpha_0, \alpha_0 ) = ( \alpha_1, \alpha_1 ) = 3$, 
$( \alpha_2, \alpha_2 ) = 1$, $( \alpha_0, \alpha_1 ) = ( \alpha_1, \alpha_0 )
 = ( \alpha_1, \alpha_2 ) = ( \alpha_2, \alpha_1 ) =  -\frac{3}{2}$,
$\langle h_i, \alpha_j \rangle$ is $(i,j)$-element of Cartan matrix and 
set $q_i = q^{(\alpha_i,\alpha_i)}$, $t_i = q^{h_i}$ for $i=1,2$.
Canonical central element is $c = h_0 + 2h_1 +h_2$.
Then the corresponding quantized universal enveloping algebra $\U{\gone}$ is the
associative algebra over $\Q(q)$ generated by $e_i$, $f_i$, $t_i$, 
$t_i^{-1}$ $(i = 0,1,2)$ satisfying the relations:
\[
t_it_j = t_jt_i, \hspace{5mm} t_i^{}t_i^{-1}= t_i^{-1}= t_i^{}=1,
\]
\[
t_i^{}e_j^{}t_i^{-1}=q_i^{\langle h_i, \alpha_j \rangle } e_j,
t_i^{}f_j^{}t_i^{-1}=q_i^{-\langle h_i, \alpha_j \rangle } f_j,
\]
\[
[e_i,f_j] = \delta_{ij}\frac{t_i^{}-t_i^{-1}}{q_i-q_i^{-1}}
\]
\[
\sum_{n=0}^l(-1)^ne_i^{(n)}e_je_i^{(l-n)} =
\sum_{n=0}^l(-1)^ne_i^{(n)}e_je_i^{(l-n)} = 0, \;\;
\left (i \neq j,\;l= 1- \langle h_i, \alpha_j \rangle \right ).
\]

%% file: s_perf.tex
\subsection{Perfect crystal[9]}
Let $B$ be a classical crystal. For $b \in B$, we set $\ve{}(b) = 
\sum_i\ve{i}(b) \Lambda_i$, and $\vp{}(b) = \sum_i\vp{i}(b) \Lambda_i$.
Note that $\wt(b) = \vp{}(b) - \ve{}(b)$. Set $P^+_{cl} = \{ \lambda \in P_{cl}
\:|\: \langle h_i, \lambda \rangle \geq 0 \mbox{ for all } i \in I \}$ and
for $l \in \Z_{\geq 0}$, let $(P^+_{cl})_l = \{ \lambda \in P^+_{cl} \:|\:
\langle c, \lambda \rangle = l \}$.
\begin{definition} 
\label{definition:perfectness}
\rm (Definition 4.6.1 in \cite{KMN2})\\ 
For $l \in \Z_{\geq 0}$, we say that $B$ is a {\it perfect crystal of level $l$}
if
\renewcommand{\labelenumi}{\rm (\arabic{enumi}) }
\begin{enumerate}
\item $B \otimes B$ is connected.
\item There exists $\lambda_0 \in P_{cl}$ such that $\wt(B) \in \lambda_0 +
\sum_{i \neq 0}\Z_{\leq 0}cl(\alpha_i)$ and that $\sharp(B_{\lambda_0}) = 1$.
\label{item:uniquehighest}
\item There is a finite dimensional $U_q^{\prime}(\G)$-module with a crystal
pseudo-base $(L,B_{ps})$ such that\\
$B \cong B_{ps} / \pm 1$.
\item For any $b \in B$, we have $\langle c, \ve{}(b) \rangle \geq l$.
\item The maps $\ve{}$, $\vp{}\::\:B_l=\{b \in B \:|\: \langle c, \ve{}(b)
\rangle = l \} \rightarrow (P^+_{cl})_l$ are bijective.
\end{enumerate}
The elements in $B_l$ are called {\it minimal elements}.
\end{definition}

\begin{definition}\label{def:polar}
Let $M$ be a $\U\G$-module. A symmetric bilinear form $( , )$ on $M$ is 
called a {\it prepolarization} of $M$ if it satisfies for $u,v \in M$:
\renewcommand{\labelenumi}{\rm (\roman{enumi}) }
\begin{enumerate}
\item $(q^hu,v) = (u,q^hv)$,
\item $(e_iu,v) = (u,q_i^{-1}t_i^{-1}f_iv)$,
\item $(f_iu,v) = (u,q_i^{-1}t_ie_iv)$.
\end{enumerate}
A prepolarization is called a {\it polarization} if it is positive definite.
\end{definition}

\begin{lemma} \label{lemma:Mirreduce}
{\rm (Lemma 3.4.4 in \cite{KMN2})}\\
Let $M$ be a $\Up\G$-module in $\mbox{Mod}^f(\G,P_{cl})$. Assume that $M$ has 
a crystal pseudo-base $(L,B)$ such that
\renewcommand{\labelenumi}{\rm (\roman{enumi}) }
\begin{enumerate}
\item there exists $\lambda \in P_{cl}$ such that $\sharp(B/\{\pm 1\})_\lambda
=1$,
\item $B/\{\pm 1\}$ is connected.
\end{enumerate}
Then $M$ is an irreducible $\Up\G$-module.
\end{lemma}
Let us introduce the subalgebras ${\bf A}_Z$ and ${\bf K}_Z$ of $\Q(q)$ as 
follows:
\[
{\bf A}_Z=\{f(q)/g(q)\mid f(q),\,g(q)\in \Z[q],\, g(0)=1\},
\]
\[
{\bf K}_Z = {\bf A}_Z[q^{-1}].
\]
\begin{proposition} {\rm (Proposition 2.6.2 in \cite{KMN1})}\\
Assume that \G~ is finite dimensional and let $M$ be a finite dimensional 
integrable $\U\G$-module of $M$. Let $( , )$ be a prepolarization on $M$, and
$M_{K_Z}$ a $\U\G$-submodule of $M$ such that $(M_{K_Z},M_{K_Z})
\subset K_Z$.
Let $\lambda_1, \ldots, \lambda_m \in P_+$ %
$(\lambda_k=\sum_ja_{kj}\Lambda_j,\:a_{kj}\in \Z_{\geq 0},\;k=1,\ldots , m)$, and
we assume the following conditions.
\begin{enumerate}
\renewcommand{\labelenumi}{\rm (\roman{enumi}) }
\item $\displaystyle \dim M_{\lambda_k} \leq \sum_{j=1}^m \dim 
V(\lambda_j)_{\lambda_k}$   for $k=1, \ldots, m.$
\item There exists $u_j \in (M_{K_Z})_{\lambda_j}$ $(j=1, \ldots, m)$ 
such that $(u_j, u_k) \in \delta_{jk} + qA$ and $(e_iu_j, e_iu_j) \in 
qq^{-2(1+ \langle h_i, \lambda_j \rangle )}A$ for all $i \in I$.
\end{enumerate}
Set $L = \{ U \in M \:|\: (u,u) \in A \}$ and set $B = \{ b \in M_{K_Z} \cap
L/M_{K_Z} \cap qL \:| (b,b)_0 = 1 \}$. Then we have 
\renewcommand{\labelenumi}{\rm (\roman{enumi}) }
\begin{enumerate}
\item $( , )$ is a polarization on $M$,
\item $M \cong  \oplus V(\lambda_j)$,
\item $( , )_0$ is positive definite, and $(L,B)$ is a crystal pseudo-base of 
$M$.
\end{enumerate}
\end{proposition}
\begin{lemma}  {\rm (Lemma 4.1.2 \cite{KMN1})} \label{lemma:polarization}
A polarization $(,)$ for Kac-Moody Lie algebra with symmetrizable 
Cartan matrix satisfies following relations.
\renewcommand{\labelenumi}{\rm (\roman{enumi}) }
\begin{enumerate} 
\item Let $b$ be a global base which satisfies $e_ib =0$ and $\langle h_i, 
\wt(b)\rangle  = 1$. Then $(b, b) = (f_ib, f_ib)$.
\item Let $b$ be a global base which satisfies $e_ib=0$ and $\langle h_i, 
\wt(b)\rangle = 2$. Then $(b, b) = (f_i^{(2)}b,f_i^{(2)}b) = 
q_i^{-1}[2]_i^{-1}(f_ib, f_ib)$. 
\item Let $b$ be a global base which satisfies $e_ib=0$ and $\langle h_i, 
\wt(b)\rangle = 3$. Then $(b, b) = (f_i^{(3)}b,f_i^{(3)}b) = 
q_i^{-2}[3]_i^{-1}(f_i^{(2)}b, f_i^{(2)}) = q_i^{-1}[2]_i^{-1}(f_ib, f_ib)$. 
\end{enumerate}
\end{lemma}
The existace of polarization for Kac-Moody Lie algebra with symmetrizable 
Cartan matrix is proved in Proposition 3.4.4 \cite{K2}.

\begin{proposition} 
\label{prop:prop345}
{\rm (Proposition 3.4.5 in \cite{KMN1})}\\
Let $m$ be a positive integer and assume the following conditions:
\renewcommand{\labelenumi}{\rm (\roman{enumi}) }
\begin{enumerate}
\item $\langle h_i, l\lambda_0 + j \alpha_{i_0} \rangle \geq 0$ for 
$i \neq i_0$ and $0 \leq j \leq m$,
\item $\dim(V_l)_{l\lambda_0 + k\alpha_{i_0}} \leq \sum_{j=0}^m 
\dim V(l\lambda_0 +j\alpha_{i_0})_{l\lambda_0+ k\alpha_{i_0}}$ for 
$0 \leq  k \leq m$, where $V(\lambda)$ is an irreducible 
$\U{\G_{I\setminus{i_0}}}$-module with highest weight $\lambda$,
\item  There exists $i_1 \in I$ such that $\{i \in I \:|\: \langle h_{i_0},
\alpha_i \rangle < 0\} = \{i_1\}$,
\item $-\langle h_{i_0}, l\lambda_0 - \alpha_{i_1} \rangle \geq 0$.
\end{enumerate}
Then we have 
\[
V_l \cong  \oplus_{j=0}^m V(l\lambda_0 + j\alpha_{i_0}) \mbox{ as a }
\U{\G_{I\setminus\{i_0\}}}\mbox{-module},
\]
and $V_l$ admits a crystal  pseudobase as a $U_q^{\prime}(\G)$-module.
\end{proposition}  

%% file: fusion.tex
\subsection{Construction of the polarization of $V^1$}
\label{section:polarization}
The space $V^1 \cong V(\Lambda_1) \oplus V(0)$ is endowed with a 
$\Up\gone$-module structure as follows.
\vspace{3mm}\\
\sayouf{0}{-6}{2}{1}, \sayouf{0}{-4}{3}{1}, \sayouf{0}{-3}{4}{1}, 
\sayouf{0}{-2}{6}{1}, \\
\sayouf{0}{-1}{9}{1}\sayounokori{1}{[2]_1}{8}, \sayouf{0}{9}{1}{[2]_0},\\
\sayouf{1}{1}{2}{1}, \sayouf{1}{4}{5}{1}, \sayouf{1}{6}{8}{1}
\sayounokori{1}{[2]_2}{7}\sayounokori{1}{[2]_1}{9},\\
\sayouf{1}{8}{-6}{[2]_2}, \sayouf{1}{-5}{-4}{1}, \sayouf{1}{-2}{-1}{1},\\
\sayouf{2}{2}{3}{1}, \sayouf{2}{3}{4}{[2]_2}, \sayouf{2}{4}{6}{[3]_2}, 
\sayouf{2}{5}{7}{1}\sayounokori{[3]_2}{[2]_1}{8}, \\
\sayouf{2}{7}{-5}{[2]_2},\sayouf{2}{-6}{-4}{1}, \sayouf{2}{-4}{-3}{[2]_2}, 
\sayouf{2}{-3}{-2}{[3]_2},\vspace{3mm}\\
\sayoue{0}{6}{-2}{1}, \sayoue{0}{4}{-3}{1}, \sayoue{0}{3}{-4}{1}, 
\sayoue{0}{2}{-6}{1}, \\
\sayoue{0}{1}{9}{1}\sayounokori{1}{[2]_1}{8}, \sayoue{0}{9}{-1}{[2]_0},\\
\sayoue{1}{-1}{-2}{1}, \sayoue{1}{-4}{-5}{1}, \sayoue{1}{-6}{8}{1} 
\sayounokori{1}{[2]_2}{7}\sayounokori{1}{[2]_1}{9},\\
\sayoue{1}{8}{6}{[2]_2}, \sayoue{1}{5}{4}{1}, \sayoue{1}{2}{1}{1}\\
\sayoue{2}{-2}{-3}{1}, \sayoue{2}{-3}{-4}{[2]_2}, \sayoue{2}{-4}{-6}{[3]_2}, 
\sayoue{2}{-5}{7}{1}\sayounokori{[3]_2}{[2]_1}{8}, \\
\sayoue{2}{7}{5}{[2]_2}, \sayoue{2}{6}{4}{1}, \sayoue{2}{4}{3}{[2]_2}, 
\sayoue{2}{3}{2}{[3]_2}, \vspace{2mm}\\
$q^{h_0}=q^{-2h_1-h_2}$,\\
where \vmini{9} ~is the base of $V(0)$.\\
We put
\[
B^1 = \{ \vminiD{i}, \vmini{7}, \vmini{8}, \vminiD{\overline{i}}, 
\vmini{9} \mid i=1,\ldots 6\}.
\]
 
Let $( , )_1$ be the polarization on the $\U{G_2}$-module $V(\Lambda_1)$.
We shall define a symmetric bilinear form $( , )$ on $V^1$ by
\[
\left ( \vmini{9}, u \right ) = \left ( u, \vmini{9} \right ) = 0 \; \; u \in
V(\Lambda_1),
\]
\[
\left ( \vmini{9}, \vmini{9} \right ) = q_0 [2]_0 \left ( \vmini{1}, \vmini{1}
\right )_1,
\]
\[
(u, v)=(u, v)_1 \;\; u,v \in V(\Lambda_1).
\]
By Lemma \ref{lemma:polarization}, we have 
\[
\left ( \vminiD{a_1}, \vminiD{a_1} \right ) = 
q_2^{-2} [3]_2^{-1} \left ( \vminiD{a_2}, \vminiD{a_2} \right ) =
q_2^{-3} [3]_2^{-1} [2]_2^{-1} \left ( \vmini{7},\vmini{7} \right ) =
q_1^{-1} [2]_1^{-1} \left ( \vmini{8}, \vmini{8} \right ),
\]
where $a_1 = 1, 2, 6, \overline{6}, \overline{2}, \overline{1}$, 
$a_2 = 3, 4, 5, \overline{5}, \overline{4}, \overline{3}$.
It follows that $( , )$ is a polarization on the $U_q^{\prime}(\gone)$-module
$V^1$. 

Therfore we see that $B^1$ is a crystal base and Kashiwara operators act as 
follows.
\vspace{3mm}\\
\sayouft{0}{-6}{2}{1}, \sayouft{0}{-4}{3}{1}, \sayouft{0}{-3}{4}{1}, 
\sayouft{0}{-2}{6}{1}, \sayouft{0}{-1}{9}{1}, \sayouft{0}{9}{1}{[2]_0},\\
\sayouft{1}{1}{2}{1}, \sayouft{1}{4}{5}{1}, \sayouft{1}{6}{8}{1}
\sayouft{1}{8}{-6}{[2]_2}, \sayouft{1}{-5}{-4}{1}, \sayouft{1}{-2}{-1}{1},\\
\sayouft{2}{2}{3}{1}, \sayouft{2}{3}{4}{[2]_2}, \sayouft{2}{4}{6}{[3]_2}, 
\sayouft{2}{5}{7}{1}, \sayouft{2}{7}{-5}{[2]_2}, \sayouft{2}{-6}{-4}{1},\\
\sayouft{2}{-4}{-3}{[2]_2}, \sayouft{2}{-3}{-2}{[3]_2},\vspace{3mm}\\
\sayouet{0}{6}{-2}{1}, \sayouet{0}{4}{-3}{1}, \sayouet{0}{3}{-4}{1}, 
\sayouet{0}{2}{-6}{1}, \sayouet{0}{1}{9}{1}, \sayouet{0}{9}{-1}{[2]_2},\\
\sayouet{1}{-1}{-2}{1}, \sayouet{1}{-4}{-5}{1}, \sayouet{1}{-6}{8}{1},
\sayouet{1}{8}{6}{[2]_2}, \sayouet{1}{5}{4}{1}, \sayouet{1}{2}{1}{1}\\
\sayouet{2}{-2}{-3}{1}, \sayouet{2}{-3}{-4}{[2]_2}, \sayouet{2}{-4}{-6}{[3]_2}, 
\sayouet{2}{-5}{7}{1}, \sayouet{2}{7}{5}{[2]_2}, \sayouet{2}{6}{4}{1}, \\
\sayouet{2}{4}{3}{[2]_2}, \sayouet{2}{3}{2}{[3]_2}, \vspace{2mm}\\
By checking of conditions of Definition \ref{definition:perfectness}, 
we have following proposition.
\begin{proposition}
\label{prop:level1B1}
$B^1$ is a perfect crystal of level $1$.
\end{proposition}

\subsection{Decomposition of the tensor product}
The tensor product $V^1 \otimes V^1$ is decomposed in the following way:
\begin{equation}
\label{eqn:decomp}
\begin{array}{ll}
(V(\Lambda_1) \oplus V(0)) \otimes (V(\Lambda_1) \oplus V(0))&\\
&\hspace{-30mm} \cong V(2\Lambda_1) \oplus V(3\Lambda_2) \oplus V(2\Lambda_2) 
\oplus V(\Lambda_1)^{\oplus 3} \oplus V(0)^{\oplus 2}.
\end{array}
\end{equation}
\begin{lemma}\label{lemma:hiest}\rm
Following vectors $u_{2\Lambda_1}$, $u_{3\Lambda_2}$, $u_{2\Lambda_2}$, 
$u_{\Lambda_1}^1$, $u_{\Lambda_1}^2$, $u_{\Lambda_1}^3$, $u_0^1$ and $u_0^2$ 
are the highest weight vectors with the weights $2\Lambda_1$, $3\Lambda_2$, 
$2\Lambda_2$, $\Lambda_1$, $\Lambda_1$, $\Lambda_1$, $0$ and $0$ respectively.
\begin{eqnarray*}
u_{2 \Lambda_1} & = & \vmini{1} \otimes \vmini{1},\\
u_{3 \Lambda_2} & = & \vmini{1} \otimes \vmini{2}
- q^3 \vmini{2} \otimes \vmini{1},\\
u_{2 \Lambda_2} & = & \vmini{1} \otimes \vmini{5} 
- q^3 \vmini{2} \otimes \vmini{4}
+ \frac{[ 2 ]_2}{[ 3 ]_2} q^4 \vmini{3} \otimes \vmini{3}
- q^7 \vmini{4} \otimes \vmini{2}
+ q^{10} \vmini{5} \otimes \vmini{1},\\
u^1_{\Lambda_1} & = & \vmini{1} \otimes \vmini{9},\\
u^2_{\Lambda_1} & = & \vmini{9} \otimes \vmini{1},\\
u^3_{\Lambda_1} & = & \vmini{1} \otimes \vmini{8}
- [ 2 ]_1 q^6 \vmini{2} \otimes  \vmini{6} 
+ \frac{[ 2 ]_1}{[ 3 ]_2} q^5 \vmini{3} \otimes \vmini{4}
- \frac{[ 2 ]_1}{[ 3 ]_2} q^6 \vmini{4} \otimes \vmini{3}\\
 & & + [ 2 ]_1 q^9 \vmini{6} \otimes \vmini{2}
- q^{12} \vmini{8} \otimes \vmini{1},\\
u^1_0 & = & \vmini{9} \otimes \vmini{9},\\
u^2_0 & = & \vmini{1} \otimes \vmini{-1}
- q^3 \vmini{2} \otimes  \vmini{-2} 
+ \frac{1}{[ 3 ]_2} q^2 \vmini{3} \otimes \vmini{-3}
- \frac{1}{[ 3 ]_2} q^3 \vmini{4} \otimes \vmini{-4}\\
&&+ \frac{1}{[ 3 ]_2} q^6 \vmini{5} \otimes \vmini{-5}
+ q^{6} \vmini{6} \otimes \vmini{-6}
- \frac{q^6}{[ 2 ]_2 [ 3 ]_2} \vmini{7} \otimes \vmini{7}
- \frac{q^6}{[ 2 ]_1} \vmini{8} \otimes \vmini{8}\\
&&+ \frac{q^8}{[ 3 ]_2} \vmini{-5} \otimes \vmini{5}
+ q^{12} \vmini{-6} \otimes \vmini{6}
- \frac{q^{11}}{[ 3 ]_2} \vmini{-4} \otimes \vmini{4}
+ \frac{q^{12}}{[ 3 ]_2} \vmini{-3} \otimes \vmini{3}\\
&&- q^{15} \vmini{-2} \otimes \vmini{2}
+ q^{18} \vmini{-1} \otimes \vmini{1}
- \frac{q^6}{[ 2 ]_1 [ 2 ]_3} \vmini{7} \otimes \vmini{8}
- \frac{q^6}{[ 2 ]_1 [ 2 ]_3} \vmini{8} \otimes \vmini{7}\\
&&- \frac{q^6}{[ 2 ]_1^2} \vmini{8} \otimes \vmini{9}
- \frac{q^6}{[ 2 ]_1^2} \vmini{9} \otimes \vmini{8}. 
\end{eqnarray*}
\end{lemma}

\subsection{Calculation of the $R$-matrix}
Let $V^1_x$ be $U_q^{\prime}(\G)$-module $\Q[x, x^{-1}] \otimes V^1$ with the 
actions of $e_i$, $f_i$, and $t_i$ given by $x^{\delta_{i0}}e_i$, 
$x^{-\delta_{i0}}f_i$, and $t_i$, respectively. 
The $R$-matrix for $V^1$ is an invertible ${\bf Q}[x, x^{-1},y,y^{-1}]\otimes
\U{G_2}$-linear map
\[
R(x,y): V_x\otimes V_y \rightarrow V_y \otimes V_x
\]
satisfying following properties.
\renewcommand{\labelenumi}{\rm (\arabic{enumi}) }
\begin{enumerate}
\item $R(x,y) \in \Q(q)[x/y,y/x]\otimes \mbox{End}_{Q(q)}(V^1\otimes V^1)$.
\item $(R(y,z)\otimes 1)(1 \otimes R(x,z))(R(x,y)\otimes 1)= 
(1\otimes R(x,y))(R(x,z)\otimes 1)(1 \otimes R(y,z))$.
\item $R(x,y)R(y,x)\in \Q(q)[x/y,y/x]$.
\end{enumerate}
By the decomposition \ref{eqn:decomp}, $\U{G_2}$-linearlity and
$R=R(x,y)$, we have
\begin{equation}
\begin{array}{l}
R(u_{2 \Lambda_1}) = a^{2 \Lambda_1} u_{2 \Lambda_1},
R(u_{3 \Lambda_2}) = a^{3 \Lambda_2} u_{3 \Lambda_2},
R(u_{2 \Lambda_2}) = a^{2 \Lambda_2} u_{2 \Lambda_2},\\
\displaystyle R(u^i_{\Lambda_1}) = \sum^3_{j=1} a^{\Lambda_1}_{ij} 
u^j_{\Lambda_1}$ $(i=1,2,3),
\displaystyle R(u^i_0) = \sum^2_{j=1}a^0_{ij} u^j_0$ $(i=1,2).
\end{array}
\label{eqn:Rmatsep}
\end{equation}
Consider the highest weight vectors given in Lemma \ref{lemma:hiest}.
Then we have in $V_x^1 \otimes V_y^1$:
\renewcommand{\labelenumi}{\rm \arabic{enumi}. }
\begin{enumerate}
\item $f_0^{(2)} f_1 f_2^{(3)} f_1 u^1_{\Lambda_1} = [ 2 ]_1 
q^{-3} x^{-1} y^{-1} u_{2 \Lambda_1}$,
\item $f_0^{(2)} f_1 f_2^{(3)} f_1 u^2_{\Lambda_1} = [ 2 ]_1 
q^{-3} x^{-1} y^{-1} u_{2 \Lambda_1}$,
\item $f_0^{(2)} f_1 f_2^{(3)} f_1 u^3_{\Lambda_1} = [ 2 ]_1 
(x - q^6 y) (x + q^{12} y) q^{-6} x^{-2} y^{-2} u_{2 \Lambda_1}$,
\item $e_1e^{(3)}_2e^{(2)}_1e^{(3)}_2e_0e^{(2)}e^{(3)}_2e_1e_0
u^1_{\Lambda_1} = [ 2 ]_1y(x+y)u_{2\Lambda_1}$,
\item $e_1e^{(3)}_2e^{(2)}_1e^{(3)}_2e_0e^{(2)}e^{(3)}_2e_1e_0
u^2_{\Lambda_1} = [ 2 ]_1q^{-6}x(x+y)u_{2\Lambda_1}$,
\item $e_1e^{(3)}_2e^{(2)}_1e^{(3)}_2e_0e^{(2)}e^{(3)}_2e_1e_0
u^3_{\Lambda_1} = [ 2 ]_2q^5(q^4+1)(x-q^6y)u_{2\Lambda_1}$,
\item $f_0u^1_{\Lambda_1} = [ 2 ]_1 q^{-6}y^{-1}u_{2\Lambda_1}$,
\item $f_0u^2_{\Lambda_1} = [ 2 ]_1 x^{-1}u_{2\Lambda_1}$,
\item $f_0u^3_{\Lambda_1} = 0$,
\item $f_0^{(2)}u^1_0 = [ 2 ]^2_1q^{-3}x^{-1}y^{-1}u_{2 \Lambda_1}$,
\item $f_0^{(2)}u^2_0 = q^{-12}(x^2+q^{30}y^2)x^{-2}y^{-2}u_{2 \Lambda_1}$,
\item $f^{(2)}_1f^{(6)}_2f^{(4)}_1f^{(6)}_2f^{(2)}_1f^{(2)}_0f_1f^{(2)}_2
f_1f_2f_0f_1f^{(3)}_2f_1f_0u^1_0 \\
 = [ 2 ]^2_1(q^4+q^2+1)q^{-8}(q^6y^2+x^2)xyu_{2\Lambda_1}$,
\item $f^{(2)}_1f^{(6)}_2f^{(4)}_1f^{(6)}_2f^{(2)}_1f^{(2)}_0f_1f^{(2)}_2
f_1f_2f_0f_1f^{(3)}_2f_1f_0u^2_0 =$\\
$[ 2 ]_2q^{-10}(q^4+q^2+1)((q^{22}+q^{18})y^2+(q^6+1)
(q^{16}-q^{14}+q^{12}-2q^{10}+q^8-2q^6+q^4-q^2+1)xy
+(q^4+1)x^2)xyu_{\Lambda_1}$,
\item $f^{(2)}_0f_1f^{(3)}_2f^{(2)}_1f^{(3)}_2u_{3\Lambda_2}
 = q^{-3}(x-q^6y)(x+y)x^{-2}y^{-2}u_{2\Lambda_1}$,
\item $f^{(2)}_0f_1f^{(3)}_2f^1f^{(2)}_2u_{2\Lambda_2}
 = q^{-8}(q^4+q^2+1)(x-q^6y)(x-q^{10}y)x^{-2}y^{-2}u_{2\Lambda_1}$.
\end{enumerate}
From these equations and the relation $[R,\Delta(X)]=0$ $(X\in \U{G_2})$, 
we obtain\\
$a_{2 \Lambda_1}(q^6y,\;x,\;(q^4+1)[ 2 ]_2 q^4 (x-q^6y))=
(q^6x,\;y,\;(q^4+1)[ 2 ]_2 q^4 (y-q^6x))\,{}^t\!( a^{\Lambda_1}_{ij} )$,
\vspace{3mm}\\
$a_{2 \Lambda_1}(x,\;y,\;0)=(y,\;x,\;0)\,{}^t\!( a^{\Lambda_1}_{ij} )$,
\vspace{3mm}\\
$a_{2 \Lambda_1}([ 2 ]^2_1 q^6,\;(x^2+q^{30}y^2)q^{-3}x^{-1}y^{-1}) =
([ 2 ]^2_1 q^6,\;(y^2+q^{30}x^2)q^{-3}x^{-1}y^{-1})\,{}^t\! ( a^0_{ij} )$,
\vspace{3mm}\\
$a_{2 \Lambda_1}(s_1(x,y),\;s_2(x,y))=(s_1(y,x),\;s_2(y,x))\,{}^t\!(a^0_{ij})$,
\vspace{3mm}\\
$(x-q^6y)a_{2 \Lambda_1} = a_{3 \Lambda_2}(y-q^6x)$,\vspace{3mm}\\
$(x-q^6y)(x-q^{10}y)a_{2 \Lambda_1} 
= a_{2 \Lambda_2}(y-q^6x)(y-q^{10}x)$,\\
where\\
$s_1(x,y)=[ 2 ]_1 (q^4-a^2+1)(x^2+q^6y^2)$\vspace{3mm}\\
$s_2(x,y)=(q^{22}+q^{18})y^2+(q^6+1)
(q^{16}-q^{14}+q^{12}-2q^{10}+q^8-2q^6+q^4-q^2+1)xy+(q^4+1)x^2$.\\
Let $P_{2\Lambda_1}$, $P_{3\Lambda_2}$, $P_{2\Lambda_2}$, $P_{u_{\Lambda_1}^i}$
$(i=1,2,3)$ and $P_{u_0^i}$ $(i=1,2)$ be the projection from $V^1 \otimes V^1$
to $V(2\Lambda_1)$, $V(3\Lambda_2)$, $V(2\Lambda_2)$, $\U{G_2} 
u^i_{\Lambda_1}$ $(i=1,2,3)$ and  $\U{G_2} u^i_0$ $(i=1,2)$ respectively. 
\begin{proposition} \label{prop:rmatrix} \rm 
Put $z = xy^{-1}$. Since the $R$-matrix $R(x,y)$ depends only on $x/y$, we denote 
$R(z)=R(x,y)$. Then we have
\begin{eqnarray*}
R(z)(u_{2\Lambda_1})&=&(1-q^{12}z)(1-q^{10}z)(1-q^8z)(1-q^6z)u_{2\Lambda_1},\\ 
R(z)(u_{3\Lambda_2})&=&(1-q^{12}z)(1-q^{10}z)(1-q^8z)(z-q^6)u_{3\Lambda_2},\\ 
R(z)(u_{2\Lambda_2})&=&(1-q^{12}z)(z-q^{10})(1-q^8z)(z-q^6)u_{2\Lambda_2},\\ 
R(u^i_{\Lambda_1}) &=& \sum^3_{j=1} a^{\Lambda_1}_{ij} 
u^j_{\Lambda_1}\;\;(i=1,2,3),\\
R(u^i_0) &=& \sum^2_{j=1}a^0_{ij} u^j_0\;\;(i=1,2).
\end{eqnarray*}
\end{proposition}  
Here $(a_{ij}^{\Lambda_1})$ and $(a_{ij}^0)$ are given by \vspace{2mm}\\
$ a^{\Lambda_1}_{11}=(1-q^{12}z)(q^6-1)(q^2+1)z
(-(q^{16}-q^{14}+q^{12}-q^{10}-q^6)z-(q^4-q^2+1))$,\vspace{2mm}\\
$ a^{\Lambda_1}_{12}=(1-q^{12}z)q^6(1-z)
(q^{12}z^2+(q^{12}-q^6-q^4-q^2)z+1)$,\vspace{2mm}\\
$ a^{\Lambda_1}_{13}=(1-q^{12}z)q^3(q^6-1)z(z-1)$,\vspace{2mm}\\
$ a^{\Lambda_1}_{21}=(1-q^{12}z)q^6(1-z)
(q^{12}z^2-(q^{10}+q^8+q^6-1)z+1)$,\vspace{2mm}\\
$ a^{\Lambda_1}_{22}=(1-q^{12}z)(q^6-1)(q^2+1)z((-q^{16}+q^{14}-q^{12})z
+(q^{10}+q^6-q^4+q^2-1))$,\vspace{2mm}\\
$ a^{\Lambda_1}_{23}=(1-q^{12}z)q^3(q^6-1)(z-1)z$,\vspace{2mm}\\
$ a^{\Lambda_1}_{31}=(1-q^{12}z)q^9(q^{12}-1)(q^4+1)(q^2+1)z(1-z)(z-q^6)$,
\vspace{2mm}\\
$ a^{\Lambda_1}_{32}=(1-q^{12}z)q^3(q^{12}-1)(q^4+1)(q^2+1)(1-z)(q^6-z)$,
\vspace{2mm}\\
$ a^{\Lambda_1}_{33}=(1-q^{12}z)(z-q^6)
(q^6z^2+(q^{18}-q^{12}-q^{10}-q^8-q^6+1)z+q^{12})$,\vspace{2mm}\\
$ a^0_{11}=q^{30}z^4-q^{24}(q^4+1)(q^2+1)z^3+
(q^{36}-q^{30}+q^{22}+q^{20}+2q^{18}+q^{16}+q^{14}-q^6+1)z^2
-q^6(q^4+1)(q^2+1)z+q^6$,\vspace{2mm}\\
$ a^0_{12}=-q^3(q^{12}-1)(q^6+1)(1-z)(1+z)z$,\vspace{2mm}\\
$ a^0_{21}=-q^3\frac{q^6-1}{q^4-q^2+1}(1-z)(1+z)((q^{22}+q^{18})+
(q^{40}-q^{38}+q^{36}-q^{34}-q^{30}-q^{26}-q^{20}-q^{14}-q^{10}-q^6
+q^4-q^2+1)z+(q^{22}+q^{18})z^2)$,\vspace{2mm}\\
$a^0_{22}=q^{30}-q^{24}(q^4+1)(q^2+1)z+
(q^{36}-q^{30}+q^{22}+q^{20}+2q^{18}+q^{16}+q^{14}-q^6+1)z^2
-q^6(q^4+1)(q^2+1)z^3+q^6z^4$.

\subsection{Fusion construction}
Let $l$ be a positive integer and $\mbox{\frak S}_l$ the symmetric group of 
order $l$.
Let $s_i$ be the simple reflection (the permutation of $i$ and $i+1$). 
Let $l(w)$ be the length of $w \in \mbox{\frak S}_l$.
Then for any $w \in \mbox{\frak S}_l$, we can define $R_w(x_1, \ldots, x_l):
V_{x_1} \otimes \cdots \otimes V_{x_l} \rightarrow V_{x_{w(1)}} \otimes \cdots
\otimes V_{x_{w(l)}}$ as follows:
\[
R_1(x_1, \ldots, x_l)=1.
\]
\[
R_{s_i}(x_1, \ldots, x_l) = \biggl ( \otimes_{j < i} {\rm id}_{V_{x_j}} 
\biggr ) \otimes R(x_i,x_{i+1}) \otimes \biggl ( \otimes_{j > i+1} 
{\rm id}_{V_{x_j}} \biggr ).
\]
For $w$, $w'$ with $l(ww') = l(w)+l(w')$,
\[
R_{ww'}(x_1, \ldots, x_l) = R_{w'}( x_{w(1)}, \ldots, x_{w(l)} ) \circ
R_w(x_1, \ldots, x_l).
\]
Fix $r \in \Z_{>0}$. For each $l \in \Z_{>0}$, we put 
\[
R_l = R_{w_0} \left ( q^{r(l-1)},q^{r(l-3)}, \ldots, q^{-r(l-1)} \right ):
\]
\[
V_{q^{r(l-1)}}\otimes V_{q^{r(l-3)}}\otimes \cdots \otimes V_{q^{-r(l-1)}}
\rightarrow V_{q^{-r(l-1)}}\otimes V_{q^{-r(l-3)}}\otimes \cdots 
\otimes V_{q^{r(l-1)}},
\]
where $w_0 \in \mbox{\frak S}_l$ is the permutation given by $i \mapsto l+1-i$.
Then $R_l$ is a $U_q^{\prime}(\G)$-linear homomorphism. 
We define 
\[
V_l = {\rm Im}R_l.
\]
Hence taking $\Lambda_1-2\Lambda_0$ as $\lambda_0$ in Definition 
\ref{definition:perfectness}. 
By Proposition \ref{prop:rmatrix}, we have
\[
\varphi(z) = (1-q^{12}z)(1-q^{10}z)(1-q^8z)(1-q^6z).
\]
Puting
\[
r=3.
\]
We see that
\[
\varphi(q^{2kr}) \mbox{ does not vanish for any } k>0.
\]
By Proposition \ref{prop:rmatrix}, we see
\begin{eqnarray*}
R(q^{2r})(u_{2\Lambda_1})&=&(1-q^{18})(1-q^{16})(1-q^{14})(1-q^{12}),\\
R(q^{2r})(u_{\Lambda_1}^1)&=&\sum_{j=1}^3 a_{1j}^{\Lambda_1}u_{\Lambda_1}^j,\\
R(q^{2r})(u_0^1)&=&\sum_{j=1}^2 a_{1j}^0u_0^j,\\
R(q^{2r})(u_{2\Lambda_2})&=&0,\\
R(q^{2r})(u_{\Lambda_1}^1-u_{\Lambda_1}^2)&=&0,\\
R(q^{2r})(u_{\Lambda_1}^3)&=&0,\\
R(q^{2r})(q^3(q^{12}-q^6+1)u_0^1-(q^6+1)u_0^2)&=&0.\makebox[5cm]{}
\end{eqnarray*}
Then $N = \ker \,R(q^{2r})$ contains $u_{3\Lambda_2}$, 
$u_{2\Lambda_2}$, $u^1_{\Lambda_1} - u^2_{\Lambda_1}$, $u^3_{\Lambda_1}$ and
$q^3(q^{12}-q^6+1)u^1_0-(q^6+1)u^2_0$. 
Therefore by (\ref{eqn:decomp}), we have 
\begin{eqnarray*}
N&\cong& \U{G_2}u_{3 \Lambda_2} \oplus \U{G_2} u_{2 \Lambda_2}\oplus 
\U{G_2}u^3_{\Lambda_1}\\
&&\oplus \U{G_2}(u^1_{\Lambda_1}-u^2_{\Lambda_1})
\oplus  \U{G_2}(q^3(q^{12}-q^6+1)u^1_0-(q^6+1)u^2_0).
\end{eqnarray*}
Then we have
\begin{eqnarray*}
\dim(V^{\otimes 2}/N)&=&\dim\left (\bigoplus^2_{j=0}V(j\Lambda_1)\right )\\ 
&=&\sharp \left \{\vminiD{b} \otimes \vminiD{a}\;\left |\;\vyokoD{a}{b}
\in B(2 \Lambda_1)\right .\right \}\\
&&+\sharp \left\{ \vminiD{a} \otimes \vmini{9} \; \left |\; 
\vminiD{a}\in B(\Lambda_1)\right . \right \}
+\sharp \left \{\vmini{9}\otimes \vmini{9} \right \}
\end{eqnarray*} 
where $a,b \in \{ 1, \ldots, 6, 0_1, 0_2, \overline{6}, \ldots, \overline{1}
\}$.\\
Hence,\\
\makebox[2cm]{} 
$\displaystyle \dim \left (V^{\otimes l}\left /\sum V^i
\otimes N \otimes V^{\otimes (l-2-i)} \right . \right )
$ \vspace{-2mm}
\begin{eqnarray*}
\makebox[15mm]{}&=& 
\sharp \left \{\vminiD{b_n}\otimes \cdots \otimes \vminiD{b_1} \otimes
\vmini{9} \otimes \cdots \otimes \vmini{9}\;\left |\;
\begin{array}{l}
\vyokoD{b_i}{b_{\!i\!+\!1}}\in B(2 \Lambda_1)\\
\vyokoD{\cdot}{b_1}\in B(\Lambda_1)
\end{array}
, \;0\leq i<n \leq l 
\right . \right \}\\
&=&\dim\left ( \bigoplus^l_{j=0}V(j\Lambda_1)\right ).
\end{eqnarray*} 
By (3.3.10) of \cite{KMN2},
$V_l$ is the quotient of 
$V^{\otimes l} \left /\sum^{l-2}_{i=0} V^{\otimes i}\otimes N
\otimes V^{(l-2-i)} \right .$. 

Therefore we have,
\[
\dim(V_l)_\lambda \leq \sum^{l}_{j=0}
\dim V(j(\Lambda_1-2\Lambda_0))_\lambda.
\]
We set $i_0 = 0$, then $\langle h_1,\lambda_0 \rangle =1$,  
$\langle h_2,\lambda_0 \rangle =0$, 
$\{i_1\}=\{1\}$, $\langle h_0,\lambda_0 \rangle =-2$.
Therefore applying Proposition \ref{prop:prop345}, we obtain 
the following results.
\begin{proposition}
\rm \mbox{} 
\begin{enumerate}
\item $V_l$  has a crystal pseudobase.
\item $V_l \cong \bigoplus_{j=0}^l (V(j(\Lambda_1-2 \Lambda_0)))$ \hspace{5mm}
 as a $\U{G_2}$-module. 
\end{enumerate}
\end{proposition}

%% file: prelim.tex
The algebra $\U{\gone}$ is $\Q(q)$-algebra generated by 
$\{e_i, f_i, t_i, t_i^{-1}\; (i \in I = \{0, 1, 2\}) \}$. 
Let $\U{(\gone)_J}$ be a $\Q(q)$-algebra generated by $\{e_i, f_i,t_i, t_i^{-1}
\;(i \in J \subset I)\}$.\\
The algebras $\U{(\gone)_{\{1,2\}}}$ and $\U{(\gone)_{\{0,1\}}}$ are 
isomorphic to $\U{G_2}$ and $\U{A_2}$ respectively. 
Let $J_0=\{1,2\}$ and $J_2=\{1,0\}$ be the index set of $A_2$ and 
$G_2$, respectively. We define $\imath_i:J_i \rightarrow I$ by 
$\imath_i(j)=j$ $(i=0,2)$.
In order to show that the crystal $B^l$ is perfect, we have to show the 
following conditions.
\renewcommand{\labelenumi}{\rm (\arabic{enumi}) }
\begin{enumerate}
\item There exists a crystal base of $\U{G_2}$ isomorphic to 
$\imath^*_0(B^l)$.\\
There exists a crystal base of $\U{A_2}$ isomorphic to $\imath^*_2(B^l)$.
Crystal base $B^l$ is connected.
\item For any $b \in B^l$, $\langle c, \vp{}(b) \rangle \geq l$.
\item The maps $\ve{}$ and $\vp{}$:$ B_l = \{ b \in B^l \mid \langle c, \ve{}(b)
\rangle = l \} \longrightarrow (P_{cl}^+)_l = \{ \lambda \in \sum \Z_{\geq 0}
\Lambda_i \mid \langle c, \lambda \rangle = 1 \}$ are bijective.
\item Crystal graph of $B^l \otimes B^l$ is connected.
\end{enumerate}
\renewcommand{\labelenumi}{\rm \arabic{enumi}. }
At first we construct crystal $B^l$ which satisfies condition (1)
(\S\ref{section:construction}, \S\ref{section:commutative}). Then we prove
conditions (2)(3) in \S\ref{section:minimal}, and (4) in \S\ref{section:connect}.

\subsection{Preparation}
\label{section:preparation}
Let $B^{A_2}(\lambda)$ (resp. $B^{G_2}(\lambda)$) be a crystal base of 
$U_q(A_2)$ (resp. $U_q(G_2)$) with the highest weight $\lambda$.
We write $\BB{\bB{b_1}{}\cdots \bB{b_n}{}}$ instead of 
Young tableau \crystalbase{b_1}{b_n}. We introduce following notations:
\begin{eqnarray*}
\bB{a}{k}&=&\overbrace{\;\bB{a}{}\cdots\:\bB{a}{}}^{k}\;\;(a=1,\ldots,6,0_1,
0_2,\overline{6},\ldots,\overline{1}),\vspace{2mm}\\
\bB{\CC}{k}&=& \tf{1}^k \bb{6}{k} \; = \;
\left \{
\begin{array}{ll}
\bb{6}{\frac{k}{2}}\bb{-6}{\frac{k}{2}}, & k \equiv 0 \pmod{2},
\vspace{2mm}\\
\bb{6}{\left [ \frac{k}{2} \right ]}\bb{8}{} 
\bb{-6}{\left [ \frac{k}{2} \right ]}, & k \equiv 1 \pmod{2},
\vspace{2mm}
\end{array}
\right . \vspace{2mm}\\
\bB{\CW}{k}&=& \tf{2}^{k}\bb{2}{k} \; = \;
\left \{
\begin{array}{ll}
\bb{2}{\frac{2k}{3}}\bb{6}{\frac{k}{3}}& k \equiv 0 \pmod{3},
\vspace{2mm}\\
\bb{2}{2\left [ \frac{k}{3} \right ]}\bb{3}{}
\bb{6}{\left [ \frac{k}{3} \right ]}& k \equiv 1 \pmod{3},
\vspace{2mm}\\
\bb{2}{2\left [ \frac{k}{3} \right ]+1}\bb{4}{}
\bb{6}{\left [ \frac{k}{3} \right ]}& k \equiv 2 \pmod{3},
\vspace{2mm}
\end{array}
\right . \vspace{2mm}\\
\bB{\CWb}{k}&=& \te{2}^{k}\bb{-2}{k} \; = \;
\left \{
\begin{array}{ll}
\bb{-6}{\frac{k}{3}}\bb{-2}{\frac{2k}{3}}& k \equiv 0 \pmod{3},
\vspace{2mm}\\
\bb{-6}{\left [ \frac{k}{3} \right ]}\bb{-3}{}
\bb{-2}{2\left [ \frac{k}{3} \right ]}& k \equiv 1 \pmod{3},
\vspace{2mm}\\
\bb{-6}{\left [ \frac{k}{3} \right ]}\bb{-4}{}
\bb{-2}{2\left [ \frac{k}{3} \right ]+1}& k \equiv 2 \pmod{3}.
\vspace{2mm}
\end{array}
\right . 
\end{eqnarray*}
Then we easily have
\begin{eqnarray}
\bB{\CC}{k} & = &\bb{6}{}\bB{\CC}{k-2}\bb{-6}{}, 
\label{eqn:hyouki1}\vspace{2mm}\\
\bB{\CW}{k} & = &\bb{2}{2}\bB{\CW}{k-3}\bb{6}{}, \vspace{2mm}\\
\bB{\CWb}{k} & = &\bb{-6}{}\bB{\CWb}{k-3}\bb{-2}{2}, 
\label{eqn:hyouki2}\vspace{2mm}\\
\tf{1}\bB{\CC}{k} & = &\bB{\CC}{k-1}\bb{-6}{}, \vspace{2mm}\\
\te{2}\bB{\CW}{k} & = &\bb{2}{}\bB{\CW}{k-1}, \vspace{2mm}\\
\te{2}\bB{\CWb}{k} & = &\bb{-6}{}\bB{\CWb}{k-2}\bb{2}{}. 
\label{eqn:hyouki3}\vspace{2mm}
\end{eqnarray}
\begin{proposition}
For $\bB{\CC}{k}$, $\bB{\CW}{k}$, $\bB{\CWb}{k}$ $(k\in {\bf Z}_{\geq 0})$, 
We have
\begin{eqnarray}
\mbox{\rm Red}_2\Bigl (\bB{\CC}{k}\Bigr )&=& u_0, 
\label{eqn:hyouki4}\vspace{2mm}\\
\mbox{\rm Red}_1\Bigl (\bB{\CW}{k}\Bigr )&=& u_0,
\label{eqn:hyouki5} \vspace{2mm}\\
\mbox{\rm Red}_1\Bigl (\bB{\CWb}{k}\Bigr )&=& u_0. 
\label{eqn:hyouki6}\vspace{2mm}
\end{eqnarray}
\end{proposition}
\proof We prove (\ref{eqn:hyouki4}).
By Definition \ref{def:uui} and (\ref{eqn:psig2}) we have
\[
\uu{2}\Bigl (\bB{\CC}{2k} \Bigr )= \uu{2}\Bigl (\bb{6}{k}\bb{-6}{k} \Bigr )=
u_+^{3k} u_-^{3k} ,
\]
\[
\mbox{\rm Red}_2 \Bigl (\bB{\CC}{2k}\Bigr )=u_0,
\]
\[
\uu{2}\Bigl(\bB{\CC}{2k+1}\Bigr )= 
\uu{2}\Bigl (\bb{6}{k}\bb{8}{}\bb{-6}{k} \Bigr ) = 
u_+^{3k} u_0 u_-^{3k}, 
\]
\[
\mbox{\rm Red}_2 \Bigl (\bB{\CC}{2k+1}\Bigr )= u_0.
\]
In a similar way, we can prove (\ref{eqn:hyouki5}), (\ref{eqn:hyouki6}).
\hfill$\Box$\\
We use the following notations for $b \in B^{A_2}(m\Lambda_1+n\Lambda_2)$,
\begin{eqnarray}
b &=& \younglambda \; \in B^{A_2}(m\Lambda_1 + n\Lambda_2)\nonumber\\
&=& \BB{\bB{b^{m+n}_1,b^{m+n}_2}{}\cdots\bB{b^{m+1}_1,b^{m+1}_2}{}
\bB{b^m_1}{}\cdots\bB{b^1_1}{}}.
\label{eqn:bbA2}
\end{eqnarray}
The highest weight element is given by  
\[
\youngHigh \;=\; \BB{\bB{1,2}{n}\bb{1}{m}}.
\] 

\begin{proposition}
Fix $l\in \Z_{>0}$. Take $i$, $j$, $p$ such that 
$0 \leq i \leq \left [ \frac{l}{2} \right ]$, 
$i\leq j \leq l-i$, $0 \leq p \leq j$.
Let $\overline{b}\in B^{A_2}((i+j)\Lambda_1+l\Lambda_2)$ be the highest weight 
element. Let $\underline{b}\in B^{A_2}((i+j)\Lambda_1+l\Lambda_2)$ be the 
lowest weight element. Then we have
\begin{eqnarray}
\ve{2}\left ( \tf{1}^{l+j}\tf{2}^l\overline{b} \right )&=&0, 
\label{eqn:A2char1}\\
\vp{2}\left ( \tf{1}^{l+j}\tf{2}^l\overline{b} \right )&=&j, 
\label{eqn:A2char1p}\\
\tf{2}^p\tf{1}^{l+j}\tf{2}^l\overline{b}&=&
\tf{1}^{l+j-p}\tf{2}^{l+p}\tf{1}^p\overline{b},
\label{eqn:A2char2}\\
\tf{1}^{l+j-p}\tf{2}^{l+p}\tf{1}^p\overline{b}&=&
\te{1}^i\te{2}^{l-p}\underline{b}.
\label{eqn:A2char3}
\end{eqnarray}
\end{proposition}
\proof
By Proposition \ref{prop:A2crystal}, we have 
\[
\overline{b}=\BB{\bB{1,2}{l}\bb{1}{i+j}},
\]
\[
\underline{b}=\BB{\bB{2,3}{l}\bb{3}{i+j}}.
\]
We prove (\ref{eqn:A2char1}), (\ref{eqn:A2char1p}).
Using Proposition \ref{proposition:vplus1}, we have
\begin{eqnarray*}
\uu{2} ( \overline{b} )&=& u_0^{i+j} ( u_0 u_+)^l,\\
\mbox{Red}_2 \left (\overline{b} \right )&=&u_+^l.
\end{eqnarray*}
Then we have 
\[
\tf{2}^l\overline{b}=\BB{\bB{1,3}{l}\bb{1}{i+j}}.
\]
Similarly, we have
\[
\tf{1}^{l+j}\tf{2}^l\overline{b}=\BB{\bB{1,3}{i}\bB{2,3}{l-i}\bb{2}{i+j}}.
\]
Using Proposition \ref{proposition:vplus1}, we have
\begin{eqnarray*}
\uu{2} \biggl ( \BB{\bB{1,3}{i}\bB{2,3}{l-i}\bb{2}{i+j}}\biggr )&=&
u_+^{i+j} (u_+ u_-)^{l-i} (u_0 u_-)^i, \\
\mbox{\rm Red}_2\biggl (\BB{\bB{1,3}{i}\bB{2,3}{l-i}\bb{2}{i+j}}\biggr ) &=&  
u_+^j.
\end{eqnarray*}
Therefore, we have
\[
\ve{2}\biggl (\BB{\bB{1,3}{i}\bB{2,3}{l-i}\bb{2}{i+j}}\biggr ) = 0,
\]
\[
\vp{2}\biggl (\BB{\bB{1,3}{i}\bB{2,3}{l-i}\bb{2}{i+j}}\biggr ) = j.
\]
Relations (\ref{eqn:A2char2}), (\ref{eqn:A2char3}) can be proved by using 
\[
\tf{1}^{l+j-p}\tf{2}^{l+p}\tf{1}^p\overline{b}=
\BB{\bB{1,3}{i}\bB{2,3}{l-i}\bb{2}{i+j-p}\bb{3}{p}}.
\]
\hfill$\Box$



\begin{proposition}\label{proposition:kakanA2a}
\rm Let $\overline{b}$ be the highest weight element of the crystal base 
$B^{A_2}(k\Lambda_1+j\Lambda_2)$ of $\U{A_2}$. Then for $0 \leq q < p \leq j$
we have,
\[
\tilde{f}_1 \tilde{f}_2 \left ( \tilde{f}_1^q \tilde{f}_2^p\overline{b} \right )=
\tilde{f}_2 \tilde{f}_1 \left ( \tilde{f}_1^q \tilde{f}_2^p\overline{b} \right ).
\]
\end{proposition}
\proof
By (\ref{eqn:bbA2}), we have
\[
\overline{b}=\BB{\bB{1,2}{j}\bb{1}{k}}.
\]
Using Proposition \ref{proposition:vplus1}, we have
\begin{eqnarray*}
\uu{2}(\overline{b})&=&u_0^k (u_0 u_+)^j,\\
\mbox{\rm Red}_2\left (\overline{b}\right )&=&u_+^j.
\end{eqnarray*}
Therefore we have
\[
\tf{2}^p\overline{b}=\BB{\bB{1,2}{j-p}\bB{1,3}{p}\bb{1}{k}}.
\]
Using Proposition \ref{proposition:vplus1} again, we have
\[
\uu{1}\left ( \tf{2}^p\overline{b}\right )=
 u_+^k (u_+ u_-)^p ( u_+ u_-)^{j-p}.
\]
Thus we have
\[
\tf{1}^q\tf{2}^p\overline{b} = 
\left \{
\begin{array}{ll}
\BB{\bB{1,2}{j-p}\bB{1,3}{p}\bb{1}{k-q}\bb{2}{q}} & (k > q), \vspace{2mm}\\
\BB{\bB{1,2}{j-p}\bB{1,3}{p-q+k}\bB{2,3}{q-k}\bb{2}{k}} & (k \leq q).
\end{array}
\right .
\]
In a similar way  we have,  
\[
\tf{2}\tf{1}\left ( \tf{1}^q\tf{2}^p\overline{b} \right ) = 
\tf{1}\tf{2}\left ( \tf{1}^q\tf{2}^p\overline{b} \right ) = 
\left \{
\begin{array}{ll}
\BB{\bB{1,2}{j-p}\bB{1,3}{p}\bb{1}{k-q-1}\bb{2}{q}\bb{3}{}} & (k\geq q),
 \vspace{2mm}\\
\BB{\bB{1,2}{j-p}\bB{1,3}{p-q+k-1}\bB{2,3}{q-k+1}\bb{2}{k-1}\bb{3}{}} & 
(k < q).
\end{array}
\right .
\]
In a similar way, we have the following Proposition.
\begin{proposition}\label{proposition:kakanA2b} 
Let $\overline{b}$ be the highest weight element of 
$B^{A_2}(k\Lambda_1+j\Lambda_2)$.
For $0 \leq p \leq j$, $0 \leq s \leq k$, we have 
\begin{equation}
\tilde{e}_1\tilde{f}_2^s\tilde{f}_1^{p+s+1}\tilde{f}_2^{p}\overline{b} =
\tilde{f}_2^s\tilde{f}_1^{p+s}\tilde{f}_2^p\overline{b} \neq
\tilde{f}_2^{s+1}\tilde{f}_1^{p+s}\tilde{f}_2^{p-1}\overline{b},
\label{eqn:kakanA2b1}
\end{equation} 
\begin{equation}
\tilde{e}_1\tilde{f}_2^{s+1}\tilde{f}_1^{p+s+1}\tilde{f}_2^{p}\overline{b} =
\tilde{f}_2^{s+2}\tilde{f}_1^{p+s}\tilde{f}_2^{p-1}\overline{b} \neq
\tilde{f}_2^{s+1}\tilde{f}_1^{p+s}\tilde{f}_2^p\overline{b}.
\label{eqn:kakanA2b2}
\end{equation} 
\end{proposition}
\begin{proposition}\label{proposition:phizero}
Let $\overline{b}$ be the highest weight element of crystal base 
$B^{A_2}(k\Lambda_1+j\Lambda_2)$. 
For an element $\tf{1}^q\tf{2}^p \overline{b}$ $(p \leq q \leq p+k)$,
we have
\begin{equation}
\varphi_2\left (\tilde{f}_1^q\tilde{f}_2^p b\right )=j+q-2p.
\end{equation}
\end{proposition}
\proof 
Using Proposition \ref{proposition:vplus1} we have,
\[
\mbox{\rm Red}_2\left (\tf{1}^q\tf{2}^p\overline{b} \right )= u_+^{j+p-2q}.
\]


\begin{proposition}\label{lemma:phizogen}\rm
For $b \in B^{G_2}(l\Lambda)$ such that $\tf{2}b \neq 0$, there exists unique $k' \in 
\{0, \ldots,  \vp{2}(b)-1\}$  such that 
\[
\vp{1}(\tf{2}^{(k+1)}(b))=\left \{
\begin{array}{ll}
\vp{1}(\tf{2}^{(k)}(b))&k\leq k',\\
\vp{1}(\tf{2}^{(k)}(b))+1 \;&k > k'.
\end{array}
\right .
\]
where $k = 0, \ldots, \vp{2}(b)-1$.
\end{proposition}
\proof
Put $b=b_1\otimes \cdots \otimes b_l$, 
and $r_1r_2\cdots r_{\ve{1}(b)+\vp{1}(b)}=\mbox{Red}_\fNs(b)$.
Let $k_1$ be $r_{\ve{1}(b)}$ if $\ve{1}(b) \neq 0$, $0$ if $\ve{1}(b)=0$.
By Proposition \ref{proposition:politensor}, there exists integers $k_2$ 
$\in \Z_{\geq 0}$ such that 
\[
\tf{2}b = b_1 \otimes \cdots \otimes \tf{2}b_{k_2} \otimes  \cdots \otimes b_l.
\]
First we consider $\uu{1}(\tf{2}b_{k_2})$, 
using $u_+$, $u_-$ in Defnition \ref{def:uui}.
\begin{equation}
\mbox{if }b_{k_2}=\bb{2}{}\!\!,\,\bb{5}{}\!\!,\,\bb{-6}{}\!\!,\,\bb{-4}{}\!\!,\; 
\mbox{ then }\uu{1}(b_{k_2})=u_-^n,\, \uu{1}(\tf{2}b_{k_2})=u_-^{n-1}\;(n=1,2).
\label{eqn:udecre}
\end{equation}
\begin{equation}
\mbox{if }b_{k_2}=\bb{3}{}\!\!,\,\bb{4}{}\!\!,\,\bb{7}{}\!\!,\,\bb{-3}{}\!\!,\; 
\mbox{ then }\uu{1}(b_{k_2})=u_+^{n'},\, \uu{1}(\tf{2}b_{k_2})=u_+^{n'+1}\;(n=0,1).
\label{eqn:uincre}
\end{equation}


We consider $\uu{1}(\tf{2}b)$. By Remark \ref{rem:uveuvp}, 
we have $\mbox{\rm Red}_1(b)= u_-^{\ve{1}(b)}u_+^{\vp{1}(b)}$. 
\vspace{3mm} 

\noindent The case of  $1 \leq k_2 \leq k_1$.\\
If (\ref{eqn:udecre}), we see that 
\begin{equation}
\mbox{\rm Red}_1(\tf{2}b)=u_-^{\ve{1}(b)-1}u_+^{\vp{1}(b)}.
\label{eqn:uminisone}
\end{equation}
If (\ref{eqn:uincre}), we see that the increased $u_+$ is paired with a neighboring 
$u_-$ and reduced. Then we have (\ref{eqn:uminisone}).\vspace{3mm}

\noindent The case of  $k_1+1 \leq k_2 \leq l$.\\
 If (\ref{eqn:udecre}), 
the decreased $u_-$ is paired with $u_+$ in $\uu{1}(b)$.
Then we see 
\begin{equation}
\mbox{\rm Red}_1(\tf{2}b)=u_-^{\ve{1}(b)}u_+^{\vp{1}(b)+1}.
\label{eqn:uplusone}
\end{equation}
If (\ref{eqn:uincre}), we see (\ref{eqn:uplusone}). 

Therefore, we have $\vp{1}(\tf{2}b)$ equals $\vp{1}(b)$ or $\vp{1}(b)+1$. 

Put $\tf{2}b = b' = b'_1 \otimes \cdots \otimes b'_l$. If $\tf{2}b' \neq 0$, we have
\[
\tf{2}b'=b'_1 \otimes \cdots \otimes \tf{2}b_{k'_2} \otimes \cdots \otimes b'_l,
\]
where $k'_2 \in \Z_{\geq 0}$.
Here the action of $\tf{2}$ is given by changing the leftmost $u_+$ to $u_-$ in 
$\mbox{Red}_2(b)$ and $\mbox{Red}_2(b') = u_-^{\ve{1}(b')}u_+^{\vp{1}(b')}$.
Then we see that $k'_2 \geq k_2$.

Thus we have the proposition.

%% file: const5.tex
\subsection{Construction of crystal base $B^l$ of $\U\gone$}
\label{section:construction}
\subsubsection{Structure of $B^l$}
As we noted earlier, $\U{(\gone)_{\{1,2\}}}$ is isomorphic to $\U{G_2}$ and
$\U{(\gone)_{\{1,0\}}}$ is isomorphic to $\U{A_2}$. 
Since the structure of $\U{A_2}$ is simpler than
that of $\U{G_2}$, we define the actions of $\tf{0}$ on $\U{G_2}$ by expoiting 
operators on $\U{A_2}$.
In view of isomorphism between $\U{(\gone)_{\{1,0\}}}$ and $\U{A_2}$,
we write $\{1,0\}$ as index set of roots of $\U{A_2}$.
The affine crystal $B^l$ for $\U\gone$ is constructed on $\CG^l$ with $\tfg{0}$.
\begin{definition}
Fix $l \in \Z_{>0}$. We define crystal base $\CG^l$ of $\U{G_2}$ and crystal
base $\CA$ of $\U{A_2}$ as follows.
\begin{eqnarray*}
\CG^l & = & \displaystyle\bigoplus_{n = 0}^{l} B^{G_2}(n \Lambda_1), 
\vspace{3mm}\\
{\cal A}_i & = & \displaystyle\bigoplus_{i \leq j_1,j_0 \leq l-i} B^{A_2},
(j_1 \Lambda_1 + j_0 \Lambda_0 ) \vspace{3mm}\\
{\cal A} & = & \displaystyle\bigoplus_{i = 0}^{\left [ \frac{l}{2} \right ]}
{\cal A}_i. 
\end{eqnarray*}
We denote $\tf{i}$, $\te{i}$, $\vp{i}$, $\ve{i}$, $\wt$, $\wt_i$ $(i=1,0)$ on 
$\CA$ by $\tfa{i}$, $\tea{i}$, $\vp{i}^\CA$, $\ve{i}^\CA$, $\wt^\CA$, 
$\wt_i^\CA$ respectively.
In a similar way, we denote $\tf{i}$, $\te{i}$, $\vp{i}$, $\ve{i}$, $\wt$, 
$\wt_i$ $(i=1,2)$ on $\CG$ by $\tfg{i}$, $\teg{i}$, $\vp{i}^\CG$, $\ve{i}^\CG$, 
$\wt^\CG$, $\wt_i^\CG$ respectively.
\end{definition}
\begin{proposition}\label{proposition:dimension}
\[
\sharp \CG^l = \sharp \CA 
\]
\end{proposition}
\proof The proposition amounts to show
\begin{equation}
\label{eqn:dimension}
\sum_{n=0}^l \sharp B^{G_2}(n\Lambda_1) =
\sum_{i=0}^{[\frac{l}{2}]}\sum_{i\leq j_1,j_2 \leq l-i}
\sharp B^{A_2}(j_1\Lambda_1+j_2\Lambda_2).
\end{equation}
Let $\Delta^+$ be the set of positive roots, put 
$\delta = \frac{1}{2}\sum_{\beta \in \Delta^+}\beta$,
and let $d(\Lambda)$ be the dimension of finite dimensinal irreducible 
representation with highest weight $\Lambda$.
Then we have 
\[
d(\Lambda) = \mathop{\mbox{\Large$\Pi$}}_{\beta \in \Delta^+}
\frac{(\Lambda+\delta, \beta)}{(\delta, \beta)}
\]
by \cite{W}. Therefore, we have 
\[
\sharp B^{G_2}(n\Lambda_1) = (n+1)\left (\frac{1}{2}n+1 \right )
\left (\frac{2}{3}n+1 \right )\left (\frac{3}{4}n+1 \right )
\left (\frac{3}{5}n+1 \right ),
\]
\[
\sharp B^{A_2}(n_1\Lambda_1+n_2\Lambda_2) = 
\frac{1}{2}(n_1+1)(n_2+1)(n_1+n_2+2).
\]
By calculation, we see that 
\begin{equation}
\label{eqn:dimGl1Amul}
\sharp B^{G_2}(l\Lambda_1) = \sum^{\left [ \frac{l}{2} \right ]}_{i=0} 
\left ( \sum^{l-i-1}_{k=i} \sharp B^{A_2}(k\Lambda_1+(l-i)\Lambda_2) + 
\sum^{l-i}_{j=i}\sharp B^{A_2}((l-i)\Lambda_1+j\Lambda_0)\right ).
\end{equation}
Thus, we have (\ref{eqn:dimension}). \hfill$\Box$

By Proposition \ref{proposition:dimension}, 
we can construct one-to-one correspondence between $\CG^l$ and $\CA$ as sets.

Let $\hwvl{i}{j_1}{j_0}\in \CA_i$ be the heighest weight element with 
weight $j_1\Lambda_1+j_0\Lambda_0$. 
Let $\lwv{l}{i}{(-j_1,-j_0)} \in \CA_i$ be the lowest weight 
element with weight $-j_1\Lambda_1-j_0\Lambda_0$. 
Elements $\hwvl{i}{j_1}{j_0}$ and $\lwv{l}{i}{(-j_0,-j_1)}$ 
satisfy 
\[
\hwvl{i}{j_1}{j_0}= 
\ea{0}{j_0}\ea{1}{j_1+j_0}\ea{0}{j_1}\lwv{l}{i}{(-j_0,-j_1)}.
\]

For a Lie algebra of type $A_2$, it is known that 
all elements $b$ in the crystal base $B^{A_2}(j_i\Lambda_1+j_0\Lambda_0)$ 
are uniquely expressed as
\begin{equation}   
b= \fa{0}{r}\fa{1}{q}\fa{0}{p}\hwvl{i}{j_1}{j_0} \;
(0 \leq p \leq j_0, p \leq q \leq p +j_1, 0 \leq r \leq j_0 +q -2p),
\label{eqn:Atwohyouki}
\end{equation}
It is possible that there exists crystal bases with same highest weight in 
$\CA_i$ and $\CA_{i'}$  $(i \neq i')$.
To distinguish them,  we define crystal base $B^i_{(j_1,j_0)} \subset \CA_i$ by
\begin{equation}
B^i_{(j_1,j_0)}=\left \{ \left .
\fa{0}{r}\fa{1}{q}\fa{0}{p}\hwvl{i}{j_1}{j_0}\;\right | \;
0\leq p\leq j_0,\;p\leq q\leq p+j_1,\;0\leq r \leq j_0+q-2p  \right \}.
\end{equation}

By calculation, we have the following proposition.
\begin{proposition}
\label{prop:reductBkj}
For $B^i_{(k,j)}$ $(i < k,j \leq l-i)$, we have
\begin{equation}
\begin{array}{l}
\displaystyle
B^i_{(k,j)} \biggm \backslash \left (
\left \{ \fa{1}{q}\fa{0}{p}\hwvl{i}{k}{j} \left | {\scriptstyle i< k,j 
\leq l-i \atop \scriptstyle 0\leq p \leq j,0 \leq q \leq p+k}  \right . 
\right \} \sqcup
\left \{ \ea{1}{q}\ea{0}{p}\lwv{l}{i}{(-j,-k)} \left | 
{\scriptstyle i< k,j \leq l-i \atop \scriptstyle 0\leq p \leq j,0 \leq q \leq 
p+k} \right . \right \}
\right )\vspace{2mm}\\
\displaystyle \makebox[3mm]{}\cong B^{A_2}((k-1)\Lambda_1+(j-1)\Lambda_0).
\end{array}
\label{eqn:reductBkj}
\end{equation}
\end{proposition}
For $l\in {\bf Z}_{>0}$, we define 
\begin{eqnarray*}
\CA^{(l)}_i &=& \left ( \bigsqcup_{i \leq k \leq l-i}
 B^i_{(k,i)} \right )  \bigsqcup
\left ( \bigsqcup_{i < j \leq l-q}  B^i_{(i,j)} \right )\\
&&\bigsqcup\left ( \bigsqcup_{\scriptstyle i< k,j \leq l-i \atop
0\leq p \leq j,0 \leq q \leq p+k} 
\fa{1}{q}\fa{0}{p}\hwvl{i}{k}{j} \right )
\bigsqcup\left ( \bigsqcup_{\scriptstyle i< k,j \leq l-i \atop
0\leq p < k ,0 \leq q \leq p+j} 
\ea{1}{q}\ea{0}{p}\lwv{l}{i}{(-j,-k)}\right ),\\
\CA^{(l)} &=& \bigsqcup_{0 \leq i \leq \left [ \frac{l}{2} \right ]} 
\CA^{(l)}_i.
\end{eqnarray*}
\begin{proposition}
\label{prop:lLambda}
For $l\in {\bf Z}_{>0}$, we have
\[
\sharp B^{G_2}(l\Lambda_1) = \sharp \CA^{(l)}.
\]
\end{proposition}

\proof 
By Proposition \ref{prop:reductBkj}, we have 
\begin{eqnarray*}
\sharp B^{A_2}\left ((l-i)\Lambda_1+(l-k)\Lambda_0 \right ) &=&
\sharp B^i_{(k,i)}\\ 
&&+\sum^{l-i-k}_{m=1}\sharp\left \{ \fa{1}{q}\fa{0}{p}\hwvl{i}{k+m}{i+m} 
\left | {0 \leq p \leq i+m \atop 0 \leq q \leq p+k+m} \right . \right \}\\
&&+\sum^{l-i-k}_{m=1} \sharp \left \{ \ea{1}{q}\ea{0}{p} 
\lwv{l}{i}{(-i-m,-k-m)} 
\left | {0 \leq p < k+m \atop 0 \leq q \leq p+i+m} \right . \right \},\\
\sharp B^{A_2}\left ((l-j)\Lambda_1+(l-i)\Lambda_0 \right ) &=&
\sharp B^i_{(i,j)} \\ 
&&+\sum^{l-i-j}_{m=1} 
\sharp \left \{ \fa{1}{q}\fa{0}{p} \hwvl{i}{i+m}{j+m} 
\left | {0 \leq p \leq j+m \atop 0 \leq q \leq p+i+m} \right . \right \}\\
&&+\sum^{l-i-j}_{m=1} \sharp \left \{ \ea{1}{q}\ea{0}{p} 
\lwv{l}{i}{(-j-m,-i-m)} 
\left | {0 \leq p < i+m \atop 0 \leq q \leq p+j+m} \right . \right \}.
\end{eqnarray*}
By (\ref{eqn:dimGl1Amul}), we have
\[
\sum^{\left [ \frac{l}{2}\right ]}_{i=0} \left (
\sum^{l-i}_{k=i}\sharp B^{A_2}((l-k)\Lambda_1+(l-i)\Lambda_0)+
\sum^{l-i}_{j=i+1}\sharp B^{A_2}((l-i)\Lambda_1+(l-j)\Lambda_0) \right )
= \sharp B^{G_2}(l\Lambda_1).
\]
\hfill$\Box$

By Proposition \ref{prop:reductBkj} and Proposition \ref{prop:lLambda}, 
we have the following lemma.
\begin{lemma}
\label{lemma:subset}
For $l\in {\bf Z}_{>0}$ we have
\[
\CG^l \setminus \CA^{(l)} \cong \CG^{l-1}.
\]
as crystal base of $\U{A_2}$.
\end{lemma}

\subsubsection{Operators $E_\CA$ and $F_\CA$ on $\CA$}
\label{sec:opeEF}
We are going to define operators $E_\CA$ and $F_\CA$ on $\CA$ which satisfy
properties 
\renewcommand{\labelenumi}{\rm (C\arabic{enumi}) }
\begin{enumerate}
\item for $b,\;b' \in \CA$, $E_\CA b=b'$, if and only if $F_\CA b'=b$,
\item for $b\in \CA$, $E_\CA \tfa{0}(b) = \tfa{0}E_\CA(b)$,
\item for $b \in \CA$, $\mbox{max}\{m \mid F_\CA^mb \neq 0\}-
\mbox{max}\{m' \mid E_\CA^{m'}b \neq 0 \}= -2\wt_1^\CA(b) -\wt_0^\CA(b)$.
\end{enumerate}
Operators $E_\CA$ and $F_\CA$ are counterpart of $\teg{2}$ and $\tfg{2}$,
respectively.
In view of Lemma \ref{eqn:reductBkj} and (C2) we define $E_\CA$ and $F_\CA$ 
inductively. 
If $l=0$, $\CA$ is trivial. 
If $l>0$, we use following relation for induction
\begin{equation}
E_\CA \left ( \fa{1}{q}\fa{0}{p}\hwvl{i}{k}{j} \right ) =
\tea{0}\left ( E_\CA \left (\fa{1}{q-1}\fa{0}p\hwv{l-1}{i}{k-1}{j-1}
\right ) \right ). 
\label{eqn:indsayou}
\end{equation}
Put $\CA^{(l)}_+ = \{ b \mid \te{0}(b)=0 \}$. 
Using Proposition \ref{proposition:vplus1}, we see that $\CA_+^{(l)}$ is given 
by
\[
\CA^{(l)}_+=\left \{ \fa{1}{q}\fa{0}{p}\hwvl{i}{k}{j}\left | 
0\leq i\leq \left [ \frac{l}{2} \right ],\: i\leq k,j \leq l-i,\: 
0\leq p \leq j,\: p\leq q \leq p+k \right . \right \}.
\]
By {\rm (C1)} and {\rm (C2)}, it is sufficient to define 
the operator $E_\CA$ for $b\in \CA_+^{(l)}$.
\begin{definition}
\label{definition:itinokettei}
We define $E_\CA$ by
\[
E_\CA \left ( \fa{1}{q}\fa{0}{p}\hwvl{i}{k}{j} \right ) =
\left \{
\begin{array}{ll}
\fa{1}{q}\fa{0}{p}\hwvl{i}{k-1}{j},&\mbox{\saba 1}\vspace{2mm}\\
\fa{1}{q+1}\fa{0}{p+1}\hwvl{i}{k}{j+1},&\mbox{\saba 2}\vspace{2mm}\\
\fa{1}{q+1}\fa{0}{p}\hwvl{i+1}{k+1}{j-1},&\mbox{\saba 3}\vspace{2mm}\\
\fa{1}{q+1}\fa{0}{p}\hwvl{i}{k+1}{j-1},&\mbox{\saba 4}\vspace{2mm}\\
\fa{1}{q+1}\fa{0}{p}\hwvl{i-1}{k+1}{j-1},&\mbox{\saba 5}\vspace{2mm}\\
0,&\mbox{\saba 6}\vspace{2mm}\\
\tea{0}\left ( E_\CA \left (\fa{1}{q-1}\fa{0}p\hwv{l-1}{i}{k-1}{j-1}
\right ) \right ). &\mbox{\rm (\ref{eqn:indsayou})}
\end{array}
\right . 
\]
Conditions for \saba 1-\saba 6 and {\rm (\ref{eqn:indsayou})} are as 
follows:
\renewcommand{\labelenumi}{\rm (\arabic{enumi}) }
\renewcommand{\labelenumii}{\rm (\roman{enumii}) }
\begin{enumerate}
\item If $i < k \leq l-i$, $i < j \leq l-i$,
\begin{enumerate} 
\item if $0\leq q \leq j-1-\left [ \frac{j-k}{3} \right ]$, 
$p=\mbox{min}(q,j)$, then the action of $E_\CA$ is given by \saba 1,
\item if $j<l-i$, $j-\left [ \frac{j-k}{3} \right ] \leq q \leq j+k$, 
$p=\mbox{min}(q,j)$, then the action of $E_\CA$ is given by \saba 2,  
\item if $j=l-i$, $j-\left [ \frac{j-k}{3} \right ] \leq q \leq j+k$, 
$p=\mbox{min}(q,j)$, then the action of $E_\CA$ is given by \saba 6,
\item if $0\leq p \leq j-1$, $p+1 \leq q \leq p+k$, then we use induction 
\mbox{\rm (\ref{eqn:indsayou})}.
\end{enumerate}
\item If $k=i$, $i+2 \leq j \leq l-i$,
\begin{enumerate}
\item if $0\leq p \leq j-1-\left [ \frac{j-i}{3} \right ]$,
$p \leq q \leq p+i$ then the action of $E_\CA$ is given by \saba 3,
\item if $j<l-i$,$j-\left [ \frac{j-i}{3} \right ] \leq p \leq j$, 
$p\leq q \leq p+i$ then the action of $E_\CA$ is given by \saba 2,
\item if $j=l-i$, 
$j-\left [ \frac{j-i}{3} \right ] \leq p \leq j$, $p\leq q \leq p+i$ 
then the action of $E_\CA$ is given by \saba 6.
\end{enumerate}
\item If $k=i$, $j=i+1$,
\begin{enumerate}
\item if $0 \leq p \leq i$, $p \leq q \leq p+i$ 
then the action of $E_\CA$ is given by \saba 4,
\item if $p=i+1$, $p \leq q \leq p+i$ 
then the action of $E_\CA$ is given by \saba 2.
\end{enumerate}
\item If $k=i$, $j=i$,
\begin{enumerate}
\item if $0\leq p \leq i-1$, $p \leq q \leq p+i$ then 
the action of $E_\CA$ is given by \saba 5,
\item if $p=i$, $i \leq q \leq 2i$ 
then the action of $E_\CA$ is given by \saba 2.
\end{enumerate}
\item If $i+1 \leq k \leq l-i$, $j=i$,
\begin{enumerate}
\item if $0 \leq p \leq i$, 
$p\leq q \leq p-1-\left [ \frac{i-k}{3} \right ]$ 
then the action of $E_\CA$ is given by \saba 1,
\item if $p=i$, $p-\left [\frac{i-k}{3}\right ]\leq q \leq p+k$ 
then the action of $E_\CA$ is given by \saba 2,
\item if $0 \leq p\leq i-1$, $p-\left [ \frac{i-k}{3} \right ] \leq q \leq p+k$
then the action of $E_\CA$ is given by \saba 5.
\end{enumerate}
\end{enumerate}
\end{definition}
By (C1) and definition of $E_\CA$, we see the action of $F_\CA$ on 
$\CA_+^{(l)}$ is given by 
\[
F_\CA \left ( \fa{1}{q}\fa{0}{p}\hwvl{i}{k}{j} \right ) =
\left \{
\begin{array}{ll}
\tf{1}^q\tf{0}^p\hwvl{i}{k+1}{j},&\mbox{\saba 1}'\vspace{2mm}\\
\tf{1}^{q-1}\tf{0}^{p-1}\hwvl{i}{k}{j-1},&\mbox{\saba 2}'\vspace{2mm}\\
\tf{1}^{q-1}\tf{0}^p\hwvl{i-1}{k-1}{j+1},&\mbox{\saba 3}'\vspace{2mm}\\
\tf{1}^{q-1}\tf{0}^p\hwvl{i}{k-1}{j+1},&\mbox{\saba 4}'\vspace{2mm}\\
\tf{1}^{q-1}\tf{0}^p\hwvl{i+1}{k-1}{j+1},&\mbox{\saba 5}'\vspace{2mm}\\
0,&\mbox{\saba 6}'\vspace{2mm}\\
\tea{0}\left ( F_\CA \left (\fa{1}{q-1}\fa{0}p\hwv{l-1}{i}{k-1}{j-1}
\right ) \right ).&(\ref{eqn:indsayou})'
\end{array}
\right . 
\]
As in the case of $E_\CA$, we use induction (\ref{eqn:indsayou})$'$. 
Conditions for \saba 1$'$-\saba 6$'$ are as follows:
\renewcommand{\labelenumi}{\rm (\arabic{enumi}) }
\renewcommand{\labelenumii}{\rm (\roman{enumii}) }
\begin{enumerate}
\item If $i < k \leq l-i$, $i < j \leq l-i$,
\begin{enumerate} 
\item if $k<l-i$, $0\leq q \leq j-1-\left [ \frac{j-k-1}{3} \right ]$, 
$p=\mbox{min}(q,j)$, then the action of $F_\CA$ is given by \saba 1$'$,
\item if $k=l-i$, $0\leq q \leq j-1-\left [ \frac{j-k-1}{3} \right ]$, 
$p=\mbox{min}(q,j)$, then  the action of $F_\CA$ is given by \saba 6$'$,
\item if $j-\left [ \frac{j-k-1}{3} \right ] \leq q \leq j+k$, 
$p=\mbox{min}(q,j)$, then the action of $F_\CA$ is given by \saba 2$'$,  
\item if $0\leq p \leq j-1$, $p+1 \leq q \leq p+k$, then we use induction 
\mbox{\rm (\ref{eqn:indsayou})$'$}.
\end{enumerate}
\item If $k=i$, $i+1 \leq j \leq l-i$,
\begin{enumerate}
\item if $0\leq p \leq j-1-\left [\frac{j-k-1}{3}\right ]$, $q=p$ then 
the action of $F_\CA$ is given by \saba 1$'$,
\item if $0\leq p \leq j-1-\left [\frac{j-k-1}{3}\right ]$, 
$p+1\leq q \leq p+k+1$ then the action of $F_\CA$ is given by \saba 3$'$,
\item if $j-\left [ \frac{j-k-1}{3} \right ] \leq p \leq j$, 
$p\leq q \leq p+k$, then the action of $F_\CA$ is given by \saba 2$'$.
\end{enumerate}
\item If $k=i$, $j=i$,
\begin{enumerate}
\item if $0\leq p \leq i$, $p=q$ 
then the action of $F_\CA$ is given by \saba 2$'$,
\item if $0\leq p \leq i$, $p+1 \leq q \leq p+i$ 
then the action of $F_\CA$ is given by \saba 3$'$.
\end{enumerate}
\item If $k=i+1$, $j=i$,
\begin{enumerate}
\item if $0\leq p \leq i$, $p=q$ 
then the action of $F_\CA$ is given by \saba 2$'$,
\item if $0\leq p \leq i$, $p+1 \leq q \leq p+i+1$ 
then the action of $F_\CA$ is given by \saba 4$'$.
\end{enumerate}
\item If $i+2 \leq k \leq l-i$, $j=i$,
\begin{enumerate}
\item if $0 \leq p \leq i$, 
$p\leq q \leq p-1-\left [ \frac{i-k-1}{3} \right ]$ 
then the action of $F_\CA$ is given by \saba 1$'$,
\item if $0 \leq p \leq i$, 
$p-\left [ \frac{i-k-1}{3} \right ] \leq q \leq p+k$ 
then the action of $F_\CA$ is given by \saba 5$'$.
\end{enumerate}
\end{enumerate}
For $b\in \CA$, we put 
\[
\ve{\CA}(b)=\mbox{max}\{m \mid F_\CA^mb \neq 0\},
\]
\[
\vp{\CA}(b)=\mbox{max}\{m' \mid E_\CA^{m'}b \neq 0 \}.
\]
We verify that $E_\CA$ is well-defined. 
By the definition of $F_\CA$, we see that for $b, b' \in \CA^{(l)}_+$,
\begin{equation}
\mbox{ if } E_\CA(b)= E_\CA(b') \neq 0, \mbox{ then } b=b'.
\label{eqn:injectEA}
\end{equation}
In order to verify (C2), we prove that 
$E_\CA(\tfa{0}(b))\neq 0$ if and only if $\tfa{0}(E_\CA(b))\neq 0$.
For this, it is sufficient to prove that for $b\in A_+^{(l)}$, 
\begin{equation}
\vp{0}(b)=\vp{0}(E_\CA b).
\label{eqn:phi0EAb}
\end{equation} 
Using Proposition \ref{proposition:phizero}, we have following formula in the 
case of \saba 1 of Definition \ref{definition:itinokettei},
\[
\vp{0}\left ( \tf{1}^q\tf{0}^p\hwvl{i}{k}{j} \right )=j+q-2p,
\]
\[
\vp{0}\left ( E_\CA\tf{1}^q\tf{0}^p\hwvl{i}{k}{j} \right )=
\vp{0}\left ( \tf{1}^q\tf{0}^p\hwvl{i}{k-1}{j} \right )=j+q-2p.
\]
In a similar way we can prove other cases.
Therefore we see that for $b\in \CA$ such that $E_\CA(b) \neq 0$, 
$\tfa{0}b \neq 0$ we have $E_\CA\left (\tfa{0}b \right ) \neq 0$.
Similarly we can show that for $b\in \CA$ such that $E_\CA(b) = 0$ or
$\tfa{0}b = 0$ we have $E_\CA\left (\tfa{0}b \right ) = 0$.
\vspace{3mm}\\  
We verify (C3). By calcualtion, we see  
\begin{equation}
\wt^{\CA} \left ( \tf{1}^q\tf{0}^p\hwvl{i}{k}{j} \right ) =
(k+p-2q)\Lambda_1 + (j-2p+q)\Lambda_0.
\label{eqn:weightCA}
\end{equation}
If the action of $E_\CA$ is given by \saba 1, we have 
\[
\wt^{\CA} \left ( E_\CA \tf{1}^q\tf{0}^p\hwvl{i}{k}{j} \right ) =
\wt^{\CA} \left ( \tf{1}^q\tf{0}^p\hwvl{i}{k-1}{j} \right ) =
(k-1+p-2q)\Lambda_1 + (j-2p+q)\Lambda_0.
\]
Therefore, we have
\begin{equation}
-2\wt_1^\CA(E_\CA b)-\wt_0^\CA(E_\CA b)=
-2\wt_1^\CA(b)-\wt_0^\CA(b)+2.
\label{eqn:wt2fueru}
\end{equation}
Similarly we can show (\ref{eqn:wt2fueru}) when the action of $E_\CA$ is 
given by \saba2 - \saba 5.\vspace{3mm}\\
We prove (C3) for $b=\fa{1}{q}\fa{0}{p}\hwvl{i}{l-i}{j} \in \CA_+^{(l)}$ where
$i\leq j \leq l-i$, $0 \leq p \leq j$, $p\leq q \leq p+l-i$. 
By Definition \ref{definition:itinokettei}, $b$ satisfies $F_\CA(b)=0$ and
\begin{equation}
\ve{\CA}(b)=2l-2i+j-3q.
\label{eqn:phi2EAb}
\end{equation}
\renewcommand{\labelenumi}{\arabic{enumi}. }
\begin{enumerate}
\item If $0 \leq p=q \leq \left [ \frac{i+j}{2} \right ]$,\\
The action of $E_\CA$ on $\tf{1}^p\tf{0}^p\hwvl{i}{k}{j}$ $(i<k\leq l-i)$ 
is given by \saba 1. Then we have
\[
E_\CA^{l-2i}(b)=\tf{1}^p\tf{0}^p\hwvl{i}{i}{j}.
\]
If $\left [ \frac{j-i}{2} \right ]>0$,
the action of $E_\CA$ on $\tf{1}^p\tf{0}^p\hwvl{i}{i}{j}$ 
is given by \saba 3.
We put $x=\frac{i+j}{2}$, $\hat{x}=\frac{i+j+1}{2}$.
Then we have
\[
E_\CA^{l-3i+[x]}(b)=
\left \{
\begin{array}{ll}
\tf{1}^{p+\hat{x}-i-1}\tf{0}^p\hwvl{\hat{x}-1}{\hat{x}-1}{\hat{x}}& 
\mbox{$j-i$ is odd},\\
\tf{1}^{p+x-i}\tf{0}^p\hwvl{x}{x}{x}&\mbox{$j-i$ is even}.
\end{array}
\right .
\]
If $i-j$ is odd, the action of $E_\CA$ on $\tf{1}^{p+\hat{x}-i-1}\tf{0}^p
\hwvl{\hat{x}-1}{(\hat{x}-1)}{\hat{x}}$ is given by \saba 4. 
Then we have
\[
E_\CA^{l-3i+[x]}(b)=
\tf{1}^{p+\hat{x}-i}\tf{0}^p\hwvl{\hat{x}-1}{\hat{x}}{(\hat{x}-1)}.
\]
If $j-p-[\hat{x}-i]>0$, then the action of 
$E_\CA$ on $\tf{1}^{p+x-i}\tf{0}^p\hwvl{x}{x}{x}$,
$\tf{1}^{p+\hat{x}-i}\tf{0}^p\hwvl{\hat{x}-1}{\hat{x}}{\hat{x}-1}$ 
are given by \saba 5.
Then we have
\[
E_\CA^{l-2i+j-p}(b) = \tf{1}^j\tf{0}^p\hwvl{p}{j+i-p}{p}.
\]
Moreover, the action of $E_\CA$ on $\tf{1}^j\tf{0}^p\hwvl{p}{j+i-p}{p}$ 
is given by \saba 2. Then we have
\begin{equation}
E_\CA^{2l-2i+j-3p}(b) = \tf{1}^{l+j-2p}\tf{0}^{l-p}\hwvl{p}{j+i-p}{l-p}.
\label{eqn:Eave1}
\end{equation}
Finally, the action of 
$E_\CA$ on $\tf{1}^{l+j-2p}\tf{0}^{l-p}\hwvl{p}{j+i-p}{l-p}$  
is given by  \saba 6. Then we have
\[
E_\CA^{2l-2i+j-3p+1}(b) = 0.
\]
\item If $\left [ \frac{i+j}{2} \right ] < p=q < \frac{2j+i}{3}$,  
we can calculate in a similar way:
\begin{equation}
\begin{array}{rcll}
E_\CA^{l-2i}(b)&=&\tilde{f}_1^p\tilde{f}_0^p\hwvl{i}{i}{j}, &\saba 1
\vspace{1mm}\\
E_\CA^{l-i+2j-3p}(b)&=&\tilde{f}_1^{2j+i-2p}\tilde{f}_0^p
\hwvl{2j+2i-3p}{2j+2i-3p}{-j-i+3p}, &\saba 3\vspace{1mm}\\
E_\CA^{2l-2i+j-3p}(b)&=&\tilde{f}_1^{l+j-2p}\tilde{f}_0^{l-i-j+p}
\hwvl{2j+2i-3p}{2j+2i-3p}{l-2j-2i+3p}, &\saba 2\vspace{1mm}\\
E_\CA^{2l-2i+j-3p+1}(b)&=&0. &\saba 6
\end{array}
\label{eqn:Eave2}
\end{equation}
\item If $\frac{2j+i}{3} < p=q \leq j$ or $j=p < q$, 
we have
\begin{equation}
\begin{array}{rcll}
E_\CA^{l-i+2j-3q}(b)&=&\tilde{f}_1^q\tilde{f}_0^p\hwvl{i}{3q-2j}{j}, 
&\saba 1\vspace{1mm}\\
E_\CA^{2l-2i+j-3q}(b)&=
&\tilde{f}_1^{l-i-j+q}\tilde{f}_0^{l-i-j+p}\hwvl{i}{3q-2j}{l-i},&\saba 2
\vspace{1mm}\\
E_\CA^{2l-2i+j-3q+1}(b)&=&0. &\saba 6
\end{array}
\label{eqn:Eave3}
\end{equation}
\item If $p<j=i$, $p<q<y+j$, we have
\begin{equation}
\begin{array}{rcll}
E_\CA^{l-2i-3(p-q)}(b)&=&\tilde{f}_1^q\tilde{f}_0^p
\hwvl{i}{i+3(q-p)}{i}, &\saba 1\vspace{1mm}\\
E_\CA^{l-i-3q+2p}(b)&=&\tilde{f}_1^{q+i-p}\tilde{f}_0^p
\hwvl{p}{2i+3q-4p}{l-p}, &\saba 5\vspace{1mm}\\
E_\CA^{2l-i-3q}(b)&=&\tilde{f}_1^{l+i-3p+q}\tilde{f}_0^{l-p}
\hwvl{p}{2i+3q-4p}{l-p}, &\saba 2\vspace{1mm}\\
E_\CA^{2l-i-3q+1}(b)&=&0.&\saba 6
\end{array}
\end{equation}
\item Other cases are reduced to $\CG^{l-1}$, using (\ref{eqn:indsayou}).
\end{enumerate}

Then for
$b=\fa{1}{q}\fa{0}{p}\hwvl{i}{l-i}{j} \in \CA_+^{(l)}$ where
$i\leq j \leq l-i$, $0 \leq p \leq j$, $p\leq q \leq p+l-i$, we have 
\[
-2\wt_1^\CA(b)-\wt_0^\CA(b)=-(2l-2i+j-3q).
\]
By (\ref{eqn:phi2EAb}), we see 
\[
\vp{\CA}(b)-\ve{\CA}(b)=-(2l-2i+j-3q).
\]
Therefore by (\ref{eqn:wt2fueru}), we verify inductively that $E_\CA$ satisfies
(C3).

\subsubsection{An involution on $\CA$}
\label{section:constsymmetry}
Let $\lwv{l}{i}{(j_1,j_0)}$ be the
lowest weight element in $\CA_i$
with weight $j_1\Lambda_1 + j_0 \Lambda_0$.
We define a map
$\mbox{\rm C}_\CA : \CA \rightarrow \CA$ by: 
\begin{equation}
\begin{array}{rcl}
\syme{\fa{0}{r}\fa{1}{q}\fa{0}{p}\hwvl{i}{k}{j}} &=& 
\ea{0}{r}\ea{1}{q}\ea{0}{p}
\lwv{l}{i}{(-k,-j)}\vspace{2mm}\\
&=& \fa{0}{j+q-2p-r}\fa{1}{k+j-q}\fa{0}{k-q+p}\hwvl{i}{j}{k}.
\end{array}
\label{eqn:gyaku}
\end{equation}
where $0 \leq p \leq j$, $p\leq q \leq p+k$, $0\leq r \leq j+q-2p$.
It is easy to see that $\mbox{C}_\CA$ is an involution, 
\[
\syme{\syme{b}}=b \:\:(b\in \CA).
\]
\begin{proposition}
\label{prop:symEF} 
For $b \in \CA$, we have
\begin{equation}
\syme{E_\CA b}= F_\CA(\syme{b}).
\label{eqn:houkousei}
\end{equation}
\end{proposition}
\proof
It is sufficient to prove (\ref{eqn:houkousei}) for $b \in \CA$ such that 
$\ve{0}^\CA(b)=0$.\\
We assume that for $b=\fa{1}{q}\fa{0}{p}\hwvl{i}{k}{j} \in \CA$, $k$, $j$, $p$,
$q$ satisfy $i+1\leq k \leq l-i$, $i+1 \leq j \leq l-i$, 
$0\leq q \leq j-1-\left [ \frac{j-k}{3} \right ]$, $p=\mbox{min}(q,j)$.
Since this is the condition (1)(i) in Definition 
\ref{definition:itinokettei}, the action of $E_\CA$ is given by \saba 1. 
Then we have 
\[
\syme{E_\CA b} = \fa{0}{j+q-2p}\fa{1}{k+j-q-1}\fa{0}{k-q+p-1}
\hwvl{i}{j}{k-1}.
\]
We consider the action of $F_\CA$ on $\syme{b}$. 
By applying the involution, we have
\[
\syme{b}=\fa{0}{j+q-2p}\fa{1}{k+j-q}\fa{0}{k-q+p}\hwvl{i}{j}{k}.
\]
Thus the action of $F_\CA$ on $\syme{b}$ is given by \saba 2$'$.
Therefore we have
\[
F_\CA \left ( \syme{b} \right ) =
\fa{0}{j+q-2p}\fa{1}{k+j-q-1}\fa{0}{k-q+p-1}\hwvl{i}{j}{k-1}.
\]
In a similar way, we have \vspace{2mm}\\
if the action of $E_\CA{b}$ is given by \saba 1, then the action of 
$F_\CA\left ( \syme{b} \right )$ is given by \saba 2$'$,\vspace{2mm}\\
if the action of $E_\CA{b}$ is given by \saba 2, then the action of 
$F_\CA\left ( \syme{b} \right )$ is given by \saba 1$'$,\vspace{2mm}\\
if the action of $E_\CA{b}$ is given by \saba 3, then the action of 
$F_\CA\left ( \syme{b} \right )$ is given by \saba 5$'$,\vspace{2mm}\\
if the action of $E_\CA{b}$ is given by \saba 4, then the action of 
$F_\CA\left ( \syme{b} \right )$ is given by \saba 4$'$,\vspace{2mm}\\
if the action of $E_\CA{b}$ is given by \saba 5, then the action of 
$F_\CA\left ( \syme{b} \right )$ is given by \saba 3$'$,\vspace{2mm}\\
if the action of $E_\CA{b}$ is given by \saba 6, then the action of 
$F_\CA\left ( \syme{b} \right )$ is given by \saba 6$'$.\vspace{2mm}\\
Thus we have $\syme{E_\CA(b)}=F_\CA\left (\syme{b} \right )$.
\hfill$\Box$\vspace{5mm}

For $a \in \Z$, we define $a_+$ by
\[
a_+ = \left \{
\begin{array}{ll}
a&(a>0),\\
0&(a\leq 0).
\end{array} \right .
\]

We set $B_C$, $B_W$, $B_U$, $B_R \subset \CA^{(l)}$ by
\[
B_C=\left \{ \sayouA{q}{p}\hwvl{i}{l-i}{j}\:\left |\:
0\leq i \leq \left[\frac{l}{2} \right ],\:i\leq j\leq l-i,
0\leq q \leq p \leq j \right . \right \}, 
\]
\[
B_W=\left \{ \sayouA{q}{p}\hwvl{i}{l-i}{j}\:\left |\:
0\leq i \leq \left[ \frac{l}{2} \right ], \: j = i, 
\:0 \leq p \leq j, \:p < q \leq y+j \right . \right \}, 
\]
\[
B_U=\left \{ \sayouA{q}{p}\hwvl{i}{l-i}{j}\:\left |\:
0\leq i \leq \left[\frac{l}{2} \right ],\:i\leq j\leq l-i, \: p = j, \:
j < q \leq y + 2j -i \right . \right \}, 
\]
\[
B_R=\left \{ \sayouA{q}{p}\hwvl{i}{l-i}{j}\:\left |\:
\begin{array}{l}
0\leq i \leq \left[\frac{l}{2} \right ],\:i < j \leq l-i,\\
0 \leq q < j, \:p < q \leq y+j-(i-p)_+ 
\end{array}
\right . \right \}, 
\]
where $y=y\left ( \fa{1}{q}\fa{0}{p}\hwvl{i}{l-i}{j} \right )=\left [ 
\frac{l-i-j}{3} \right ]$. 

In view of Definition \ref{definition:itinokettei}, 
the conditions  $\ve{0}^{\CA}(b) = 0$, $\vp{\CA}(b) = 0$  
on  $b=\sayouA{q}{p}\hwvs$ 
$(0 \leq i \leq \left [ \frac{l}{2} \right ], i \leq j \leq l-i, 
0 \leq p \leq j )$ are rephrased as
\[
\left \{
\begin{array}{ll}
0 \leq q \leq p - \left [ \frac{-l+i+j}{3} \right ] +(j-i)&
(p \leq i,\:l-i-j \equiv 0 \pmod{3}),\\
0 \leq q \leq p -1 - \left [ \frac{-l+i+j}{3} \right ] +(j-i)&
(p \leq i,\:l-i-j \not\equiv 0 \pmod{3}),\\
0 \leq q \leq j - \left [ \frac{-l+i+j}{3} \right ] & 
(p > i,\:l-i-j \equiv 0 \pmod{3}),\\
0 \leq q \leq j -1 - \left [ \frac{-l+i+j}{3} \right ] & 
(p > i,\:l-i-j \not\equiv 0 \pmod{3}).
\end{array}
\right .
\]
Summarizing, we have
\[
0 \leq q \leq \left [ \frac{l-i-j}{3} \right ] + j - (i-p)_+ =
y+j -(i-p)_+.
\]
Thus we have
\[
\begin{array}{l}
\displaystyle
\left \{ b \in \CA \left | \ve{0}^{\CA}(b)=0, 
\vp{\CA}(b)=0 \right . \right \} \\
\displaystyle \; \; = \left \{ \sayouA{q}{p}\hwvl{i}{l-i}{j} \left | 
0 \leq i \leq \left [ \frac{l}{2}\right ], i \leq j \leq l-i, 
0 \leq p \leq j, p \leq q \leq y+j-(i-p)_+ \right . \right \}.
\end{array}
\]
Hence, we see easily that 
\[
\left \{ b \in {\cal A} \left | \ve{0}^{\cal A}(b)=0, 
\vp{\CA}(b)=0 \right . \right \} \subset
B_C \cup B_W \cup B_U \cup B_R.
\]
By Proposition \ref{proposition:kakanA2a}, we have 
\begin{equation}
\left \{ \left . \fa{0}{r}\fa{1}{p}\fa{0}{p}\hwvl{i}{l-i}{j} \right |
 0 \leq p \leq j, 0 \leq r \leq j-p \right \} = B_C.
\label{eqn:kakanA2c}
\end{equation}

\begin{proposition}
\label{prop:inBC}
For $b\in B_W \cup B_U$, we have 
\[
\syme{E_\CA^{\ve{\CA}(b)}(b)}\in B_C.
\] 
\end{proposition}

\proof 
By calculation for $b = \sayouA{q}{p}\hwvl{i}{l-i}{i} \in B_W$ 
$(0 \leq p < i, p < q \leq y+i)$ we have
\[
\syme{E_\CA^{\ve{\CA}(b)}b}
= \sayouA{i+2q-2p}{2i+3q-4p}\hwvl{p}{l-p}{2i+3q-4p}. 
\]
Since $p < q$, $p < i$, we have
\[
2i+3q-4p -(i+2q-2p) = (i-p)+(q-p) > 0.
\]
Then we have
\[
\syme{E_\CA^{\ve{\CA}(b)}b}\in B_C.
\]
In a similar way,  for $b \in B_U$ we have
\[
\syme{E_\CA^{\ve{\CA}(b)}b}\in B_C.
\]
\hfill$\Box$

Similar to Proposition \ref{prop:inBC}, we have following proposition:
\begin{proposition}
\label{prop:inBR}
For $b\in B_R$, we have
\[
\syme{E_\CA^{\ve{\CA}(b)}(b)} \in B_R.
\]
\end{proposition}

\subsubsection{Definition of $\Phi$ on $\CA$ and $\tfg{0}$ on $\CG^l$}
\label{sec:defPhi}
In this section, we define one-to-one map $\Phi\,:\,\CA \rightarrow \CG^l$
and Kashiwara operator $\tfg{0}$ on $\CG^l$ expoiting operators on $\CA$.

Since roots $\alpha_0$ and $\alpha_2$ of algebra $\U{\gone}$ are 
orthogonal, the Kashiwara operator $\tfg{0}$ commutes with $\tfg{2}$. 
\vspace{2mm}
\renewcommand{\labelenumi}{\rm (D\arabic{enumi}) }
\begin{enumerate}
\item  $\tfg{0}\tfg{2}(b)=\tfg{2}\tfg{0}(b)$,  
$\teg{0}\teg{2}(b)=\teg{2}\teg{0}(b)$, 
for $b \in \CG^l$.
\end{enumerate}
In view of Proposition \ref{proposition:dimension}, we are going to construct 
the one-to-one map 
\[
\Phi\::\: \CA \rightarrow \CG^l,
\]
with the following properties; 
\renewcommand{\labelenumi}{\rm (E\arabic{enumi}) }
\begin{enumerate}
\item $\teg{1}\Phi(b) = \Phi\left (\tea{1}b \right )$, for $b\in \CA$,
\item $\tfg{2}\Phi(b) = \Phi \left ( F_\CA b \right )$, 
 for $b\in \CA$,
\item $\wt_1^\CG\left (\Phi(b) \right)=\wt_1^\CA(b)$, for $b\in \CA$,
\item $\wt_2^\CG\left (\Phi(b) \right )= -2\wt_1^\CA(b)-\wt_0^\CA(b)$,
for $b\in \CA$,
\item $\teg{0}\Phi(b)=0$ (resp. $\tfg{0}\Phi(b)=0$) 
if and only if 
$\tea{0}(b)=0$ (resp. $\tfa{0}(b)=0$), for $b \in \CA$.
\end{enumerate}
By (E1) and (E2), we see easily
\renewcommand{\labelenumi}{\rm (E\arabic{enumi})$'$ }
\begin{enumerate}
\item $\tfg{1}\Phi(b) = \Phi\left (\tfa{1}b \right )$, for $b \in \CA$,
\item $\teg{2}\Phi(b) = \Phi\left ( E_\CA b \right )$, for $b \in \CA$.
\end{enumerate}
We call an element $b\in \CA$  {\it $\CA$-terminal} if $F_\CA(b)=0$. 
An {\it $E_\CA$-string} is a sequence of elements in $\CA$, $b$ , $E_\CA(b)$, 
$E_\CA^2(b)$, $\ldots$, $E_\CA^n(b)$, where $b$ is $\CA$-terminal and 
$E_\CA^{n+1}(b)=0\;(n\in {\bf Z}_{\geq 0})$. 
We first define $\Phi$ for $\CA$-terminal elements using (E3)(E4) mainly.
Next we define the the action of $\Phi$ for elements of an $E_\CA$-string using 
(E2). We verify (E1) completely in \S\ref{section:commutative}.

\begin{definition}
\label{fact:URelement}
We define
\begin{eqnarray}
\Phi\left ( \fa{0}{p}\hwvl{0}{l}{l}\right ) &=& 
\BB{\bb{6}{p}\bb{-2}{l-p}}\;\;( 0\leq p \leq l),\\
\Phi\left ( \ea{0}{p}\lwv{l}{0}{(-l,-l)} \right ) &=&
\BB{\bb{2}{l-p}\bb{-6}{p}}\;\;( 0\leq p \leq l).
\end{eqnarray}
\end{definition}
We verify the properties (E1)--(E5).

By the definition of $\hwvl{0}{l}{l}$, 
\[
\wt^\CA \left ( \hwvl 0ll \right ) = l\Lambda_1+l\Lambda_0.
\]
By (E3) and (E4), we see that
$\Phi \left (\hwvl 0ll \right )$ is the element with weight 
$l\Lambda_1-3l\Lambda_2$.
We can see easily that $b=\BB{\bb{6}{k}\bb{-2}{l-k}}$ $(0 \leq k \leq l)$
is the only element such that $\wt_2^\CG(b)=-3l$.
Moreover $b=\BB{\bb{-2}{l}}$ is the only element such that $\wt_1^\CG(b)=l$.
The highest weight element is only element in $B^i_{(k,j)}$ that satisfies 
$\wt(b)=k'\Lambda_1+j'\Lambda_0$ $(k' \geq k, j' \geq j)$.
By the definition of $\CA$, the element $b\in \CA$ with
$\wt(b)=l\Lambda_1+l\Lambda_0$ is $\hwvl{0}{l}{l}$.  
Therefore, we have
\begin{equation}
\Phi \left ( \hwvl 0ll \right ) = \BB{\bb{-2}{l}}.
\label{eqn:i0hashi}
\end{equation}
We verify (E1).
Since $\tea{1}\hwvl{0}{l}{l}=0$, we verify  $\teg{1}\BB{\bb{-2}{l}}=0$.
Using Proposition \ref{proposition:vplus1}, 
\[
\uu{1}\biggl ( \BB{\bb{-2}{l}} \biggr ) = u_+^l.
\]
Then we have 
\[
\teg{1}\BB{\bb{-2}{l}} = 0
\]
In a similar way, we verify (E2). 
We consider $\fa{0}{p}\hwvl{0}{l}{l}$. 
We see that the element $b \in \CG^l$ which satisfy following 
formula is unique: 
\[
\wt^\CG(b)=
(l+k)\Lambda_1-3l\Lambda_2.
\]
Thus, we can define
\begin{equation}
\Phi \left ( \fa{0}{p}\hwvl0ll \right ) = \BB{\bb{6}{p}\bb{-2}{l-p}}
\mbox{ ($0 \leq p \leq l$)}.
\label{eqn:i0hashi0p}
\end{equation}
Similarly, by the calculation of weight, we can define uniquely 
\[
\Phi\left ( \ea{0}{p}\lwv{l}{0}{(-l,-l)} \right ) =
\BB{\bb{2}{l-p}\bb{-6}{p}}\;\;( 0\leq p \leq l).
\]
\hfill$\Box$

By (E5) and (\ref{eqn:i0hashi0p}), we define the action of $\teg{0}$ and 
$\tfg{0}$ for $\BB{\bb{6}{p}\bb{-2}{l-p}}$ $(0 \leq p \leq l)$ by
\begin{equation}
\tfg{0}\BB{\bb{6}{p}\bb{-2}{l-p}}=\left \{
\begin{array}{ll}
\BB{\bb{6}{p+1}\bb{-2}{l-p-1}}& 0 \leq p \leq l-1 \vspace{2mm}\\
0&p=l
\end{array}
\right .
\label{eqn:f0gUR}
\end{equation}
\begin{equation}
\teg{0}\BB{\bb{6}{p}\bb{-2}{l-p}}=\left \{
\begin{array}{ll}
\BB{\bb{6}{p-1}\bb{-2}{l-p+1}}& 1 \leq p \leq l \vspace{2mm}\\
0&p=0
\end{array}
\right .
\label{eqn:e0gUR}
\end{equation}

In view of Definition \ref{definition:itinokettei}, \ref{fact:URelement} 
and (E2)$'$, following definition is led.
\begin{definition}
\label{fact:Uelement}
We define 
\begin{equation}
\Phi\left ( \fa{0}{p}\hwvl{0}{k}{l} \right ) = 
\eg{2}{l-k}\;\BB{\bb{6}{p}\bb{-2}{l-p}}.
\label{eqn:Uelement}
\end{equation}
where $0 \leq k,p \leq l$.
\end{definition}
We verify the properties (E1)--(E5). The property (E2) is obvious.
We verify (E1)(E3).
We see 
\[
\ve{1}^\CA \left (\fa{0}{p}\hwvl{0}{k}{l} \right )=0,
\]
\[
\vp{1}^\CA \left (\fa{0}{p}\hwvl{0}{k}{l} \right )=p+k.
\]
On the other hand,
\[
\begin{array}{ll}
\eg{2}{l-k}\BB{\bb{6}{p}\bb{-2}{l-p}}\\
&\makebox[-2cm]{}= \left \{
\begin{array}{ll}
\BB{\bb{2}{(l-k)/3}\bb{6}{p-(l-k)/3}\bb{-2}{l-p}}&
0 \leq l-k \leq 3p, \;\;l-k\equiv 0 \pmod{3}\\
\BB{\bb{2}{[(l-k)/3]}\bb{4}{}\bb{6}{p-1-[(l-k)/3]}\bb{-2}{l-p}}&
0 < l-k < 3p, \;\;l-k\equiv 1 \pmod{3}\\
\BB{\bb{2}{[(l-k)/3]}\bb{3}{}\bb{6}{p-1-[(l-k)/3]}\bb{-2}{l-p}}&
0 < l-k < 3p, \;\;l-k\equiv 2 \pmod{3}\\
\BB{\bb{2}{p}\bB{\CWb}{l-k-3p}\bb{-2}{k+2p}}&
l-k>3p
\end{array}
\right .
\end{array}
\]
Using Proposition\ref{proposition:vplus1}, we have
\begin{equation}
\ve{1} \biggl ( \eg{2}{l-k}\;\BB{\bb{6}{p}\bb{-2}{l-p}} \biggr ) = 0,
\end{equation}
\begin{equation}
\vp{1} \biggl ( \eg{2}{l-k}\;\BB{\bb{6}{p}\bb{-2}{l-p}} \biggr ) = p+k.
\label{eqn:ve2Rk}
\end{equation}
By calculation, we can verify (E4). 
We verify (E5).
By (\ref{eqn:e0gUR}) and (D1), we must have 
\begin{equation}
\teg{0}\biggl ( \eg{2}{l-k}\;\BB{\bb{-2}{l}} \biggr ) = 0
\label{eqn:ve0Rk}
\end{equation}
\hfill$\Box$\vspace{5mm}

By Definition \ref{fact:Uelement} and (E1)$'$, 
following definition is led.
\begin{definition}
\label{def:hamidashi}
We define
\[
\Phi\biggl( \fa{1}{q}\fa{0}{l}\hwvl{0}{k}{l} \biggr)=
\fg{1}{q} \biggl( \Phi\biggl( \fa{0}{l}\hwvl{0}{k}{l} \biggr)
\;\; (0\leq p \leq l+j).
\]
\end{definition}
We can verify (E1),(E2), (E3), (E4).
We see easily
\[
\tfa{0}\left ( \fa{1}{q}\fa{0}{l}\hwvl{0}{k}{l} \right )=0 \;\;
(0 \leq q \leq l).
\]
By (\ref{eqn:A2char1}), we see
\[
\tea{0}\left ( \fa{1}{q}\fa{0}{l}\hwvl{0}{k}{l} \right )=0 \;\;
(l \leq q \leq l+j).
\]
In order to satisfy (E4), we define 
\[
\tfg{0}\left (\Phi\left ( \fa{1}{q}\fa{0}{l}\hwvl{0}{k}{l} \right ) \right )=0 
\;\;(0 \leq q \leq l),
\]
\[
\teg{0}\left(\Phi\left ( \fa{1}{q}\fa{0}{l}\hwvl{0}{k}{l} \right )\right )=0 
\;\;(l \leq q \leq l+j).
\]

For $b \in B^i_{(k,j)}$ we define $y(b) \in \Z_{\geq 0}$ by 
\[
y(b)= \left [\frac{l-i-j}{3} \right ].
\]
We often write $y$ instead of $y(b)$ for simplicity.
\begin{definition}
\label{fact:allhighest}
For $\hwvl{i}{l-i}{j}$, we define
\begin{equation}
\Phi\left ( \hwvl{i}{l-i}{j}\right ) = \left \{
\begin{array}{ll}
\BB{\bb{6}{y}\bB{\CC}{y+i}\bb{-2}{y+i}\bb{-2}{j-i}}&
l-i-j \equiv 0 \pmod{3}, 
\vspace{2mm}\\
\BB{\bb{6}{y+1}\bB{\CC}{y+i-1}\bb{4}{}\bb{-2}{y+i}\bb{-2}{j-i}}&
\left (
\begin{array}{l}
l-i-j \equiv 1 \pmod{3}\\
y + i > 0
\end{array}
\right ), \vspace{2mm}\\ 
\BB{\bb{-5}{}\bb{-2}{j-i}} & 
\left (
\begin{array}{l}
l-i-j \equiv 1 \pmod{3}\\
y + i = 0
\end{array}
\right ),
\vspace{2mm}\\ 
\BB{\bb{6}{y+1}\bB{\CC}{y+i}\bb{-3}{}\bb{-2}{y+i}\bb{-2}{j-i}} &
(l-i-j \equiv 2 \pmod{3}).
\end{array} \right .
\label{eqn:highests}
\end{equation}
\end{definition}
We verify the properties (E1)--(E5). We verify (E1) and (E3).
We see easily
\[
\tea{1}\left (\hwvl{i}{l-i}{j} \right ) = 0,\:
\tfg{1}\left (\hwvl{i}{l-i}{j} \right ) = l-i.
\]
By proposition \ref{proposition:vplus1}, we have
\begin{eqnarray*}
\uu{1}\biggl (\BB{\bb{6}{y}\bB{\CC}{y+i}\bb{-2}{y+j}}\biggr )&=& 
u_+^{y+j}u_-^{y+i}u_+^{y+i}u_+^{2y},\vspace{2mm}\\
\mbox{Red}_1\biggl (\BB{\bb{6}{y}\bB{\CC}{y+i}\bb{-2}{y+j}}\biggr )&=&
u_+^{l-i}.
\end{eqnarray*}
Then we have
\[
\teg{1}\biggl ( \BB{\bb{6}{y}\bB{\CC}{y+i}\bb{-2}{y+j}} \biggr )=0,
\]
\[
\tfg{1}\biggl ( \BB{\bb{6}{y}\bB{\CC}{y+i}\bb{-2}{y+j}} \biggr) =l-i.
\]
In a similar way, we have 
\[
\begin{array}{ll}
\teg{1}\biggl (\BB{\bb{6}{y+1}\bB{\CC}{y+i-1}\bb{-4}{}\bb{-2}{y+j}}\biggr )=0,&
\tfg{1}\biggl (\BB{\bb{6}{y+1}\bB{\CC}{y+i-1}\bb{-4}{}\bb{-2}{y+j}}\biggr )=l-i,
\vspace{2mm}\\
\teg{1}\biggl (\BB{\bb{-5}{}\bb{-2}{j-i}}\biggr )=0,&  
\tfg{1}\biggl (\BB{\bb{-5}{}\bb{-2}{j-i}}\biggr )=l-i,  
\vspace{2mm}\\
\teg{1}\biggl ( \BB{\bb{6}{y+1}\bB{\CC}{y+i}\bb{-3}{}\bb{-2}{y+j}}\biggr )=0,&
\teg{1}\biggl ( \BB{\bb{6}{y+1}\bB{\CC}{y+i}\bb{-3}{}\bb{-2}{y+j}}\biggr )=l-i,
\end{array}
\]
By calculation, we can verify (E2) and (E4).
We verify (E5). We put $b=\BB{\bb{6}{y}\bB{\CC}{y+i}\bb{-2}{y+j}}$.
By calculation and Lemma \ref{fact:Uelement}, we have
\begin{eqnarray*}
\eg{2}{2l-2i+j}b &=& \BB{\bb{2}{y(b)}\bB{\CC}{y(b)+i}\bb{-6}{y(b)+j}},\\
\eg{1}{l+j}\eg{2}{2l-2i+j}b &=& \BB{\bb{2}{y(b)}\bb{6}{2y(b)+i+j}},\\
\fg{2}{l-i-j}\eg{1}{l+j}\eg{2}{2l-2i+j}b &=& \BB{\bb{6}{l}}.
\end{eqnarray*}
By Definition \ref{def:hamidashi}, we must have 
\[
\Phi \left ( \fa{1}{l+j}\fa{0}{l}\hwvl{0}{i+j}{l} \right ) =
\BB{\bb{2}{y}\bB{\CC}{y+i}\bb{-6}{y+j}},
\]
\[
\teg{0}\,\BB{\bb{2}{y}\bB{\CC}{y+i}\bb{-6}{y+j}}=0.
\]
By (D1), we must have 
\[
\teg{0}\BB{\bb{6}{y}\bB{\CC}{y+i}\bb{-2}{y+j}}=0.
\]
Similarly, we can prove other cases. 
\hfill $\Box$\vspace{5mm}

\begin{definition}
\label{fact:f0pR}
We set $b=\hwvs$. We define the action of $\Phi$ as follows:\\
If $0\leq p\leq i$,
\[
\Phi\left ( \fa{0}{p} b \right ) \,=\, \left \{
\begin{array}{ll}
\BB{\bb{1}{p}\bb{6}{y}\bB{\CC}{y+i-p}\bb{-2}{y+j}}&
(l-i-j \equiv 0 \pmod{3}),\vspace{2mm}\\
\BB{\bb{1}{p}\bb{6}{y+1}\bB{\CC}{y+i-p-1}\bb{-4}{}\bb{-2}{y+j}}&
(y+i>p,\;l-i-j \equiv 1 \pmod{3}),\vspace{2mm}\\ 
\BB{\bb{1}{i}\bb{-5}{}\bb{-2}{j}}&
(y +i = p,\;l-i-j \equiv 1 \pmod{3}),\vspace{2mm}\\
\BB{\bb{1}{p}\bb{6}{y+1}\bB{\CC}{y+i-p}\bb{-3}{}
\bb{-2}{y+j}}&(l-i-j \equiv 2 \pmod{3}),
\end{array} \right .
\]
if $i<p\leq j$, 
\[
\Phi\left ( \fa{0}{p} b \right ) \,=\,
\left \{
\begin{array}{ll}
\BB{\bb{1}{i}\bb{6}{p-i+y}\bB{\CC}{y}\bb{-2}{y+j-p+i}}&
(l-i-j\equiv 0 \pmod{3}), \vspace{2mm}\\
\BB{\bb{1}{i}\bb{6}{p-i+y+1}\bB{\CC}{y-1}\bb{-4}{}
\bb{-2}{y+j-p+i}}&(y > 0, l-i-j\equiv 1 \pmod{3}), 
\vspace{2mm}\\ 
\BB{\bb{1}{i}\bb{6}{p-i}\bb{-5}{}\bb{-2}{j-p+i}}&
(y = 0, l-i-j\equiv 1 \pmod{3}), \vspace{2mm}\\
\BB{\bb{1}{i}\bb{6}{p-i+y+1}\bB{\CC}{y}\bb{-3}{}\bb{-2}{y+j-p+i}}&
(l-i-j\equiv 2 \pmod{3}).
\end{array} \right .
\]
\end{definition}
By calculation, we can verify (E1)--(E4). 
In order to satisfy (E5), we define $\teg{0}$ on 
$\Phi\left ( \fa{1}{p}\fa{0}{p}\hwvl{i}{l-i}{j} \right )$ 
$(0 \leq p \leq j)$ and $\tfg{0}$ on 
$\Phi\left ( \fa{1}{q}\fa{0}{j}\hwvl{i}{l-i}{j} \right )$ $(0 \leq q \leq j)$.

Using Proposition \ref{prop:inBC}, \ref{prop:inBR} and (E2),
it is enough to define the action of $\Phi$ on $B_C \cup B_R$. 

By (E1)$'$, we are led to the following definition.
\begin{definition}
\label{lemma:f1qf0p}
We define
\[
\Phi \left ( \fa{1}{q}\fa{0}{p}\hwvl{i}{l-i}{j} \right ) =
\left ( \tfg{1} \right )^q \left ( \Phi \left ( \fa{0}{p}\hwvl{i}{l-i}{j}
\right ) \right ),
\]
where $(0\leq p \leq j,\;0 \leq q \leq y+j-(i-p)_+)$.   
\end{definition}
Thus we have defined $\Phi(b)$ where $b \in B_C \cup B_R$.
\begin{definition}
\label{def:symg}
For $b=\BB{\bB{b_1}{}\cdots\bB{b_n}{}} \in B^{G_2}(n\Lambda_1)$ 
$(0 \leq n \leq l)$, we define 
the involution on $\CG^l$ by
\[
\symg{b}=\BB{\bB{\overline{b_n}}{}\cdots\bB{\overline{b_1}}{}},
\]
where if $b_i = \overline{m}$ $(m=1,\ldots, 6)$ then 
$\bB{\overline{b_i}}{}=\bB{m}{}$, 
if $b_i=0_1, 0_2$, then $\bB{\overline{b_i}}{}=\bB{b_i}{}$.
\end{definition}
\begin{remark}
Let $\overline{b}$ (resp. $\underline{b}$) be the highest (resp. lowest) weight 
element of $B^{G_2}(l\Lambda_1)$. 
By Proposition \ref{proposition:vplus1}, \ref{proposition:g2graph} and
\ref{proposition:gtwocrystal}, we have 
\[
\overline{b}=\BB{\bb{1}{l}},
\]
\[
\underline{b}=\BB{\bb{-1}{l}}.
\]
By Defintion \ref{def:symg}, we have
\[
\mbox{C}_\CG\bigl ( \overline{b} \bigr )=\underline{b}.
\]
Using Proposition \ref{proposition:vplus1} again, we see 
\[
\symg{\tf{i_N}\cdots \tf{i_1}\overline{b}}=
\te{i_N}\cdots\te{i_1}\underline{b}.
\] 
\end{remark}
\begin{definition}
\label{def:extea}
We define $\Phi$ for $b \in B_C \cup B_R$
\begin{equation}
\Phi(E_\CA^m b)=\eg{2}{m} \left ( \Phi (b) \right )\;\;
(0 \leq m \leq \ve{\CA}(b)).
\label{eqn:exte2}
\end{equation}
We define $\Phi$ for $b \in B_W \cup B_U$, 
\begin{equation}
\Phi(E_\CA^m b)=\symg{\Phi \left ( \syme{E_\CA^m b}\right )}\;\;
(0 \leq m \leq \ve{\CA}(b)).
\label{eqn:barsym}
\end{equation}
For $b \in B_R$, we define 
\begin{equation}
\Phi\left ( E_\CA^m \left ( \fa{0}{\vp{0}(b)}(b)\right ) \right )
= \symg{ \Phi \left ( F_\CA^{\ve{\CA}(b)-m}\syme{b} \right )}.
\end{equation}
\end{definition}
\begin{remark}
By Definition \ref{def:extea},
we see that the relation (\ref{eqn:exte2}) satisfies {\rm (E2)}.
By (\ref{eqn:exte2})  and Proposition \ref{prop:inBC}, 
we see that $\syme{E_\CA^mb}$ is already defined. 
By  (\ref{eqn:barsym}) ,We see $\syme{\syme{E_\CA^mb}}=b$ 
$(0 \leq m \leq \ve{\CA}(b))$.\vspace{3mm}
\end{remark}
We have defined $\Phi(b)$ where $b\in \CA^{(l)}$.
\begin{proposition}
For $b\in \CA^{(l)}$,
\[
\symg{\Phi(b)}=\Phi(\syme{b}).
\]
\end{proposition}
\proof
We put $b=\sayouA{q}{j}\hwvs$. By (\ref{eqn:Eave1}), (\ref{eqn:Eave2}) and
(\ref{eqn:Eave3}), 
\begin{equation}
E_\CA^{\ve{\CA}(b)}\left ( \ea{0}{j-q}b \right ) = \left \{
\begin{array}{ll}
\fa{1}{l+j-2q}\fa{0}{l-q}\hwvl{q}{j+i-q}{l-q} & \mbox{ if }
\displaystyle 0 \leq q \leq \left [ \frac{i+j}{2} \right ],
\vspace{3mm}\\
\fa{1}{l+j-2q}\fa{0}{l-i-j+q}\hwvl{2i+2j-3q}{2i+2j-3q}{l-2i-2j+3q} & 
\mbox{ if }
\displaystyle \left [ \frac{i+j}{2} \right ] < q \leq \frac{2j+i}{3},
\vspace{3mm}\\
\fa{1}{l-i-j+q}\fa{0}{l-i-j+q}\hwvl{i}{-2j+3q}{l-i} & \mbox{ if }
\displaystyle \frac{2j+i}{3} < q \leq j.
\end{array}
\right .
\label{eqn:saki}
\end{equation}
If $0 \leq q \leq \left [ \frac{i+j}{2} \right ]$,
\begin{eqnarray*}
\Phi \left ( \fa{0}{j}\hwvl{i}{l-i}{j} \right ) &=&
\BB{\bb{1}{i}\bb{6}{j-i+y}\bB{\CC}{y}\bb{-2}{y+i}},
\vspace{2mm}\\
\Phi \left ( \sayouA{q}{j}\hwvl{i}{l-i}{j} \right ) &=&
\left \{
\begin{array}{ll}
\BB{\bb{1}{i}\bb{6}{j-i+y}\bB{\CC}{y}\bb{-2}{y+i-q}\bb{-1}{q}} & 
(i > q), \vspace{2mm}\\
\BB{\bb{1}{i}\bb{6}{j-q+y}\bB{\CC}{y+q-i}\bb{-2}{y}\bb{-1}{i}}&
(i \leq q), \vspace{2mm}
\end{array}
\right . \vspace{2mm}\\
\Phi \left ( E_\CA^{\ve{\CA}(b)}\sayouA{q}{j}\hwvl{i}{l-i}{j} \right ) &=&
\left \{
\begin{array}{ll}
\BB{\bb{1}{i}\bb{2}{j-i+y}\bB{\CC}{y}\bb{-6}{y+i-q}\bb{-1}{q}} & 
(i > q), \vspace{2mm}\\
\BB{\bb{1}{i}\bb{2}{j-q+y}\bB{\CC}{y+q-i}\bb{-6}{y}\bb{-1}{i}}&
(i \leq q), \vspace{2mm}
\end{array}
\right . \vspace{2mm}\\
\symg{\Phi \left ( E_\CA^{\ve{\CA}(b)}\sayouA{q}{j}\hwvl{i}{l-i}{j} 
\right )}&=& \left \{
\begin{array}{ll}
\BB{\bb{1}{q}\bb{6}{i-q+y}\bB{\CC}{y}\bb{-2}{y+j-i}\bb{-1}{i}} & 
(i > q), \vspace{2mm}\\
\BB{\bb{1}{i}\bb{6}{y}\bB{\CC}{y+q-i}\bb{-2}{y+j-q}\bb{-1}{i}}&
(i \leq q). \vspace{2mm}
\end{array}
\right . \vspace{2mm}
\end{eqnarray*}
On the other hand,
\begin{eqnarray*}
\Phi \left ( \fa{0}{i}\hwvl{q}{l-q}{i+j-q} \right ) &=&
\left \{
\begin{array}{ll}
\BB{\bb{1}{q}\bb{6}{i-q+y}\bB{\CC}{y}\bb{-2}{y+j}} & (i > q),\vspace{2mm}\\
\BB{\bb{1}{i}\bb{6}{y}\bB{\CC}{y+q-i}\bb{-2}{y+j+i-q}} &(i \leq q),\vspace{2mm}
\end{array} \right . \vspace{2mm}\\
\Phi \left ( \syme{E_\CA^{\ve{\CA}(b)}b} \right ) &=&
\Phi \left ( \sayouA{i}{i}\hwvl{q}{l-q}{i+j-q} \right ) \vspace{2mm}\\
&=& \left \{
\begin{array}{ll}
\BB{\bb{1}{q}\bb{6}{i-q+y}\bB{\CC}{y}\bb{-2}{y+j-i}\bb{-1}{i}} & 
(i > q),\vspace{2mm}\\
\BB{\bb{1}{i}\bb{6}{y}\bB{\CC}{y+q-i}\bb{-2}{y+j-q}\bb{-1}{i}} & 
(i \leq q).\vspace{2mm}
\end{array} \right . \vspace{2mm}
\end{eqnarray*}
By (\ref{eqn:barsym}), we have 
\[
\symg{\Phi\left ( E_\CA^m b \right )}=\Phi \left (\syme{E_\CA^mb}\right )
\;\;(0 \leq m \leq \ve{0}(b))
\]
In a similar way, we can calculate other cases. \hfill$\Box$

\begin{definition}
\label{def:f0teigi}
For $b \in \CG^l$, we define $\tfg{0}(b)$ and $\teg{0}(b)$ by
\[
\tfg{0}(b)=\Phi \tfa{0} \Phi^{-1}(b),
\]
\[
\teg{0}(b)=\Phi \tea{0} \Phi^{-1}(b).
\]
We define $\ve{0}^\CG(b)$ and $\vp{0}^\CG(b)$ by
\[
\ve{0}^\CG(b)=\mbox{\rm max}\left \{n \left | \eg{0}{n}b \neq 0 \right . 
\right \},
\]
\[
\vp{0}^\CG(b)=\mbox{\rm max}\left \{n \left | \fg{0}{n}b \neq 0 \right .
\right \}.
\]
\end{definition}
By (C2), it is obvious that $\tfg{0}=\Phi\tfa{0}\Phi^{-1}$,
$\teg{0}=\Phi\tea{0}\Phi^{-1}$ satisfy (D1).
\begin{remark}
By Difinition \ref{def:f0teigi}, for $b \in \CA$ we have
\[
\ve{0}^\CG \left ( \Phi (b) \right ) = \ve{0}^\CA(b),
\]
\[
\vp{0}^\CG \left ( \Phi (b) \right ) = \vp{0}^\CA(b).
\]
\end{remark}

\begin{proposition}
The action of $\tfg{0}$ is unique.
\end{proposition}
\proof Let $\Phi$, $\Phi'$ be one-to-one maps $\CA \rightarrow \CG^l$ which 
satisfy conditions (E1)(E2)(E3) in \S \ref{sec:defPhi}.
Here, by (E2) we have $\Phi F_\CA = \tfg{2}\Phi$, 
$\Phi' F_\CA = \tfg{2}\Phi'$.
We see that if for $b \in \CA$ 
\[
\Phi(\tfa{0}b) \neq \Phi'(\tfa{0}b), \;\;\mbox{ for some } b\in \CA
\] 
we have
\[
\Phi \left ( \tfa{0}F_\CA^{\vp{\CA}(b)}(b) \right ) \neq
\Phi' \left ( \tfa{0}F_\CA^{\vp{\CA}(b)}(b) \right ). 
\]
Therefore it is sufficient to verify that $\Phi(b)$ is unique, for $b \in \CA$ 
such that $F_{\CA}(b)=0$. By (E1), $\Phi\tfa{1}=\tfg{1}\Phi$,
then it is enough to verify that 
\[
\Phi \left ( \fa{0}{p}\hwvl{i}{l-i}{j} \right )\;\;
\left ( 0 \leq i \leq \left [ \frac{l}{2} \right ],\: 0 \leq j \leq l-i,\:
0\leq p \leq j \right )
\]
is unique. This is obvious by Proposition \ref{fact:f0pR}.\hfill $\Box$
\vspace{5mm}

The affine crystal $B^l$ which is constructed with $\CG$ and $\tfg{0}$ 
satisfies Proposition \ref{proposition:construction}.

%% file: commuta.tex
\subsection{Proof of commutativity of $\Phi\tfa{1}=\tfg{1}\Phi$}
\label{section:commutative}
In this section we prove (E1) in \S\ref{sec:defPhi},
namely for $b \in \CA$
\[
\begin{array}{ccc}
b & \longrightarrow & \Phi(b)\\
\downarrow & & \downarrow\\
\tfa{1} b & \longrightarrow & \tfg{1}\left ( \Phi (b) \right ).
\end{array}
\]
In the proof of commutativity we use $\tea{1}$, $\teg{1}$ instead of 
$\tfa{1}$, $\tfg{1}$ respectively. 
Consider an $E_\CA$-string $b$, $E_\CA(b)$, $\ldots$, $E_\CA^{\vp{\CA}(b)}(b)$,
where $b\in \CA$ is $\CA$-terminal. For any $b\in \CA$, we can denote   
\[
b=E_\CA^{n'} b',
\]
where $n' \in \Z_{\geq 0}$, $b' \in \CA$ is an $\CA$-terminal element, namely 
$F_\CA(b')=0$.    
Similarly, we can denote
\[
\tea{1}b =  E_\CA^{n''}b'',
\]
where $n''\in \Z_{\geq 0}$, $b'' \in \CA$ is an $A$-terminal element.
In order to show (E1), we will verify 
\begin{equation}
\teg{1}\left ( \Phi \left ( E_\CA^{n'} b' \right ) \right )=
\Phi \left ( E_\CA^{n''}b'' \right ).
\label{eqn:e1phi}
\end{equation}
Let us define $\CA_+$, $\CA_-$ by
\[
\CA_+=\bigoplus_{i = 0}^{\left [ \frac{l}{2} \right ]}
\bigoplus_{i \leq k \leq j \leq l-i} B^i_{(k,j)},
\]
\[
\CA_-=\bigoplus_{i = 0}^{\left [ \frac{l}{2} \right ]}
\bigoplus_{i \leq j < k \leq l-i} B^i_{(k,j)}.
\]
For $b'\in \CA_+$ such that $E_\CA(b')\in \CA_-$, we see that the action of 
$E_\CA$ is given by \saba 3 or \saba 4 or \saba 5.
We take $b = \sayouA{q}{p}\hwvs \in B_C$ $(0 \leq q \leq p \leq j)$, 
which is an $\CA$-terminal element of an $E_\CA$-string. 
By consulting the proof of (C3) in \S\ref{sec:opeEF}, 
we see that for $b$ the action of $E_\CA$ has following relations,
using sequence depending $b$, $0 \leq n_1 \leq n_2 \leq \cdots \leq n_6 =
\ve{\CA}(b)$. \\
the action of $E_\CA$ on $E_\CA^{n'}b$ $(0 \leq n' < n_1)$ is given 
by \saba 1, if $n_1 \neq 0$,\\
the action of $E_\CA$ on $E_\CA^{n'}b$ $(n_1 \leq n' < n_2)$ is given 
by \saba 3, if $n_1<n_2$,\\
the action of $E_\CA$ on $E_\CA^{n'}b$ $(n_2 \leq n' < n_3)$ is given 
by \saba 4, if $n_2<n_3$,\\
the action of $E_\CA$ on $E_\CA^{n'}b$ $(n_3 \leq n' < n_4)$ is given 
by \saba 5, if $n_3<n_4$,\\
the action of $E_\CA$ on $E_\CA^{n'}b$ $(n_4 \leq n' < n_5)$ is given 
by \saba 2, if $n_4<n_5$,\\
the action of $E_\CA$ on $E_\CA^{n_6}b$ is given by \saba 6.\\
Then  we see that the number $\bar{n}(b) \in  \Z_{\geq 0}$ which satisfy 
following condition is at most one: 
\[
E_\CA^{\bar{n}(b)}b\in \CA_+,
\]
\[
E_\CA^{\bar{n}(b)+1}b \in \CA_-.
\]
If there does not exist such $\bar{n}(b)$, we put
\[
\bar{n}(b)=\ve{\CA}(b)=2l-2i+j-3q.
\]
By Proposition \ref{prop:symEF} and induction (\ref{eqn:indsayou}) 
it is sufficient to prove (\ref{eqn:e1phi})
for $b'=E_\CA^{n}b$ $(0 \leq n \leq \bar{n}(b), b \in B_C)$.

\subsubsection{Actions of $\tea{1}$}
\begin{proposition}
\label{prop:eone}
For $\sayouA{q}{p}\hwvs \in B_C$ $(0 \leq q \leq p \leq j)$,
we have following relations.\vspace{3mm}\\
Put $b_1=\sayouA{q}{p}\hwvs$ $(3q\geq i+2j-\left [\frac{j-i}{2} \right ]$, 
$p+2q>i+2j)$.
We have
\begin{eqnarray}
\tea{1}\left ( E_\CA^{n} b_1 \right ) &=&
E_\CA^{n} \left ( \sayouA{q-1}{p}\hwvl{i}{l-i}{j} \right ),
\label{eqn:eone1aa}\\ 
\tea{1}\left ( E_\CA^{l-i-j+3(j-q)+1} \; b_1 \right ) &=&
E_\CA^{l-i-j+3(j-q)}\left ( \sayouA{q}{p+1}\hwvl{i}{l-i}{j+1} \right ),
\label{eqn:eone1c}\vspace{3mm}\\
\tea{1}\left ( E_\CA^{l-i-j+3(j-q)+2} \; b_1 \right ) &=&
E_\CA^{l-i-j+3(j-q)}\left ( \sayouA{q+1}{p+2}\hwvl{i}{l-i}{j+2} \right ),
\label{eqn:eone1d}\vspace{3mm}\\
\tea{1}\left ( E_\CA^{l-i-j+3(j-q+1)+n'} \; b_1 \right ) &=&
E_\CA^{l-i-j+3(j-q)+n'}\left ( \sayouA{q+2}{p+3}\hwvl{i}{l-i}{j+3} 
\right ), \label{eqn:eone1ee}
\end{eqnarray}
where $0 \leq n \leq l-j-i+3(j-q)$, $0 \leq n' \leq 3(y-1)$, 
$y=\left [ \frac{l-i-j}{3}\right ]$. 
\vspace{5mm}

\noindent Put $b_2=\sayouA{q}{p}\hwvs$
$( p-q \leq \left [ \frac{j-i}{2} \right ]$, $p+2q\leq i+2j )$.
We have
\begin{eqnarray}
\tea{1}\left ( E_\CA^{n} \; b_2 \right ) &=&
E_\CA^{n} \left ( \sayouA{q-1}{p}\hwvl{i}{l-i}{j} \right ),
\label{eqn:eone2aa}\\ 
\tea{1}\left ( E_\CA^{l-2i+[(j-i)/2 ]-n'} \; b_2 \right )
 &=& E_\CA^{l-2i+[(j-i)/2] -3-n'}
\left ( \sayouA{q}{p}\hwvl{i+1}{l-(i+1)}{j-1} \right ),
\label{eqn:eone2cc} 
\end{eqnarray}
where $0 \leq n \leq l-2i+p-q$, $0 \leq n' \leq \left [ \frac{j-i}{2} \right ]$. 
\vspace{5mm}

\noindent Put $b_3=\sayouA{q}{p}\hwvs$  
$( 3q<i+2j-\left [ \frac{j-i}{2} \right ]$, 
$\left [ \frac{j-i}{2} \right ] < p-q )$, then we have
\begin{eqnarray}
\tea{1}\left ( E_\CA^{n} \; b_3 \right ) &=&
E_\CA^{n}\left ( \sayouA{q-1}{p}\hwvl{i}{l-i}{j} \right ),
\label{eqn:eone3b} 
\end{eqnarray}
where $0 \leq n \leq l-2i+\left [ \frac{j-i}{2} \right ]$. \vspace{5mm}
\end{proposition}
\proof
We consider the action of $\tea{1}$ for $\fa{0}{r}\fa{1}{q}\fa{0}{p}
\hwvl{i}{k}{j}$ $(0 \leq p \leq j, p \leq q \leq p+k, 0 \leq r \leq j+q-2p)$.
By Proposition \ref{proposition:kakanA2b}, if $q=p$ we have 
\[
\tea{1}\fa{0}{r}\fa{1}{p}\fa{0}{p}\hwvl{i}{k}{j}=
\fa{0}{r}\fa{1}{p-1}\fa{0}{p}\hwvl{i}{k}{j}
\]
where $0 \leq r \leq j-p$.
If $q>p$,
\[
\tea{1}\fa{0}{r}\fa{1}{p}\fa{0}{p}\hwvl{i}{k}{j}\neq
\fa{0}{r}\fa{1}{p-1}\fa{0}{p}\hwvl{i}{k}{j}
\]
where $q-p\leq r \leq j+q-2p$.

We consider the case of $q=p$ in part I, and the case of $q>p$
in part II. 
\paragraph{Part I}\mbox{} 
We consider $ \fa{0}{r}\sayouA{p}{p}\hwvl ikj$ 
$( 0 \leq p \leq j, 0 \leq r \leq j-p)$.\vspace{2mm}\\
We define $b, b' \in B^i_{(k,j)}$ by 
\begin{eqnarray*}
b&=&\sayouA{p}{p}\hwvl{i}{k}{j},\\
b'&=&\left \{
\begin{array}{ll}
\sayouA{p-1}{p-1}\hwvl{i}{k}{j}&(p>0)\\
0&(p=0)
\end{array}
\right . .
\end{eqnarray*}
Then by Proposition \ref{proposition:kakanA2a} we have
\[
\tea{1}\fa{0}{r} b = \fa{0}{r+1}b',
\]
 (Figure \ref{fig:pIfun}).
We define $m_b$ by
\[
m_b=j-k-3(j-p)=-2j-k+3p,
\]
We consider an $E_\CA$-string. We can express $b$ using an $\CA$-terminal 
element. By Definition \ref{definition:itinokettei}, 
we have
\begin{equation}
b = E_\CA^{l-i-k+(m_b)_+}\left (\sayouA{p-(m_b)_+}{p}\hwvl{i}{l-i}{j-(m_b)_+}
\right )
\label{eqn:bshiki1}
\end{equation}
If $m_b < 0$ the action of $E_\CA$ is given by \saba 1 or \saba 3 of Definition
\ref{definition:itinokettei}, and if $m_b \geq 0$ the action is given by 
\saba 2.\vspace{5mm}

\paragraph{The case of \saba 1:$E_\CA b \in B^i_{(k-1,j)}$.} 
\mbox{}\vspace{3mm}\\
This case is $b=E_\CA^{n'}b_1$ $(0 \leq n' < \mbox{min}\{l-2i, l-i+2j-3q\})$,
$b=E_\CA^{n'}b_2$ $(0 \leq n' < l-2i)$,
$b=E_\CA^{n'}b_3$ $(0 \leq n' < l-2i)$.

\noindent If $p=0$, we have  
\[
\tea{1}E_\CA b=0.
\]
If $p>0$, since  $m_b < 0$, we have $m_{b'} < 0$. 
Then the actions of $E_\CA$ for $b'$  is given by \saba 1. Then we have
\begin{eqnarray*}
E_\CA b&=&\fa{0}{p}\fa{1}{p}\hwvl{i}{k-1}{j},\\
E_\CA b'&=&\fa{0}{p-1}\fa{1}{p-1}\hwvl{i}{k-1}{j},
\end{eqnarray*}
(Figure \ref{fig:pIa1}).
Then we obtain
\begin{equation}
\tea{1} E_\CA \left ( \fa{0}{r} b \right ) = 
E_\CA \left ( \fa{0}{r+1} b'\right ).
\label{eqn:bshiki2}
\end{equation}
We have the case of $(0 \leq n \leq \mbox{min}\{l-2i, l-i+2j-3q\})$ of 
(\ref{eqn:eone1aa}), the case of $(0 \leq n \leq l-2i)$ of (\ref{eqn:eone2aa}), 
and the case of $(0 \leq n \leq l-2i)$ of (\ref{eqn:eone3b}),
\vspace{5mm}

\unitlength 0.8mm
\begin{figure}
\begin{center}
\begin{minipage}{72mm}
\begin{picture}(90,100)
\crystalAtwoA
\crystalAtwoB
\put(30,90){\makebox(0,0){$\hwvl{i}{k}{j}$}}
\put(15,52){\makebox(0,0){$b'$}}
\put(15,40){\makebox(0,0){$b$}}
\put(27,76){\makebox(0,0){\small$p\!\!-\!\!1$}}
\put(24,61){\makebox(0,0){\small$p\!\!-\!\!1$}}
\put(64,55){\makebox(0,0){\small$j\!\!-\!\!p$}}
\end{picture}
\caption{Part I, the location of $b$ and $b'$ in the crystal graph of $B^i_{(k,j)}$.}
\label{fig:pIfun}
\end{minipage}
\makebox[5mm]{}
\begin{minipage}{72mm}
\begin{picture}(90,100)
\crystalAtwoA
\crystalAtwoB
\put(32,90){\makebox(0,0){$\hwvl{i}{k-1}{j}$}}
\put(10,52){\makebox(0,0){$E_\CA b'$}}
\put(10,40){\makebox(0,0){$E_\CA b$}}
\put(27,76){\makebox(0,0){\small$p\!\!-\!\!1$}}
\put(24,61){\makebox(0,0){\small$p\!\!-\!\!1$}}
\put(64,55){\makebox(0,0){\small$j\!\!-\!\!p$}}
\end{picture}
\caption{Part I, the case of \saba 1, the location of $E_\CA b$ and $E_\CA b'$ 
in the crystal graph of $B^i_{(k-1,j)}$.}
\label{fig:pIa1}
\end{minipage}
\end{center}
\end{figure}

\paragraph{The case of \saba 2:$E_\CA b \in B^i_{(k,j+1)}$.}
\mbox{}\vspace{3mm}\\
This case is $b=E_\CA^{n'}b_1$ $(l-i+2j-3q \leq n' < 2l-2i+j-3q, 
l-i+2j-3q < l-2i)$.\\
We denote
\[
E_\CA b=\fa{1}{p+1}\fa{0}{p+1}\hwvl{i}{k}{j+1}.
\]
if  $m_b = 0,1,2$, since $m_{b'} < 0 \leq m_b$, then  the action of 
$E_\CA$ for $b'$ is given by \saba 1. Then we have
\[
E_\CA b'=\fa{1}{p-1}\fa{0}{p-1}\hwvl{i}{k-1}{j}.  
\]
Therefore, we have
\[
\tea{1} E_\CA \left ( \fa{0}{r} b \right ) \neq 
E_\CA \left ( \fa{0}{r+1} b'\right ).
\]
By (\ref{eqn:bshiki1}), we have 
\[
\tea{1}E_\CA\left ( \fa{0}{r} \right ) b = 
E_\CA^{l-i-k}\left ( \fa{0}{r+1}\sayouA{p}{p}\hwvl{i}{l-i}{j+1} \right ).
\]
if  $m_b \geq 3$, then the action of $E_\CA$ for $b'$ is given by \saba 2,
since $0 \leq m_{b'} < m_b$. Thus we have
\[
E_\CA b'=\fa{1}{p}\fa{0}{p}\hwvl{i}{k}{j+1},
\]
(Figure \ref{fig:pIa2}).
Therefore, we have
\begin{equation}
\begin{array}{rcl}
\tea{1}E_\CA \left ( \fa{0}{r} b \right ) &=& 
E_\CA \left ( \fa{0}{r+1} b' \right )\\
&=& E_\CA^{l-i-k+m_b-2} 
\left ( \fa{0}{r+1}\sayouA{p-m_b+2}{p}\hwvl{i}{l-i}{j-m_b+3} \right ).
\end{array}
\label{eqn:bshiki3}
\end{equation}
Therefore we have (\ref{eqn:eone1c}), (\ref{eqn:eone1d}), (\ref{eqn:eone1ee}),
if $l-i+2j-3q<l-2i$.\\
\vspace{5mm}

\paragraph{The case of \saba 3:$E_\CA b\in B^{i+1}_{(i+1,j-1).}$
.} \mbox{}\vspace{3mm}\\
This case is $b=E_\CA^{l-2i}b_1$, $b=E_\CA^{l-2i}b_2$,
$b=E_\CA^{l-2i}b_3$.\\
We denote
\[
E_\CA b=\fa{1}{p+1}\fa{0}{p}\hwvl{i+1}{k+1}{j-1},
\]
Since $m_{b'} < m_b <0$ the actions of $E_\CA$ for $b'$ is given by \saba 3.
Then we have
\[
E_\CA b'=\fa{1}{p}\fa{0}{p-1}\hwvl{i+1}{k+1}{j-1}.
\]
Here by Proposition \ref{proposition:kakanA2b} we have
\[
\tea{1}\fa{0}{r}\sayouA{p+1}{p}\hwvl{i}{k}{j} =
\fa{0}{r+1}\sayouA{p}{p-1}\hwvl{i}{k}{j}
\neq \fa{0}{r}\sayouA{p}{q}\hwvl{i}{k}{j}\;(r>0),
\]
\[
\tea{1}\sayouA{p+1}{p}\hwvl{i}{k}{j} =
\sayouA{p}{p}\hwvl{i}{k}{j}
\neq \tfa{0}\sayouA{p}{p-1}\hwvl{i}{k}{j},
\]
(Figure \ref{fig:pIa3}).
Then we have
\begin{equation}
\tea{1}\fa{0}{r}E_\CA b = \fa{0}{r+1}E_\CA b' \;\; (r > 0),\\
\end{equation}
\begin{equation}
\tea{1}E_\CA b = 
E_\CA^{l-2(i+1)}\left ( \sayouA{p}{p}\hwvl{i+1}{l-i-1}{j-1} \right ).
\end{equation}
Therefore we have the case of $n=l-2i+1$ of (\ref{eqn:eone1aa}), 
(\ref{eqn:eone2aa}), (\ref{eqn:eone3b}).
\vspace{5mm}

\begin{figure}
\begin{center}
\begin{minipage}{72mm}
\begin{picture}(90,100)
\crystalAtwoA
\crystalAtwoB
\put(32,90){\makebox(0,0){$\hwvl{i}{k}{j+1}$}}
\put(10,40){\makebox(0,0){$E_\CA b$}}
\put(30,76){\makebox(0,0){\small$p$}}
\put(27,61){\makebox(0,0){\small$p$}}
\put(64,55){\makebox(0,0){\small$j\!\!-\!\!p$}}
\end{picture}
\caption{Part I, the case of \saba 2, the localtion of $E_\CA b$ in the crystal graph of
$B^i_{(k,j+1)}$.}
\label{fig:pIa2}
\end{minipage}
\makebox[5mm]{}
\begin{minipage}{72mm}
\begin{picture}(90,100)
\crystalAtwoA
\put(20,90){\line(1,-1){45}}
\put(65,45){\line(-1,-1){20.5}}
\put(45,65){\line(-1,-1){31}}
\put(65,45){\circle*{2}}
\put(20,50){\line(1,-1){6}}
\put(26,44){\line(1,-5){1}}
\put(27,39){\line(1,-1){21}}
\put(23,43){\line(1,-1){1}}
\put(24,42){\line(6,-1){6}}
\put(30,41){\line(1,-1){16}}
\put(45.5,25.5){\circle*{2}}
\put(23,43){\circle*{2}}
\put(28,38){\circle*{2}}
\put(48,18){\circle*{2}}
\put(43,23){\circle*{2}}
\put(23,33){\circle*{2}}
\put(38,18){\circle*{2}}
\put(18,38){\line(1,-1){25}}
\put(18,38){\circle*{2}}
\put(43,13){\circle*{2}}
\put(21,31){\line(1,1){7}}
\put(36,16){\line(1,1){7}}
\put(39,9){\line(1,1){9}}
\put(30.5,30.5){\circle*{0.5}}
\put(33,28){\circle*{0.5}}
\put(35.5,25.5){\circle*{0.5}}
\put(35,90){\makebox(0,0){$\hwvl{i+1}{k+1}{j-1}$}}
\put(10,52){\makebox(0,0){$E_\CA b'$}}
\put(10,40){\makebox(0,0){$E_\CA b$}}
\put(26,75){\makebox(0,0){\small$p\!\!-\!\!1$}}
\put(25,60){\makebox(0,0){\small$p$}}
\put(68,55){\makebox(0,0){\small$j\!\!-\!\!p\!\!-\!\!1$}}
\put(26,22){\makebox(0,0){\small$j\!\!-\!\!p$}}
\end{picture}
\caption{Part I, the case of \saba 3, the location of $E_\CA b$ and $E_\CA b'$ in 
the crystal graph $B^{i+1}_{(k+1,j-1)}$.}
\label{fig:pIa3}
\end{minipage}
\end{center}
\end{figure}

\paragraph{Part II.}\mbox{} 
We consider $\tea{1}E_\CA\fa{0}{r}\sayouA{p+s}{p}\hwvl{i}{k}{j}$ 
$(i\leq k,j \leq l-i,\; k-i \geq s>0)$ \\
If  $i<k,j \leq l-i$, we use induction (\ref{eqn:indsayou}).    
Otherwise, using involution, we have 
\[
\syme{\hwvl{i}{j}{i}} = \fa{0}{i}\sayouA{i+j}{j}\hwvl{i}{i}{j}.
\]
Then it is sufficient to consider following element:
\[
\tea{1}E_\CA\fa{0}{r}\sayouA{p+s}{p}\hwvl{i}{i}{j}\;\;
(i < j \leq l-i,\; 0 < s \leq k-i ).
\]
\vspace{3mm}\\
We denote $b, b', b'' \in B^i_{(k,j)}$  by
\begin{eqnarray*}
b&=&\sayouA{p+s}{p}\hwvl{i}{k}{j},\\
b'&=&\sayouA{p+s+1}{p-1}\hwvl{i}{k}{j},\\
b''&=&\tilde{e}_1^{\cal A}b.
\end{eqnarray*}
Then we have
\begin{eqnarray*}
\tea{1}\fa{0}{r}b &=& \fa{0}{r} b''\:\:(r<s),\\
\tea{1}\fa{0}{r}b &=& \fa{0}{r+1}b'\:\:(r\geq s),
\end{eqnarray*}
(Figure \ref{fig:pIIfun}).
We rewrite $m_b$
\[
m_b=(j+s)-(k-s)+s-3((j+s)-p)=-2j-k+3p.
\]
Then we denote $b$ using the $\CA$-terminal element in an $E_\CA$-string.  
\begin{equation}
\hspace{-5mm}
b = \left \{
\begin{array}{l}
E_\CA^{s+l-2(k-s+1)}\sayouA{p}{p}\hwvl{k-s+1}{l-(k-s+1)}{j+s}
\:\mbox{ if } m_b<0,\vspace{3mm}\\
E_\CA^{s+l-2(k-s+1)+m_b}\sayouA{p-m_b}{p-m_b}
\hwvl{k-s+1}{l-(k-s+1)}{j+s-m_b},\:\mbox{ if } m_b \geq 0.
\end{array}
\right .
\label{eqn:bshiki4}
\end{equation}
\vspace{5mm}

\noindent
\paragraph{The case of \saba 2: $E_\CA b \in B^i_{(i,j+1)}$.}
\mbox{} \vspace{3mm}\\ 
This case is $b=E_\CA^{n'}b_1$ $(l-2i \leq l-i+2j-3q \leq n' < 2l-2i+j-3q)$.\\
We denote
\[
E_\CA b=\fa{1}{p+s+1}\fa{0}{p+1}\hwvl{i}{i}{j+1}.
\] 
if  $r < s$, since $m_b = m_{b''}$, then the action of $E_\CA$ for $b''$ 
is given by \saba 2. Then we have
\[
E_\CA b''=\fa{1}{p+s}\fa{0}{p+1}\hwvl{i}{i}{j+1}.
\]
Therefore we have,
\[
\tea{1}E_\CA\left ( \fa{0}{r} b \right ) = 
E_\CA\left ( \fa{0}{r} b'' \right ). 
\]
If $r \geq s$, $m_b=0,1,2$, since $m_{b'} < 0$, then 
the action of $E_\CA$ for $b'$ is given by \saba 3. Then we have
\[
E_\CA b'=\fa{1}{p+s}\fa{0}{p-1}\hwvl{i+1}{i+1}{j-1}.
\]
therefore we have
\begin{equation}
\tea{1}E_\CA\left ( \fa{0}{r} b \right ) 
= E_\CA^{s+l-2(i-s+1)} \left ( \fa{0}{r+1} \sayouA{p}{p}
\hwvl{k-s+1}{l-(i-s+1)}{j+s+1} \right ).
\label{eqn:bshiki5}
\end{equation}
If $r \leq s$, $m_b \geq 3$, than the actions of 
$E_\CA$ for $b'$ is given by \saba 2, since $m_{b'} \geq 0$. Then we have,
\[
E_\CA b'=\fa{1}{p+s}\fa{0}{p}\hwvl{i}{i}{j+1},
\]
(Figure \ref{fig:pIIa2}).
Thus we have 
\begin{equation}
\begin{array}{l}
\tea{1}E_\CA\left ( \fa{0}{r} b \right ) =
E_\CA\left ( \fa{0}{r+1} b'\right )\vspace{2mm}\\
\makebox[10mm]{}
=E_\CA^{s+l-2(i-s+1)+m-2} \left ( \fa{0}{r+1}\sayouA{p-m+2}{p-m+2}
\hwvl{i-s+1}{l-(i-s+1)}{j-m+2+s} \right )
\end{array}
\label{eqn:bshiki6}
\end{equation}
Therefore we have the case of (\ref{eqn:eone1c}), (\ref{eqn:eone1d}), 
(\ref{eqn:eone1ee}), if $l-2i \leq l-i+2j-3q$. \\
\vspace{5mm}
\begin{figure}
\begin{center}
\begin{minipage}{72mm}
\begin{picture}(90,100)
\crystalAtwoC
\put(30,90){\makebox(0,0){$\hwvl{i}{i}{j}$}}
\put(2,48){\makebox(0,0){$b'$}}
\put(-2,35){\makebox(0,0){$b$}}
\put(22,80){\makebox(0,0){\small$p\!\!-\!\!1$}}
\put(22,70){\makebox(0,0){\small$p\!\!-\!\!1$}}
\put(34,56){\makebox(0,0){\small$p$}}
\put(10,56){\makebox(0,0){\small$s$}}
\put(5,23){\makebox(0,0){\small$s\!\!-\!\!1$}}
\put(18,38){\vector(-1,0){8}}
\put(22,38){\makebox(0,0){$b''$}}
\end{picture}
\caption{Part II, the location of $b$, $b'$ and $b''$ in the crystal graph $B^i_{(i,j)}$.} 
\label{fig:pIIfun}
\end{minipage}
\makebox[5mm]{}
\begin{minipage}{72mm}
\begin{picture}(90,100)
\crystalAtwoC
\put(32,90){\makebox(0,0){$\hwvl{i}{i}{j+1}$}}
\put(-5,35){\makebox(0,0){$E_\CA b$}}
\put(18,38){\vector(-1,0){8}}
\put(25,38){\makebox(0,0){$E_\CA b''$}}
\put(25,80){\makebox(0,0){\small$p$}}
\put(25,70){\makebox(0,0){\small$p$}}
\put(35,56){\makebox(0,0){\small$p\!\!+\!\!1$}}
\put(10,56){\makebox(0,0){\small$s$}}
\put(5,23){\makebox(0,0){\small$s\!\!-\!\!1$}}
\end{picture}
\caption{Part II, the case of \saba 2, the location of $E_\CA b$ and $E_\CA b''$ in the
crystal graph of $B^i_{(i,j+1)}$.}
\label{fig:pIIa2}
\end{minipage}
\end{center}
\end{figure}

\paragraph{The case of \saba 3:$E_\CA b \in B^{i+1}_{(i+1,j-1)}$.}
\mbox{}\vspace{3mm}\\
This case is $b=E_\CA^{n'}b_1$ $(l-2i < n' < l-i+2j-3q)$,
$b=E_\CA^{n'}b_2$ $(l-2i < n' < \left [ \frac{j-i}{2}\right ])$,
$b=E_\CA^{n'}b_3$ $(l-2i < n' < \left [ \frac{j-i}{2}\right ])$.
(\ref{eqn:eone3b}).\\
Since $m_{b'} < m_{b''} = m_b < 0$, the action of $E_\CA$ for $b'$, $b''$
is given by \saba 3. Then we have
\begin{eqnarray*}
E_\CA b&=&\fa{1}{p+s+1}\fa{0}{p}\hwvl{i+1}{k+1}{j-1},\\
E_\CA b'&=&\fa{1}{p+s}\fa{0}{p-1}\hwvl{i+1}{k+1}{j-1},\\
E_\CA b''&=&\fa{1}{p+s}\fa{0}{p}\hwvl{i+1}{k+1}{j-1},
\end{eqnarray*}
(Figure \ref{fig:pIIa3}). 
With Proposition \ref{proposition:kakanA2b}, we have
\begin{equation}
\begin{array}{ll}
\tea{1}E_\CA\left ( \fa{0}{r} b \right ) =
E_\CA\left ( \fa{0}{r} b'' \right )& (r<s+1),\vspace{2mm}\\
\tea{1}E_\CA\left ( \fa{0}{r} b \right ) =
E_\CA\left ( \fa{0}{r+1} b' \right ) & (r\geq s+1).
\end{array}
\label{eqn:bshiki7}
\end{equation}
Therefore we have the case of $(l-2i < n \leq l-i+2j-3q)$ of (\ref{eqn:eone1aa}),
the case of $(l-2i< n \leq \left [ \frac{j-i}{2}\right ])$ of 
(\ref{eqn:eone2aa}) and (\ref{eqn:eone2cc}),
\vspace{3mm}\\
By (\ref{eqn:bshiki1}) -- (\ref{eqn:bshiki7}), we have Proposition 
\ref{prop:eone}. \hfill$\Box$
\vspace{5mm}
\begin{figure}
\begin{center}
\begin{minipage}{72mm}
\begin{picture}(90,100)
\crystalAtwoC
\put(32,90){\makebox(0,0){$\hwvl{i+1}{i}{j+1}$}}
\put(-1,48){\makebox(0,0){$E_\CA b'$}}
\put(-5,35){\makebox(0,0){$E_\CA b$}}
\put(22,80){\makebox(0,0){\small$p\!\!-\!\!1$}}
\put(22,70){\makebox(0,0){\small$p\!\!-\!\!1$}}
\put(33,57){\makebox(0,0){\small$p$}}
\put(8,56){\makebox(0,0){\small$s\!\!+\!\!1$}}
\put(8,24){\makebox(0,0){\small$s$}}
\put(18,38){\vector(-1,0){8}}
\put(25,38){\makebox(0,0){$E_\CA b''$}}
\end{picture}
\caption{Part II, the case of \saba 3, the location of $E_\CA b$, $E_\CA b'$ and 
$E_\CA b''$ in the crystal graph of $B^{i+1}_{(i,j+1)}$.}
\label{fig:pIIa3}
\end{minipage}
\end{center}
\end{figure}

\subsubsection{Actions of $\teg{1}$}
\begin{proposition} \label{proposition:kakanlong} \rm
For $b \in \CG$, $n \in \Z_{>0}$ if $\tfg{1}\left ( \teg{2} \right )^n b = 
\left ( \teg{2} \right )^n \tfg{1} b \neq 0$, we have
\[
\tfg{1}\left (\teg{2} \right )^kb =
\left ( \teg{2}\right )^k\tfg{1}b \; 
\; ( 0 \leq k \leq n ).
\]
\end{proposition}
\proof
We prove the proposition by the induction on $k$.
We assume 
\[
\tf{1}\left( \te{2} \right)^{k'}b=\left( \te{2} \right)^{k'}\tf{1}b\;\;(0 \leq k' < k ),
\]
then we prove 
\begin{equation}
\tf{1}\te{2} \left ( \te{2} \right )^{k-1} b =\te{2}\tf{1} \left ( \te{2} \right )^{k-1}b. 
\label{eqn:kakantf1te2}
\end{equation}
We write $\te{2}^{k-1} b = b_1 \otimes \cdots \otimes b_l$. We have
\begin{equation}
\tf{1}\left (\te{2} \right )^{k-1}b = b_1 \otimes \cdots \otimes \tf{1}b_{k_1}\otimes
\cdots \otimes b_l,
\label{eqn:sayouf1k1}
\end{equation}
\begin{equation}
\te{2}\left (\te{2} \right )^{k-1}b = b_1 \otimes \cdots \otimes \te{1}b_{k_2}\otimes
\cdots \otimes b_l,
\label{eqn:sayoue2k2}
\end{equation}
where $1 \leq k_1, k_2 \leq l$.\vspace{3mm}

\noindent The case of  $k_1 > k_2$.
Similar to Proposition \ref{lemma:phizogen}, 
by Proposition \ref{proposition:politensor}
we see that the action of $f_1$ increase the number of $u_+$ of $\mbox{\rm Red}_2(b)$.
Then the operator $\tf{1}$ does not influence the action of $\te{2}$.
In a similar way, we see that the operator $\te{2}$ does not influence the action of 
$\tf{1}$. 
Therfore we have (\ref{eqn:kakantf1te2}).

\noindent The case of $k_1 \leq k_2$.
We assume 
\begin{equation}
\tf{1}\te{2} \left ( \te{2} \right )^{k-1} b \neq
\te{2}\tf{1} \left ( \te{2} \right )^{k-1}b. 
\label{eqn:asuumeneq}
\end{equation}
By assumption, at most one of following relations are satisfied:
\begin{equation}
\tf{1}\te{2}\left ( \te{2}^{k-1} \right )b= \left \{
\begin{array}{l}
b_1\otimes \cdots \otimes \tf{1}\te{2}b_{k_2} \otimes \cdots \otimes b_l\\
b_1\otimes \cdots \otimes \te{2}b_{k_2} \otimes \cdots \otimes \tf{1}b_{k'_1}
\otimes \cdots \otimes b_l  
\end{array} \right .,
\label{eqn:sayouf1e2} 
\end{equation}
\begin{equation}
\te{2}\tf{1}\left ( \te{2}^{k-1} \right ) b= \left \{
\begin{array}{l}
b_1\otimes \cdots \otimes \te{2}\tf{1}b_{k_1} \otimes \cdots \otimes b_l\\
b_1\otimes \cdots \otimes \te{2}b_{k'_2} \otimes \cdots \otimes \tf{1}b_{k_1}
\otimes \cdots \otimes b_l
\end{array} \right .,
\label{eqn:sayoue2f1}
\end{equation}
Since $\tf{1}u_+=u_-$ and $\uu{1}(b)=u_-^{\ve{1}(b)}u_+^{\vp{1}(b)}$,
we have $k_1 < k'_1$.
In a similar way we have $k'_2 < k_2$.
If (\ref{eqn:sayouf1e2}) is satisfied, 
we see that $k_1$-th element of $\te{2}^{k'}\tf{1}\te{2}^kb$ and 
$\tf{1}\te{2}^{k'}\te{2}^kb$ $(1 \leq k' \leq n-k)$ are different.
Then we have 
\begin{equation}
\te{2}^n\tf{1}b \neq \tf{1}\te{2}^n b.
\label{eqn:mujun}
\end{equation}
This is contradiction.
In a similar way, if (\ref{eqn:sayoue2f1}) is satisfied, we have contradiction.
Thus we have (\ref{eqn:kakantf1te2}). \hfill$\Box$

We prove from (\ref{eqn:eone1aa}) to (\ref{eqn:eone3b}) exchanging $\tea{1}$
with $\teg{1}$. For example, we denote following equations exchanging
$\tea{1}$ with $\teg{1}$ from (\ref{eqn:eone1aa}) to (\ref{eqn:eone1ee}).
We put $b_1=\fa{1}{q}\fa{0}{p}\hwvs \in B_C$ 
$( 3q\geq i+2j-\left [\frac{j-i}{2} \right ]$, $p+2q>i+2j )$.\vspace{3mm}\\
We are going to show
\begin{eqnarray}
\makebox[-8mm]{}\teg{1}\left ( \Phi\left ( E_\CA^{n} \; b_1 \right ) \right )&=&
\Phi \left ( E_\CA^{n} \left ( \sayouA{q-1}{p}\hwvl{i}{l-i}{j} \right ) 
\right ),\label{eqn:eone4aa}\\ 
\makebox[-8mm]{}
\teg{1}\left ( \Phi \left ( E_\CA^{l-i-j+3(j-q)+1} \; b_1 \right ) \right ) &=&
\Phi \left ( E_\CA^{l-i-j+3(j-q)}
\left ( \sayouA{q}{p+1}\hwvl{i}{l-i}{j+1} \right ) \right ),
\label{eqn:eone4c}\vspace{3mm}\\
\makebox[-8mm]{}
\teg{1}\left ( \Phi \left ( E_\CA^{l-i-j+3(j-q)+2} \; b_1 \right ) \right ) &=&
\Phi \left ( E_\CA^{l-i-j+3(j-q)} \left ( 
\sayouA{q+1}{p+2}\hwvl{i}{l-i}{j+2} \right ) \right ),\label{eqn:eone4d}
\vspace{3mm}\\
\makebox[-8mm]{}
\teg{1}\left ( \Phi \left ( E_\CA^{l-i-j+3(j-q+1)+n'} \; b_1 \right ) \right ) 
&=& \Phi \left ( E_\CA^{l-i-j+3(j-q)+n'} \left ( \sayouA{q+2}{p+3}
\hwvl{i}{l-i}{j+3} \right ) \right ), \label{eqn:eone4ee} 
\end{eqnarray}
where $0 < n < l-j-i+3(j-q)$, $0 < n' < 3(y-1)$. 
\vspace{5mm}

Here we prove the case of $l-i-j \equiv 0 \pmod{3}$,  we can prove other 
cases similarly. \vspace{3mm}

By Proposition \ref{fact:f0pR}, we have
\[
\Phi(b_1) = \BB{\bb{1}{i}\bb{6}{j-q}\bb{6}{y}\bB{\CC}{y+q+p-j-i}
\bb{2}{y-1}\bb{-1}{j+i-p}}.
\]
\paragraph{The case of (\ref{eqn:eone4aa}), $n=0$.}
By (E1), we have  
\[
\teg{1}\left ( \Phi (b_1) \right ) = \Phi \left ( \sayouA{q-1}{p}\hwvs \right ).
\]
\paragraph{The case of (\ref{eqn:eone4aa}), $n=l-i-j+3(j-q)$.}
Put $b_4 = E_\CA^{l-j-i+3(j-q)}b$. We see 
\[
\Phi(b_4) =  \BB{\bb{1}{i}\bb{2}{j-p}\bb{2}{y}\bB{\CC}{y+q+p-j-i}
\bb{-2}{y}\bb{-1}{j+i-p}}.
\]
Using Definition \ref{def:uui}
\[
\begin{array}{lcl}
\uu{1}(b_4)&=& u_-^{j-p+i} u_+^y u_-^{y+q+p-j-i} u_+^{y+q+p-j-i} u_-^{y+j-p} 
u_+^i,\\
\mbox{Red}_1(b_4) &=& u_-^{j-p+i} u_-^{q+p-j-i} u_+^{2q+p-2j-i} u_+^i. 
\vspace{2mm}
\end{array}
\]
By Proposition \ref{proposition:vplus1}, we have
\begin{eqnarray*}
\teg{1}\left ( \Phi(b_4) \right )&=& \BB{\bb{1}{i}\bb{2}{j-p}\bb{2}{y}
\teg{1}\bB{\CC}{y+q+p-j-i}\bb{-2}{y}\bb{-1}{j+i-p}}\\
&=&\BB{\bb{1}{i}\bb{2}{j-p}\bb{2}{y}\bb{6}{}\bB{\CC}{y+q+p-j-i-1}
\bb{-2}{y}\bb{-1}{j+i-p}}.
\end{eqnarray*}
Put $b_4^{\prime}= E_\CA^{n}\left ( \sayouA{q-1}{p}\hwvs \right )$.
We see
\[
\Phi(b_4^{\prime}) = \BB{\bb{1}{i}\bb{2}{j-p}\bb{2}{y}\bb{6}{}
\bB{\CC}{y+p+q-j-i-1}\bb{-2}{y}\bb{-1}{j+i-p}}.
\]
Thus we have
\[
\teg{1}\left ( \Phi \left ( E_\CA^{n} b_1 \right ) \right ) = 
\Phi \left ( E_\CA^{n}\left ( \sayouA{q-1}{p}\hwvl{i}{l-i}{j} \right ) 
\right ).
\]
\paragraph{The case of (\ref{eqn:eone4c}).}
We set $b_5 = E_\CA^{l-j-i+3(j-q)+1}b$.
\[
\Phi(b_5) = \BB{\bb{1}{i}\bb{2}{j-p}\bb{2}{y}\bB{\CC}{y+q+p-j-i}
\bb{-3}{}\bb{-2}{y-1}\bb{-1}{j+i-p}}.
\]
We consider $\teg{1}(\Phi(b_5))$ with $u_+$ and $u_-$.
\[
\begin{array}{lcl}
\uu{1}(\Phi(b_5))&=& 
u_-^{j+i-p} u_+^{y-1} u_0 u_-^{y+q+p-j-i} u_+^{y+q+p-j-i} u_-^{y+j-p} u_+^i, 
\vspace{3mm}\\
\mbox{Red}_1\left (\Phi(b_5) \right )&=& 
u_-^{j+i-p} u_-^{p+q-j-i+1} u_+^{2q+p-2j-i} u_+^i.
\end{array}
\]
Then we have
\begin{eqnarray*}
\teg{1}\left ( \Phi(b_5) \right )&=& \BB{\bb{1}{i}\bb{2}{j-p}\bb{2}{y}
\teg{1}\bB{\CC}{y+q+p-j-i}\bb{-3}{}\bb{-2}{y-1}\bb{-1}{j+i-p}}\\
&=&\BB{\bb{1}{i}\bb{2}{j-p}\bb{2}{y}\bb{6}{}\bB{\CC}{y+q+p-j-i-1}
\bb{-3}{}\bb{-2}{y-1}\bb{-1}{j+i-p}}.
\end{eqnarray*}
We put $b_5^{\prime} = E_\CA^{l-i-j+3(j-q)}\left ( \sayouA{q}{p+1}
\hwvl{i}{l-i}{j+1} \right )$. Then we have 
\[
\Phi(b_5^{\prime}) = \BB{\bb{1}{i}\bb{2}{j-q}\bb{2}{y}\bb{6}{}
\bB{\CC}{y+p+q-j-i-1}\bb{-3}{}\bb{-2}{y-1}\bb{-1}{j+i-p}}.
\]
Thus we have
\[
\teg{1}\left ( \Phi \left ( E_\CA^{l-i-j+3(j-q)+1} b_1 \right ) \right )
= E_\CA^{l-i-j+3(j-q)}\left ( \sayouA{q}{p+1} \hwvl{i}{l-i}{j+1} \right ).
\]
Similarly, in case of (\ref{eqn:eone4d}),  $n'=0$ and $n'=2(l-i-j)+3(j-q)$ of
(\ref{eqn:eone4ee}), we have
\[
\teg{1} \left ( \Phi \left ( E_\CA^{l-i-j+3(j-q)+2} b_1 \right ) \right ) =
\Phi \left ( E_\CA^{l-i-j+3(j-q)} \left ( \sayouA{q+1}{p+2}
\hwvl{i}{l-i}{j+2} \right ) \right ),
\]
\[
\teg{1}\left ( \Phi \left ( E_\CA^{l-i-j+3(j-p)+3} b_1 \right ) \right ) =
\Phi \left ( E_\CA^{l-i-j+3(j-q)}
\left ( \sayouA{q+2}{p+3}\hwvl{i}{l-i}{j+3} \right ) \right ),
\]
\[
\teg{1}\left ( \Phi \left ( E_\CA^{2(l-i-j)+3(j-p)} b_1 \right ) \right )= 
\Phi \left ( E_\CA^{2(l-i-j)+3(j-q-1)}
\left ( \sayouA{q+2}{p+3} \hwvl{i}{l-i}{j+3} \right ) \right ).
\]
By Proposition \ref{proposition:kakanlong}, we have $\teg{1}\Phi=\Phi\tea{1}$ for
$E_\CA^nb_1$ $(0 \leq n \leq \ve{\CA}(b_1)=2(l-i-j)+3(j-p))$.
\vspace{5mm}

Similarly, we can prove $\teg{1}\Phi=\Phi\tea{1}$ for any $b \in \CA$.
\vspace{5mm}

%% file: minimal.tex
\subsection{Selection of minimal}
\label{section:minimal}
\subsubsection{Selection of minimal elements in $B^{G_2}(l\Lambda_1)$}
In this section, we only consider crystal $B^l$, so we denote $\te{1}$, 
$\te{2}$ instead of $\teg{1}$, $\teg{2}$ respectively for simplicity. 
By definition \ref{definition:perfectness}, minimal elements are 
$b \in B^l$ such that
$\left \langle c, \vp{}(b) \right \rangle=l$. 
At first we search for the element $b \in B(l\Lambda)$
such that $\left \langle c, \vp{1}(b)\Lambda_1+\vp{1}(b)\Lambda_2 \right 
\rangle = l$. Next we verify that there does not exist $b\in B(l\Lambda_1)$ 
such that $\langle c, \vp{1}(b)\Lambda_1+\vp{2}(b)\Lambda_2 \rangle < l$.
Then we verify that for $b\in B^l$ such that 
$\langle c, \vp{1}(b)\Lambda_1+\vp{2}(b)\Lambda_2\rangle=l$,
we have $\vp{0}(b)=0$
\begin{proposition}
\label{proposition:saishouchi}
For $b \in B^{G_2}(l\Lambda_1)$ such that $\ve{2}(b)=0$, we have 
\[
\min \{ 2 \vp{1}(\tf{2}^{(k)}(b))+\vp{2}(\tf{2}^{(k)}(b))\:|
\:k=0,\ldots,\vp{2}(b) \} \geq l-\left [\frac{1}{2}wt_0(b) \right ],
\]
\end{proposition}

\vspace{5mm}

\proof
We consider 
\begin{equation}
\min\{ 2\varphi_1(b')+\varphi_2(b')\:|\: b'=\tilde{f}_2^{(k)}(b),
\:k=0,\ldots,\varphi_2(b)\}.
\label{eqn:minif2}
\end{equation}
By Proposition \ref{lemma:phizogen}, we have
\[
2\vp{1}(\tf{2}^{(k+1)}(b))+\vp{2}(\tf{2}^{(k+1)}(b))-
(2\vp{1}(\tf{2}^{(k)}(b))+\vp{2}(\tf{2}^{(k)}(b)))=\left \{
\begin{array}{ll}
1& (k < k')\\
-1& (k \leq k')
\end{array} \right ..
\]
We set $b' = \tf{2}^{\vp{2}(b)}b$. Then (\ref{eqn:minif2}) is 
\[
2 \vp{1}(b)+\vp{2}(b)-\frac{\vp{2}(b)-\left (2\vp{1}(b')-\left (2\vp{1}(b)
+\vp{2}(b)\right )\right )}{2}= \vp{1}(b)+\vp{1}(b').
\]
Then we prove
\[
\vp{1}(b)+\vp{1}(b')-l+\left [ \frac{1}{2}\wt_0(b) \right ]\geq 0.
\]
We recall that 
$x_i(b)= \sharp \left \{ \left . \vminiD{b_k} = \vminiD{i} \:\right |
\: k = 1, \ldots ,n \right \}$ $(i= 1,\ldots,6,0_1,0_2)$, 
$\overline{x}_i(b)=\sharp \left \{ \left . \vminiD{b_k} = \vminiD{\overline{i}}
\:\right |\: k = 1, \ldots ,n \right \}$ $(i= 1,\ldots,6)$.  
We consider $b$ such that $\ve{2}(b)=0$. Such an element $b$ satisfy 
following conditions:
\renewcommand{\labelenumi}{\rm (\arabic{enumi}) }
\begin{enumerate}
\item $0 \leq x_1(b)+\overline{x}_1(b)\leq l$.
This is because $\tf{2}\left ( \vmini{1} \right )=0$, 
$\te{2}\left ( \vmini{1} \right ) = 0$, 
$\tf{2}\left ( \vmini{-1} \right ) = 0$, 
$\te{2}\left ( \vmini{-1} \right ) = 0$. 
\item $x_{0_2}(b)=0 \mbox{ or } 1$.
Since $\tf{2}\left ( \vmini{8} \right ) = 0$,
$\te{2}\left ( \vmini{8} \right ) = 0$ and by Proposition 
\ref{proposition:gtwocrystal} $x_{0_2}(b)<2$.
\item $\overline{x}_2(b)=\overline{x}_3(b)=\overline{x}_4(b)=
\overline{x}_5(b)=\overline{x}_{0_1}(b)=0$. Because if  
$\overline{x}_2(b) >0$ or  $\overline{x}_3(b) >0$ or $\overline{x}_4(b) >0$ or 
$\overline{x}_5(b) >0$ or  $\overline{x}_{0_1}(b) >0$ then $\ve{2}(b)>0$.
\item $x_3(b)+x_4(b)+x_5(b) = 0 \mbox{ or } 1$. By  Proposition 
\ref{proposition:gtwocrystal}. 
\item $x_6(b) \leq  \overline{x}_6(b)$. In particular, if $x_3(b)=1$ or 
$x_4(b)=1$ then $x_6(b) < \overline{x}_6(b)$ and if $x_5(b)$ then 
$x_6(b)=0$. 
\item $\displaystyle\sum_{i=1}^6 \left (x_i(b)+\overline{x}_i(b)\right )+
x_{0_1}(b)+x_{0_2}(b)=l$. 
\end{enumerate}
We start with the case of $x_5(b)=0$.\\
By calculation, we have  
$x_1(b')=x_1(b)$, $x_6(b')=x_2(b)+x_3(b)+x_4(b)+x_6(b)$, 
$x_{0_2}(b')=x_{0_2}(b)$, $\overline{x}_6(b')=x_6(b)$, 
$\overline{x}_4(b')=x_4(b)$, $\overline{x}_3(b')=x_3(b)$, 
$\overline{x}_2(b')=\overline{x}_6(b)-x_{0_2}(b)-x_6(b)-x_4(b)-x_3(b)$, 
$\overline{x}_1(b')=\overline{x}_1(b)$.
We put $\delta =\wt_0(b)\:\mod\:2$. Thus we have 
\begin{eqnarray*}
-wt_0(b) &=& 2x_1(b)+x_2(b)+x_3(b)+x_4(b)+x_6(b)
-\overline{x}_6(b)-2\overline{x}_1(b), \\
l &=& x_1(b)+x_2(b)+x_3(b)+x_4(b)+x_6(b)+x_{0_2}(b)
+\overline{x}_6(b)+\overline{x}_1(b),  \\
l-\left [\frac{1}{2}wt_0(b)\right ] &=& 2x_1(b)+\frac{3}{2}x_2(b)
+\frac{1}{2}x_3(b)+\frac{1}{2}x_4(b)+\frac{3}{2}x_6(b)
+\frac{1}{2}\overline{x}_6(b)+\frac{\delta}{2}, \\
\varphi_1(b) &=& x_1(b)+(x_{0_2}(b)+2x_6(b)+x_4(b)-x_2(b))_+, \\
\varphi_1(b') &=& x_1(b')+x_{0_2}(b')+2x_6(b')+(\overline{x}_2(b')
-\overline{x}_4(b')-\overline{x}_6(b')-x_{0_2}(b'))_+, \\
&=& x_1(b)+2x_2(b)+2x_3(b)+2x_4(b)+2x_6(b)+x_{0_2}(b) \\
&&+(\overline{x}_6(b)-x_{0_2}(b)-3x_6(b)-2x_4(b)-x_3(b))_+.
\end{eqnarray*}
\begin{eqnarray}
\varphi_1(b)+\varphi_1(b')-l+\left[ \frac{1}{2}wt_0(b)\right ] &=&
\frac{1}{2} \left | x_{0_2}(b)+2x_6(b)+x_4(b)-x_2(b) \right |\nonumber \\
&&+\frac{1}{2} \left | \overline{x}_6(b)-x_{0_2}(b)-3x_6(b)-2x_4(b)-x_3(b)
\right | -\frac{\delta}{2}\label{eqn:minimaljouken1}\\
&\geq& 0 \mbox{ (integer greater than $-\frac{1}{2}$ )}.\nonumber
\end{eqnarray}
The case of $x_5(b)=1$.
By Proposition \ref{proposition:gtwocrystal}, $x_3(b)=x_4(b)=x_6(b)=
x_{0_2}(b)=0$. Similar to the case of $x_5(b)=0$, we have 
\begin{equation}
\varphi_1(b)+\varphi_1(b')-l+\left[ \frac{1}{2}wt_0(b)\right ] =
\frac{1}{2}x_2(b)+\frac{1}{2}\overline{x}_6(b) -\frac{\delta}{2}\geq 0 
\mbox{ ( integer greater than $-\frac{1}{2}$ )}
\label{eqn:minimaljouken2}
\end{equation}
$\Box$

A minimal element is an element $b\in B^l$ such that $\langle c, \vp{}(b)
\rangle =l$. By $\langle c, \vp{}(b) \rangle = \vp{0}(b) + 2\vp{1}(b) 
+ \vp{2}(b)$ and Proposition \ref{proposition:saishouchi},  
if $\wt_0(b)>0$ we have $\langle c, \vp{}(b) \rangle > l$.
Therefore for a minimal element $b$ we must have $\wt_0(b)=0$.
\begin{lemma} we have
\[
\min\{\langle c,\varphi(b) \rangle \mid b \in B(n\Lambda_1) \}=l .
\]
\end{lemma}
\proof
By (\ref{eqn:minimaljouken1}) and (\ref{eqn:minimaljouken2}), the elements $b$ 
such that $\ve{2}(b)=0$, $\vp{1}(b)+\vp{1}(\tf{2}^{\vp{2}(b)}(b))=l$ 
satisfy following conditions.\\
If $x_5(b)=0$,
\begin{eqnarray*}
x_2(b)&=&2x_6(b)+x_{0_2}(b)+x_4(b),\\
\overline{x}_2(b)&=&3x_6(b)+x_{0_2}(b)+2x_4(b)+x_3(b).
\end{eqnarray*}
If $x_5(b)=1$,
\begin{eqnarray*}
x_2(b)&=&0,\\
x_6(b)&=&0.
\end{eqnarray*}
Then we can express $b$ such that $\ve{2}(b)=0$, 
$\vp{1}(b)+\vp{1}(\tf{2}^{\vp{2}(b)}(b))=l$ as: 
\[
b\,=\,
\left \{ 
\begin{array}{ll}
\BB{\bb{1}{m}\bb{2}{n}\bB{\CC}{n}\bb{-6}{n}\bb{-1}{m}}&(l=2m+3n),
\vspace{3mm}\\
\BB{\bb{1}{m}\bb{2}{n}\bb{4}{}\bB{\CC}{n-1}\bb{-6}{n+1}\bb{-1}{m}}&
(n>0,\,l=2m+3n+1),
\vspace{3mm}\\
\BB{\bb{1}{m}\bb{5}{}\bb{-1}{m}}&(l=2m+1), \vspace{3mm}\\
\BB{\bb{1}{m}\bb{2}{n}\bb{3}{}\bB{\CC}{n}\bb{-6}{n+1}\bb{-1}{m}}&
(l=2m+3n+2).
\end{array} \right .
\]
where $m,n \in \Z_{\geq 0}$. We put
\[
\tilde{b}_{(m)}=\tf{2}^{3n} b = \left \{
\begin{array}{ll}
\BB{\bb{1}{m}\bb{2}{n}\bB{\CC}{n}\bb{-2}{n}\bb{-1}{m}}&(l=2m+3n),
\vspace{3mm}\\
\BB{\bb{1}{m}\bb{2}{n}\bb{4}{}\bB{\CC}{n-1}\bb{-4}{}\bb{-6}{n}\bb{-1}{m}}
&(n>0,l=2m+3n+1),\vspace{3mm}\\
\BB{\bb{1}{m}\bb{7}{}\bb{-1}{m}}&(l=2m+1),
\vspace{3mm}\\
\BB{\bb{1}{m}\bb{2}{n}\bb{3}{}\bB{\CC}{n}\bb{-3}{}\bb{-2}{n}\bb{-1}{m}}&
(l=2m+3n+2),
\end{array} \right .
\]
where $m \in {\bf Z}_\geq 0 \mbox{ such taht } l-2m \geq 0$. 
We prove $\langle c, 2\vp{1}(\tilde{b}_{(m)})\Lambda_1 + \vp{2}(\tilde{b}_{(m)})
\Lambda_2 \rangle=l$ $(k=1,2,3)$.
We prove the case of $l=2m+3n$.
By Proposition \ref{proposition:vplus1}
\begin{eqnarray*}
\uu{1}(\tilde{b}_{(m)}) &=& u_-^m u_+^n u_-^n u_+^n u_-^n u_+^m,\\
\mbox{Red}_1\left ( \tilde{b}_{(m)} \right )&=& u_-^m u_+^m,
\end{eqnarray*}
we have
\[
\vp{1}(\tilde{b}_{(m)}) = m.
\]
By Proposition \ref{proposition:vplus1} 
\begin{eqnarray*}
\uu{2}(\tilde{b}_{(m)})&=& u_0 u_-^{3n} u_0 u_+^{2n} u_0,\\
\mbox{Red}_2\left (\tilde{b}_{(m)} \right )&=& u_-^{3n} u_+^{3n},
\end{eqnarray*}
we have
\[
\vp{2}(\tilde{b}_{(m)}) = 3n.
\]
Therefore,
\[
\langle c, 2\vp{1}(\tilde{b}_{(m)})\Lambda_1+\vp{2}(\tilde{b}_{(m)})\Lambda_2 
\rangle=2m+3n.
\]
In similar way, we can prove other cases.

Let $B^0_l\cong\{b \mid b\in B(l\Lambda_1),\;\langle c,\vp{1}(b)\Lambda_1+
\vp{2}(b) \rangle=l\}$, 
$(P_{cl}^+)_l^0\cong\{\lambda\in \Z\Lambda_1+\Z\Lambda_2 \mid 
\langle c,\lambda \rangle =l\}$.
Then we see following proposition easily.
\begin{proposition}
The maps $\ve{},\;\vp{}\,:\, B^0_l \rightarrow (P_{cl}^+)_l^0$ are bijective.  
\end{proposition}

\subsubsection{Existence of minimal elements on $B^{G_2}(l\Lambda_1)$}
We consider $\fa{1}{i}\fa{0}{i}\hwvl{i}{i}{i}$.
By Definition \ref{definition:itinokettei}, we have 
\[
\fa{1}{i}\fa{0}{i}\hwvl{i}{i}{i} = 
\ea{2}{l-2i}\fa{1}{i}\fa{0}{i}\hwvl{i}{(l-i)}{i}.
\]
For $b\in B^i_{k\Lambda_1+j\Lambda_0}$, we define 
$y=y(b)=\left [ \frac{l-i-j}{3} \right ]$.\\
By Proposition \ref{fact:f0pR} and Lemma \ref{lemma:f1qf0p},
\[
\begin{array}{l}
\Phi\left ( \fa{1}{i}\fa{0}{i}\hwvl{i}{(l-i)}{i} \right ) = \vspace{2mm}\\
\makebox[20mm]{}
\left \{
\begin{array}{ll}
\BB{\bb{1}{i}\bb{6}{y}\bB{\CC}{y}\bb{-2}{y}\bb{-1}{i}}
&(l-2i \equiv 0 \pmod{3}), \vspace{2mm}\\
\BB{\bb{1}{i}\bb{-5}{}\bb{-1}{i}}
&(l-2i \equiv 1 \pmod{3}, y=0), \vspace{2mm}\\
\BB{\bb{1}{i}\bb{6}{y}\bb{6}{}\bB{\CC}{y-1}\bb{-4}{}\bb{-2}{y}
\bb{-1}{i}}&(l-2i \equiv 1 \pmod{3},y>0), \vspace{2mm}\\
\BB{\bb{1}{i}\bb{6}{y}\bb{6}{}\bB{\CC}{y}\bb{-3}{}\bb{-2}{y}\bb{-1}{i}}
&(l-2i \equiv 2 \pmod{3}), \vspace{2mm}
\end{array} \right .
\end{array}
\]
we have
\[
\begin{array}{l}
\Phi\left ( \ea{2}{l-2i}\fa{1}{i}\fa{0}{i}\hwvl{i}{(l-i)}{i} \right ) = 
\vspace{2mm}\\
\makebox[20mm]{}
\left \{
\begin{array}{ll}
\BB{\bb{1}{i}\bb{2}{y}\bB{\CC}{y}\bb{-2}{y}\bb{-1}{i}}
&(l-2i \equiv 0 \pmod{3}), \vspace{2mm}\\
\BB{\bb{1}{i}\bb{7}{}\bb{-1}{i}}
&(l-2i \equiv 1 \pmod{3}, y=0), \vspace{2mm}\\
\BB{\bb{1}{i}\bb{2}{y}\bb{4}{}\bB{\CC}{y-1}\bb{-4}{}\bb{-2}{y}
\bb{-1}{i}}&(l-2i \equiv 1 \pmod{3},y>0), \vspace{2mm}\\
\BB{\bb{1}{i}\bb{2}{y}\bb{3}{}\bB{\CC}{y}\bb{-3}{}\bb{-2}{y}\bb{-1}{i}}
&(l-2i \equiv 2 \pmod{3}), \vspace{2mm}
\end{array} \right .
\end{array}
\]
Therefore we have proved existence of minimal elements.  
By (\ref{eqn:A2char1p}) we see
\[
\vp{0}\left ( \fa{1}{i}\fa{0}{i}\hwvl{i}{i}{i}\right ) = 0.
\]
Thus we have
\[
\langle c, \vp{}(\tilde{b}_{(m)}) \rangle = l,
\]
where $m\in {\bf Z}_{\geq 0}$.

\subsubsection{Selection of minimals on $B^l$}
By (\ref{eqn:indsayou}), for $n < l$, we have
\[
\varphi_0 \Bigl |_{B^{n+1}}(b) =
\varphi_0 \Bigl |_{B^n}(b) +1.
\]
Then for $b \in B^n$ such that 
$\langle c,\vp{}(b) \rangle = n$, we have $\langle c,\vp{}(b) \rangle = n+1$  
on $B^{n+1}$. 
If $l=1$, by calculation we see 
$\mbox{min}\left \{\langle c, \vp{}(b) \rangle \mid b\in B^1 \right \}= 1$.
Then we see 
\[
\min \left\{\langle c,\varphi(b) \rangle \;\left |\;b\in B^l \right .\right \}
=l,
\]
inductively.
Therefore a minimal element for level $l-1$ is also a minimal element for 
level $l$. Thus it is obvious that $\ve{}$ and $\vp{}$ are bijective.

\begin{example} Following elements are minimal $(0 \leq l \leq 7)$.
\begin{eqnarray*}
\mbox{ if } l \geq 0 && \phi,\\
\mbox{ if } l \geq 1 && \BB{\bb{7}{}},\vspace{3mm}\\ 
\mbox{ if } l \geq 2 && \BB{\bb{1}{}\bb{-1}{}} \hspace{5mm}
 \BB{\bb{3}{}\bb{-3}{}},\vspace{3mm}\\
\mbox{ if } l \geq 3 && \BB{\bb{1}{}\bb{7}{}\bb{-1}{}}\hspace{5mm}
\BB{\bb{2}{}\bb{8}{}\bb{-2}{}},\vspace{3mm}\\
\mbox{ if } l \geq 4 && \BB{\bb{1}{2}\bb{-1}{2}}\hspace{5mm} 
\BB{\bb{1}{}\bb{3}{}\bb{-3}{}\bb{-1}{}}\hspace{5mm} 
\BB{\bb{2}{}\bb{4}{}\bb{-4}{}\bb{-2}{}},\vspace{3mm}\\
\mbox{ if } l \geq 5 && \BB{\bb{1}{2}\bb{7}{}\bb{-1}{2}}\hspace{5mm} 
\BB{\bb{1}{}\bb{2}{}\bb{8}{}\bb{-2}{}\bb{-1}{}}\hspace{5mm} 
\BB{\bb{2}{}\bb{3}{}\bb{8}{}\bb{-3}{}\bb{-2}{}},\vspace{3mm}\\
\mbox{ if } l \geq 6 && \BB{\bb{1}{3}\bb{-1}{3}}\hspace{5mm} 
\BB{\bb{1}{2}\bb{3}{}\bb{-3}{}\bb{-1}{2}}\\ 
&&\BB{\bb{1}{}\bb{2}{}\bb{4}{}\bb{-4}{}\bb{-2}{}\bb{-1}{}}\hspace{5mm}
\BB{\bb{2}{2}\bb{6}{}\bb{-6}{}\bb{-2}{2}},\vspace{3mm}\\
\mbox{ if } l \geq 7 && \BB{\bb{1}{3}\bb{7}{}\bb{-1}{3}}\hspace{5mm} 
\BB{\bb{1}{2}\bb{2}{}\bb{8}{}\bb{-2}{}\bb{-1}{2}}\\ 
&&\BB{\bb{1}{}\bb{2}{}\bb{3}{}\bb{8}{}\bb{-3}{}\bb{-2}{}\bb{-1}{}}\hspace{5mm}
\BB{\bb{2}{2}\bb{4}{}\bb{8}{}\bb{-4}{}\bb{-2}{2}}.
\end{eqnarray*}
\end{example}

\subsection{Connectedness of $B^l \otimes B^l$}
\label{section:connect}
We will show connectedness of $B^l \otimes B^l$, by showing that any element of
$B^l \otimes B^l$ is connected with $\phi \otimes \phi$.
We consider decomposition of tensor products $B^{G_2}(m_1\Lambda_1)\otimes
B^{G_2}(m_2\Lambda_1)$. Each connected component has a lowest weight element.
Then, we prove the lowest elements connected with $\phi \otimes \phi$.   

For $b_1\in B^{G_2}(m_1\Lambda_1)$, $b_2\in B^{G_2}(m_2\Lambda_1)$, 
there exist a sequance $(i_1, \ldots, i_k)$ 
$(i_n \in \{1,2\}, n=1, \ldots, k)$ such that $\tf{i_1}\cdots \tf{i_2}
b_1\otimes b_2$ is a lowest weight element $b_1^{\prime} \otimes 
\BB{\bb{-1}{m_2}}\in B^{G_2}(m_1\Lambda_1)\otimes B^{G_2}(m_2\Lambda_2)$. 
We can express
\[
\tilde{f}_{i_1} \cdots \tilde{f}_{i_2} b_1 \otimes b_2 = 
b_1^{\prime} \otimes \BB{\bb{-1}{m_2}}.
\]
Then there exists $b''_1 \in B^{G_2}(m'_1\Lambda_1)$ such that 
\[
\tf{0}^{(\vp{0}(b_1^{\prime})-m_2)_+ +m_2} b_1^{\prime}\otimes 
\BB{\bb{-1}{m_2}} = b''_1 \otimes \phi. 
\]
There exist another sequence $(i_1^{\prime}, \cdots, i_{k'}^{\prime})$ 
$(i_n^{\prime} \in \{1,2\}, n=1,\ldots,k')$ such that 
$\tf{0}^{m'_1}\tf{i_1^{\prime}}\cdots\tf{i_{k'}^{\prime}}b''_1$ is lowest 
weight element $\BB{\bb{-1}{m'_1}}$.
Therefore we can express as
\[
\tf{0}^{m'_1}\tf{i_1^{\prime}}\cdots\tf{i_{k'}^{\prime}}b''_1 \otimes \phi
= \tf{0}^{m'_1} \BB{\bb{-1}{m'_1}} \otimes \phi = 
\phi \otimes \phi.
\]
Then we have 
\[
b_1 \otimes b_2 \mbox{ is connected with } \phi \otimes \phi.\vspace{5mm}
\]

By \S\ref{section:minimal}, \S\ref{section:connect} we have Theorem
\ref{theorem:perfectness}.
%

%% file: appendix.tex
\section*{Perfect crystal of level 2}
\vspace*{10mm}
\unitlength .75mm
\begin{figure}[h]
\begin{center}
\begin{picture}(154,245)
\put(130,200){
\put(0,40){\Ltwos}
\put(0,40){\RRtwos{-2}{-2}{-2}{-1}{-1}{-1}}
\put(10,30){\Cone{6}{-2}}
\put(10,30){\RRthree{-2}{6}{-1}{-1}{8}{-1}{-6}{-1}}
\put(19.5,18){\line(-1,-1){5.8}}
\put(11,8){\line(-1,-1){5.8}}
\put(3.5,8){\line(-1,-1){5.8}}
\put(3,-2){\line(-4,-3){7.8}}
\put(-4.5,-2){\line(-4,-3){8}}
\put(-10,-12){\line(-4,-3){7.8}}
\put(20,20){\Cone{6}{6}}
\put(10,10){\Ctwo{6}{6}{8}}
\put(2,0){\Ctwo{8}{6}{-6}}
\put(-6.5,0){\makebox(0,0.8){$\phi$}}
\put(-10,-10){\Ctwo{-6}{8}{-6}}
\put(-20,-20){\Cone{-6}{-6}}
\put(-10,30){\Lthree}
\put(-19,18){\line(2,-3){3.8}}
\put(-14,8){\line(1,-1){6.3}}
\put(-8,8){\line(1,-1){5.8}}
\put(-7,-2){\line(5,-3){10.4}}
\put(-1,-2){\line(5,-3){9.8}}
\put(7,-12){\line(5,-3){9.8}}
\put(20,0){\Cone{1}{6}}
\put(20,0){\RRthree{1}{1}{8}{2}{1}{-6}{1}{6}}
\put(-20,0){\Lthree}
\put(20,-20){\RRtwos{1}{1}{1}{2}{2}{2}}
\put(-20,-20){\Ltwos}
}
\put(60,200){%
\put(0,40){\Ltwos}
\put(0,40){\RRone{-3}{-2}{-3}{-1}}
\put(10,30){\RRtwo{4}{-2}{-3}{4}{-1}{5}{-1}}
\put(20,20){\Cone{4}{6}}
\put(20,20){\RRthree{4}{4}{8}{5}{4}{-6}{5}{-6}}
\put(20,0){\RRtwos{1}{4}{1}{5}{2}{5}}
\put(-10,30){\Lthree}
\put(-10,10){\Ltwo}
\put(-10,-10){\Lone}
}
\put(5,200){%
\put(0,40){\Ltwos}
\put(0,40){\Cone{-4}{-2}}
\put(10,30){\RRone{3}{-2}{3}{-1}}
\put(20,20){\RRtwos{3}{6}{3}{8}{3}{-6}}
\put(0,20){\Lone}
}
\put(130,125){%
\put(0,20){\Lone}
\put(0,20){\RRtwos{-5}{-2}{-5}{-1}{-4}{-1}}
\put(10,10){\Cone{6}{-5}}
\put(10,10){\RRthree{-5}{6}{-4}{-4}{8}{-4}{-6}{-4}}
\put(10,-10){\RRtwo{1}{-5}{3}{1}{-4}{2}{-4}}
\put(10,-30){\RRone{1}{3}{2}{3}}
\put(-10,10){\Ltwo}
\put(-20,0){\Lthree}
\put(-20,-20){\Ltwos}
}
\put(60,125){%
\put(0,20){\Lone}
\put(0,20){\RRone{7}{-2}{7}{-1}}
\put(10,10){\RRtwo{4}{-5}{7}{4}{-4}{5}{-4}}
\put(10,-10){\RRone{1}{7}{2}{7}}
\put(-10,10){\Ltwo}
\put(-10,-10){\Lone}
}
\put(5,125){%
\put(0,20){\Lone}
\put(0,20){\Cone{5}{-2}}
\put(10,10){\RRone{3}{-5}{3}{-4}}
}
\put(130,70){%
\put(0,0){\RRtwos{6}{-3}{8}{-3}{-6}{-3}}
\put(0,-20){\RRone{1}{-3}{2}{-3}}
\put(0,-40){\Cone{2}{4}}
\put(-10,-10){\Lone}
\put(-20,-20){\Ltwos}
}
\put(60,70){%
\put(0,0){\RRone{4}{-3}{5}{-3}}
\put(0,-20){\Cone{2}{-5}}
\put(-10,-10){\Lone}
}
\put(5,70){\Cone{3}{-3}}
\put(12,22){%
\put(0,20){\Lone}
\put(0,20){\RRone{8}{-2}{-6}{-2}}
\put(9,8){\line(-2,-3){4}}
\put(2,-2){\line(-5,-3){9.8}}
\put(10,10){\Cone{1}{-2}}
\put(4.5,0){\Cone{1}{-1}}
\put(-4.5,0){\Cone{2}{-2}}
\put(-10,-10){\Cone{2}{-1}}
\put(10,-10){\RRone{2}{6}{2}{8}}
\put(-9,8){\line(2,-3){4}}
\put(-2,-2){\line(5,-3){9.8}}
\put(-10,-10){\Lone}
}
\put(164,15){\line(-1,0){124}}
\put(164,15){\line(0,1){235}}
\put(-10,250){\line(0,-1){195}}
\put(-10,250){\line(1,0){174}}
\put(-15,55){\line(1,0){55}}
\put(40,55){\line(0,-1){65}}
\put(-15,-10){\line(1,0){55}}
\put(-15,-10){\line(0,1){65}}
\put(126,240){\vector(-1,0){62}}
\put(56,240){\vector(-1,0){47}}
\put(97.3,130){\arc{195.8}{4.37}{5.06}}
\put(34,140){\arc{173}{4.42}{5.03}}
\put(126,145){\vector(-1,0){62}}
\put(56,145){\vector(-1,0){47}}
\put(97,35){\arc{196}{4.37}{5.06}}
\put(150.5,163.6){\arc{196}{2.764}{3.49}}
\put(140,153.6){\arc{196}{2.764}{3.49}}
\put(200,143.6){\arc{196}{2.764}{3.49}}
\put(210.5,153.6){\arc{196}{2.764}{3.49}}
\put(95,163.6){\arc{196}{2.764}{3.49}}
\put(3.3,198){\vector(1,3){0}}
\put(118.8,188){\vector(1,3){0}}
\put(108.3,178){\vector(1,3){0}}
\put(48.3,188){\vector(1,3){0}}
\put(58.8,198){\vector(1,3){0}}
\put(119,113){\vector(1,4){0}}
\put(64.5,222.5){\vector(-4,-1){0}}
\put(64.5,127.5){\vector(-4,-1){0}}
\put(8.5,222.4){\vector(-3,-1){0}}
\put(126,70){\vector(-1,0){62}}
\put(56,70){\vector(-1,0){47}}
\put(5,72){\vector(0,1){51}}
\put(50,62){\vector(0,1){51}}
\put(110,52){\vector(0,1){51}}
\put(170,175){\arc{334}{2.4}{2.949}}
\put(200,83.3){\arc{173}{2.895}{3.49}}
\put(3,238){\vector(0,-1){4.1}}
\put(3,231){\makebox(0,0){\marumoji{1}}}
\put(3,218){\vector(0,-1){4.1}}
\put(3,211){\makebox(0,0){\marumoji{2}}}
\put(0,40){\makebox(0,0){\marumoji{1}}}
\put(0,38.5){\vector(0,-1){4.2}}
\put(0,20){\makebox(0,0){\marumoji{2}}}
\put(0,18.5){\vector(0,-1){4.2}}
\put(16,42){\vector(1,0){3.9}}
\put(22,41.5){\makebox(0,0){\marumoji{3}}}
\put(6,32){\vector(1,0){3.9}}
\put(12,31.5){\makebox(0,0){\marumoji{4}}}
\put(105,59.5){\makebox(0,0){\marumoji{3}}}
\put(107,60){\vector(1,0){9}}
\put(95,49.5){\makebox(0,0){\marumoji{4}}}
\put(97,50){\vector(1,0){9}}
\put(47.4,62){\vector(1,-1){0}}
\put(156,15){\makebox(5,9)[r]{${\cal A}_0$}}
\put(32,-10){\makebox(5,9)[r]{${\cal A}_1$}}
\put(70,15){\makebox(0,9)[r]{$\swarrow \tilde{f}_1$\hspace{3mm}$\searrow 
\tilde{f}_0$}}
\put(110,15){\makebox(0,9)[r]{other arrows: $E_\CA$ }}
\put(20,24){\makebox(3,3){\small$*$}}
\put(8,72){\makebox(3,3){\small$*$}}
\put(55,127){\makebox(3,3){\small$*$}}
\put(123,201.5){\makebox(3,3){\small$*$}}
\put(140,15){\makebox(0,9)[r]{$*$ is minimal.}}
\end{picture}
\caption{Perfect crystal of level 2}
\label{picture:level2}
\end{center}
\end{figure}